\documentclass{amsart}

\usepackage[english]{babel}
\usepackage[all]{xy}
\usepackage{latexsym}
\usepackage{eufrak}
\usepackage{amssymb}
\usepackage{fancyhdr}

\newtheorem{thm}{\sc Theorem}[section]
\newtheorem{lem}[thm]{\sc Lemma}
\newtheorem{pro}[thm]{\sc Proposition}
\newtheorem{cor}[thm]{\sc Corollary}

\theoremstyle{definition}

\newtheorem*{rem}{\sc Remark}
\newtheorem*{dei}{\sc Definition}
\newtheorem*{deo}{\sc Proof}
\newtheorem*{ex}{\sc Examples}
\newtheorem*{ex1}{\sc Example}
\newtheorem*{cex1}{\sc Counterexample}
\newtheorem*{cex}{\sc Counterexamples}

\newenvironment{proo}{\begin{trivlist} \item{\sc {Proof.}}}
  {\hfill $\square$ \end{trivlist}}

\def\blfootnote{\xdef\@thefnmark{}\@footnotetext}

\long\def\symbolfootnote[#1]#2{\begingroup%
\def\thefootnote{\fnsymbol{footnote}}\footnote[#1]{#2}\endgroup}

\newcommand{\cone}{\textrm{cone}}
\newcommand{\B}{\mathcal{B}}
\newcommand{\oC}{\bar{\mathcal{C}}}

\newcommand{\oPo}{\bar{\Po}}
\newcommand{\Bi}{\mathcal{B}i}
\newcommand{\BLi}{\mathcal{B}i\mathcal{L}ie}
\newcommand{\IBi}{\varepsilon\mathcal{B}i}
\newcommand{\G}{\mathcal{G}}
\newcommand{\N}{\mathcal{N}}

\newcommand{\kj}{{\bar{k},\, \bar{\jmath}}}
\newcommand{\oi}{{\bar{\imath}}}
\newcommand{\oj}{{\bar{\jmath}}}
\newcommand{\ok}{{\bar{k}}}
\newcommand{\ol}{{\bar{l}}}
\newcommand{\coker}{\mathop{\mathrm{coker}}}

\newcommand{\FS}{\mathcal{FS}}
\newcommand{\Sy}{\mathbb{S}}

\newcommand{\C}{\mathcal{C}om}
\newcommand{\A}{\mathcal{A}s}

\newcommand{\Po}{\mathcal{P}}
\newcommand{\F}{\mathcal{F}}
\newcommand{\ac}{\scriptstyle \textrm{!`}}
\newcommand{\K}{\mathcal{K}}

\newcommand{\Qo}{\mathcal{Q}}
\newcommand{\Li}{\mathcal{L}ie}
\newcommand{\Sc}{\mathbb{S}^c}

\newcommand{\Co}{\mathcal{C}}
\newcommand{\cqfd}{\ \hfill \square}
\newcommand{\NN}{\mathbb{N}}
\newcommand{\Ss}{\mathcal{S}}
\newcommand{\Ao}{\mathcal{A}}
\newcommand{\Bo}{\mathcal{B}}

\title{\bf A Koszul duality for props}
\author{Bruno Vallette}

\begin{document}

\begin{abstract}
The notion of prop models the operations with multiple inputs and
multiple outputs, acting on some algebraic structures like the
bialgebras or the Lie bialgebras. In this paper, we generalize the
Koszul duality theory of associative algebras and operads to
props.
\end{abstract}

\maketitle

\symbolfootnote[0]{2000 \emph{Mathematics Subject Classification.}
18D50 (16W30, 17B26, 55P48). }

\section*{Introduction}

The Koszul duality is a theory developed for the first time in 1970 by S. Priddy for
associative algebras in \cite{Priddy}. To every quadratic algebra $A$, it associates a dual
coalgebra $A^{\ac}$ and a chain complex called Koszul complex. When this
complex is acyclic, we say that $A$ is a Koszul algebra. In this case, the algebra
$A$ and its representations have many properties (\emph{cf.} A. Beilinson, V. Ginzburg and
W. Soergel \cite{BGS}).\\

In 1994, this theory was generalized to algebraic operads by V.
Ginzburg and M.M. Kapranov (\emph{cf.} \cite{GK}). An operad is an
algebraic object that models the operations with $n$ inputs (and
one output) $A^{\otimes n} \to A$ acting on a type of algebras.
For instance, there exists operads $\A$, $\C$ and $\Li$ coding
associative, commutative and Lie algebras. The Koszul duality
theory for operads has many {applications:} construction of a
``small'' chain complex to compute the homology groups of an
algebra, minimal model
of an operad, notion of algebra up to homotopy.\\

The discovery of quantum groups (\emph{cf.} V. Drinfeld
\cite{Drinfeld, D2}) has popularized algebraic structures with
products and coproducts, that's-to-say operations with multiple
inputs and multiple outputs $A^{\otimes n} \to A^{\otimes m}$ . It
is the case of bialgebras, Hopf algebras and Lie bialgebras, for
instance. The framework of operads is too narrow to treat such
structures. To model the operations with multiple inputs and
outputs, one has to use a more general algebraic object : the
props. Following J.-P. Serre in
\cite{Serre}, we call an ``algebra'' over a prop $\Po$, a $\Po$-gebra. \\

It is natural to try to generalize Koszul duality theory to props.
Few works has been done in that direction by M. Markl and A.A.
Voronov in \cite{MV} and W.L. Gan in \cite{Gan}. Actually, M.
Markl and A.A. Voronov proved Koszul duality theory for what M.
Kontsevich calls $\frac{1}{2}$-PROPs and W.L. Gan proved it for
dioperads. One has to bear in mind that an associative algebra is
an operad, an operad is a $\frac{1}{2}$-PROP, a $\frac{1}{2}$-PROP
is a dioperad and a dioperad induces a prop.

$$\textrm{Associative algebras} \hookrightarrow   \textrm{Operads} \hookrightarrow \frac{1}{2}-\textrm{PROPs}
\hookrightarrow \textrm{Dioperads} \leadsto \textrm{Props}.$$

In this article, we introduce the notion of properad, which is the
connected part of a prop and contains all the information to code
algebraic structure like bialgebras, Lie bialgebras, infinitesimal
Hopf algebras (\emph{cf.} M. Aguiar \cite{Ag1, Ag2, Ag3}). We
prove the Koszul duality theory for properads which gives Koszul
duality for quadratic props defined by connected relations (it is
the case of the three previous
examples).\\

The first difficulty is to generalize the bar and cobar constructions to properads and props. The differentials
of the bar and cobar constructions for associative algebras, operads, $\frac{1}{2}$-PROPs and dioperads are
given by the notions of `edge contraction' and `vertex expansion' introduced by M. Kontsevich
in the context of graph homology (\emph{cf.} \cite{Ko}). We define the differentials of the bar and cobar
constructions of properads using conceptual properties. This generalizes the previous constructions. \\

The preceding proofs of Koszul duality theories are based on
combinatorial properties of trees or graphs of genus $0$. Such
proofs can not be used in the context of props since we have to
work with any kind of graphs. To solve this problem, we introduce
an extra graduation induced by analytic functors and generalize
the comparison lemmas
of B. Fresse \cite{Fresse} to properads. \\

To any quadratic properad $\Po$, we associated a Koszul dual coproperad $\Po^{\ac}$ and a Koszul complex.
The main theorem of this paper is a criterion which claims that the Koszul complex is acyclic if and
only if the cobar construction of the Koszul dual $\Po^{\ac}$ is a resolution of $\Po$. In this case, we say
that the properad is Koszul and this resolution is the minimal model of $\Po$. \\

As a corollary, we get the deformation theories for gebras over
Koszul properads. For instance, we prove that the properad of Lie
bialgebras and the properad of infinitesimal Hopf algebras are
Koszul, which gives the notions of Lie bialgebra up to homotopy
and infinitesimal Hopf algebra up to homotopy. Structures of
gebras up to homotopy appear in the level of on some chain
complexes. For instance, the Hochschild cohomology of an
associative algebra with coefficients into itself has a structure
of Gerstenhaber algebra up to homotopy (\emph{cf.} J.E. McClure
and J.H. Smith \cite{McCS} see also C. Berger and B. Fresse
\cite{BF}). This conjecture is know as Deligne's conjecture. The
same kind of conjectures exist for gebras with coproducts up to
homotopy. One of them was formulated by M. Chas and D. Sullivan in
the context of String Topology (\emph{cf.} \cite{CS}). This
conjecture says that the structure of involutive Lie bialgebra (a
particular type of Lie bialgebras) on some homology groups can be
lifted to a structure of involutive Lie bialgebras up to homotopy.
\\

Working with complexes of graphes of genus greater than $0$ is not
an easy task. This Koszul duality for props provides new methods
for studying the homology of these chain complexes. One can
interpret the bar and cobar constructions of Koszul properads in
terms of graph complexes. Therefore, it is possible to compute
their (co)homology in Kontsevich's sense. For instance, the case
of the properad of Lie bialgebras leads to the computation of the
cohomology of classical graphs and the case of the properad of
infinitesimal Hopf algebras leads to the computation of the
cohomology of ribbon graphs (\emph{cf.} M. Markl and A.A. Voronov
\cite{MV}). In these cases, the homology groups are given by the
associated Koszul dual coproperad. So we get a complete
description of them as $\Sy$-bimodules. In a different context, S.
Merkulov gave recently another proof of the famous formality
theorem of M. Kontsevich using the properties of chain complexes of Koszul props in \cite{Merkulov}. \\

In the last section of this paper, we generalize the Poincar\'e
series of associative algebras and operads to properads. When a
properad $\Po$ is Koszul, we prove a functional equation between
the Poincar\'e series of $\Po$ and the Poincar\'e series of its
Koszul dual $\Po^{\ac}$.

\section*{Conventions}

Throughout the text, we will use the following conventions.\\

We work over a field $k$ of characteristic $0$.

\subsection{$n$-tuples}
 To simplify the notations, we often write $\oi$ an $n$-tuple
 $(i_1,\, \ldots ,\, i_n)$.
 The $n$-tuples considered here are $n$-tuples of positive integers.
 We denote $|\oi|$ the sum $i_1+\cdots +i_n$.
 When there is no ambiguity about the number
 of terms, the $n$-tuple $(1,\, \ldots ,\, 1)$ is denoted $\bar{1}$.

We use the notation $\oi$ to represent the ``products''
of elements indexed by the $n$-tuple
$(i_1,\, \ldots ,\, i_n)$. For instance, in the case of
$k$-modules, $V_\oi$ corresponds to
$V_{i_1}\otimes_k \cdots \otimes_k V_{i_n}$.

\subsection{Partitions}
A \emph{partition of an integer $n$} is an ordered increasing
sequence $(i_1,\, \ldots ,\, i_k)$ of
positive integers such that $i_1+\cdots +i_k=n$.

We denote $[n]$ the set $\{1,\, \ldots ,\, n\}$. A
\emph{partition of the set $[n]$} is a set
$\{I_1,\, \ldots ,\, I_k\}$ of disjoint subsets of $[n]$ such
that $I_1\cup \cdots \cup I_k = [n]$.

\subsection{Symmetric group and block permutations}
We denote $\Sy_n$ the symmetric groups of permutations of $[n]$.
Every permutation $\sigma$  of $\Sy_n$
is written with the $n$-tuple $(\sigma(1),\,\ldots ,\, \sigma(n))$.
We denote $\Sy_\oi$ the subgroup $\Sy_{i_1}\times \cdots
\times \Sy_{i_n}$ of $\Sy_{|\oi|}$. From every permutation $\tau$ of
$\Sy_n$ and every $n$-tuple
$\oi=(i_1,\, \ldots,\, i_n)$, one associates a permutation $\tau_\oi$
of $\Sy_{|\oi|}$,
called a \emph{block permutation of type $\oi$}, defined by
\begin{eqnarray*}
\tau_\oi=\tau_{i_1,\ldots,\, i_n} &=&
(i_1+\cdots+i_{\tau^{-1}(1)-1}+1,\, \ldots ,\,
i_1+\cdots+i_{\tau^{-1}(1)},\,
\ldots ,\, \\
&& i_1+\cdots+i_{\tau^{-1}(n)-1}+1,\, \ldots ,\,
i_1+\cdots+i_{\tau^{-1}(n)}).
\end{eqnarray*}

\subsection{Gebra}

Following the article of J.-P. Serre \cite{Serre}, we call
\emph{gebra} any algebraic structure like
algebras, coalgebras, bialgebras, etc ...

\section{Composition products}

We describe here an algebraic framework which is used to represent
the operations with multiple inputs and multiple outputs acting on
algebraic structures. We introduce composition products which
correspond to the compositions of such operations.

\subsection{$\Sy$-bimodules}

\begin{dei}[$\Sy$-bimodule]
\index{$\Sy$-bimodule} An \emph{$\Sy$-bimodule} is a collection
$\left( \Po(m,\, n)\right)_{m,\, n \in \mathbb{N}}$, of
$k$-modules $\Po(m,\, n)$ endowed with an action of the symmetric
group $\Sy_m$ on the left and an action of the symmetric group
$\Sy_n$ on the right such that these two actions are compatible.
\end{dei}

A morphism between two $\Sy$-bimodules $\Po$, $\Qo$ is a
collection of morphisms $f_{m,\, n} \, :\, \Po(m,\, n) \to
\Qo(m,\, n)$ which commute with the action of $\mathbb{S}_m$ on
the left and with the action of $\mathbb{S}_n$ on the right.

The $\Sy$-bimodules and their morphisms form a category denoted
$\Sy$-biMod. \\

The $\Sy$-bimodules are used to code the operations acting on
different types of gebras. The module $\Po(m,\,n)$ represents the
operations with $n$ inputs and $m$ outputs.
$$\Po(m,\, n)\otimes_{\Sy_n} A^{\otimes n} \to A^{\otimes m}.$$

\begin{dei}[Reduced $\Sy$-bimodule]
\label{reducedSbiMod} A \emph{reduced $\Sy$-bimodule} is an
$\Sy$-bimodule $\big( \Po(m,$ $n)\big)_{m,\, n \in
\mathbb{N}}$ such that $\Po(m,\,0)=\Po(0,\, n)=0$, for every $m$
and $n$ in $\NN$.
\end{dei}

\subsection{Composition product of $\Sy$-bimodules}

We define a product in the category of $\Sy$-bimodules in order to
represent to composition of operations.

\begin{dei}[Directed graphs]
A \emph{directed graph} is a non-planar graph where the
orientations of the edges are given by a global flow (from the top
to the bottom, for instance). We suppose that the inputs and the
outputs of each vertex are labelled by integers $\{1,\,\ldots ,\,
n \}$. The global inputs and outputs of each graphs are also
supposed to be indexed by integers. We denote the set of such
graphs by $\G$.
\end{dei}

Every graph studied in this article will be directed.

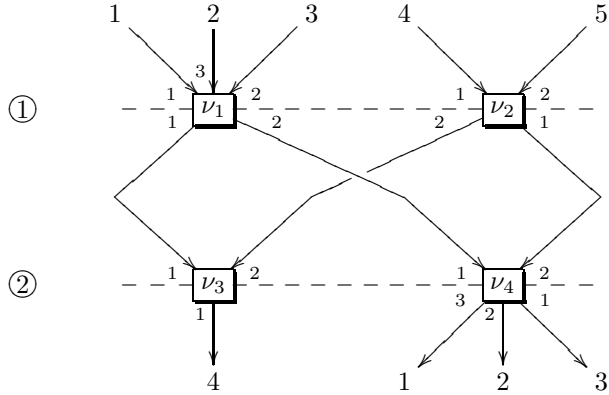
\begin{figure}[h]
$$ \xymatrix{
 & \ar[dr]_(0.7){1} 1&2 \ar[d]_(0.6){3} &3\ar[dl]^(0.7){2} &4\ar[dr]_(0.7){1} & &5\ar[dl]^(0.7){2} \\
*+[o][F-]{1} & \ar@{--}[r]& *+[F-,]{\nu_1} \ar@{-}[dl]_(0.3){1}
\ar@{-}[drr]^(0.3){2} \ar@{--}[rrr]& & & *+[F-,]{\nu_2}
\ar@{-}[dr]^(0.3){1}   \ar@{-}[dll]_(0.3){2} |(0.75) \hole  \ar@{--}[r]& \\
 & *=0{} \ar[dr]_(0.7){1} & &*=0{}\ar[dl]^(0.7){2} &*=0{}\ar[dr]_(0.7){1} & &*=0{} \ar[dl]^(0.7){2}\\
*+[o][F-]{2}& \ar@{--}[r] & *+[F-,]{\nu_3} \ar[d]_(0.3){1}
\ar@{--}[rrr]& & & *+[F-,]{\nu_4} \ar[dl]_(0.3){3} \ar[d]_(0.3){2}
\ar[dr]^(0.3){1}
\ar@{--}[r] & \\
& &4 & &1 &2 &3 } $$ \caption{Example of a 2-level graph.}
\label{2-levelled graph}
\end{figure}

\begin{dei}[2-level graphs]
When the vertices of a directed graph $g$ can be dispatched on two levels,
we say that $g$ is a \emph{2-level graph}. In this case, we denote
by $\mathcal{N}_i$ the set of vertices on the $i^{\rm th}$ level.
The set of such graphs is denoted by $\G^2$. (\emph{cf.}
Figure~\ref{2-levelled graph})
\end{dei}

We can now define a product in the category of $\Sy$-bimodules
based on 2-level graphs $\G^2$. We denote by $In(\nu)$ and
$Out(\nu)$ the sets of inputs and outputs of a vertex $\nu$ of a
graph.

\begin{dei}[Composition product $\boxtimes$]
Given two $\Sy$-bimodules $\Po$ and $\Qo$, we define their product
by the following formula

$$\Qo\boxtimes \Po := \left( \bigoplus_{g\in\G^2}
\bigotimes_{\nu \in \N_2} \Qo(|Out(\nu)|,\, |In(\nu)|) \otimes_k
\bigotimes_{\nu \in \N_1} \Po(|Out(\nu)|,\, |In(\nu)|)\right)
\Bigg/ \approx \ ,$$ where the equivalence relation $\approx$ is generated
by
$$
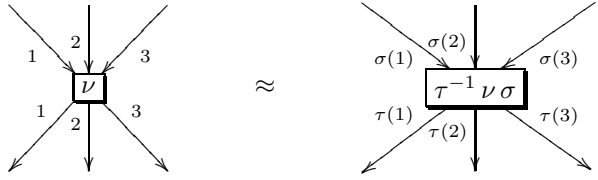
\xymatrix{\ar[dr]_(0.5){1} &\ar[d]_(0.5){2} & \ar[dl]^(0.5){3}
& & \ar[dr]_(0.5){\sigma(1)}& \ar[d]_(0.5){\sigma(2)}&
\ar[dl]^(0.5)
{\sigma(3)} \\
 & *+[F-,]{\nu} \ar[dl]_(0.4){1}\ar[d]_(0.4){2}\ar[dr]^(0.4){3}&
 &\approx & &
 *+[F-,]{\tau^{-1}\, \nu\,\sigma} \ar[dl]_(0.5){\tau(1)}\ar[d]_(0.5){\tau(2)}
 \ar[dr]^(0.5){\tau(3)}&\\
 & & & & & &\ .}$$
\end{dei}

This product is given by the directed graphs on 2 levels where the
vertices are indexed by elements of $\Qo$ and $\Po$ with respect
to the inputs and outputs.

\begin{rem}
By concision, we will often omit the equivalence relation in the
rest of the paper. It will always be implicit in the following
that a graph indexed by elements of an $\Sy$-bimodule is
considered with the equivalence relation.
\end{rem}

This product has an algebraic writing using symmetric groups.\\

Given two $n$-tuples $\oi$ and $\oj$, we use the following
conventions. The notation $\Po(\oj,\, \oi)$ denotes the product
$\Po(j_1,\, i_1)\otimes_k \cdots \otimes_k \Po(j_n,\, i_n)$ and
the notation $\Sy_{\oi}$ denotes the image of the direct product
of the groups $\Sy_{i_1}\times\cdots \times \Sy_{i_n}$ in
$\Sy_{|\oi|} $.

\begin{thm}
\label{IsoCompositionProduct}
Let $\Po$ and $\Qo$ be two $\Sy$-bimodules. Their composition
product is isomorphic to the $\Sy$-bimodule given by the formula

$$\mathcal{Q}\boxtimes \mathcal{P}(m,\, n) \cong
 \bigoplus_{N\in \mathbb{N}} \left( \bigoplus_{\ol,\, \ok,\, \oj,\, \oi} k[\mathbb{S}_m]
 \otimes_{\mathbb{S}_\ol}
\mathcal{Q}(\ol,\, \ok)\otimes_{\mathbb{S}_{\ok}} k[\Sy_{N}]
\otimes_{\mathbb{S}_\oj} \mathcal{P}(\oj,\, \oi)
\otimes_{\mathbb{S}_\oi} k[\mathbb{S}_n] \right)\Bigg/\sim  ,$$
where the direct sum runs over the $b$-tuples $\ol$, $\ok$ and the
$a$-tuples $\oj$, $\oi$ such that  $|\ol|=m$, $|\ok|=|\oj|=N$,
$|\oi|=n$ and where the equivalence relation $\sim$ is defined by
\begin{eqnarray*}
&&\theta \otimes q_1\otimes \cdots \otimes q_b \otimes \sigma
\otimes p_1 \otimes \cdots \otimes p_a \otimes\omega
  \sim \\
&&\theta \,\tau^{-1}_\ol \otimes q_{\tau^{-1}(1)}\otimes \cdots
\otimes q_{\tau^{-1}(b)} \otimes
 \tau_\ok\,
 \sigma \, \nu_\oj \otimes p_{\nu(1)} \otimes \cdots \otimes p_{\nu(a)}
 \otimes \nu^{-1}_{\oi} \,
 \omega,
\end{eqnarray*}
for $\theta \in \Sy_m$, $\omega \in \Sy_n$, $\sigma \in \Sy_N$ and
for $\tau \in \mathbb{S}_b$ with $\tau_\ok$ the
corresponding permutation by block (\emph{cf.} Conventions), $\nu
\in \mathbb{S}_a$ and $\nu_\oj$ the corresponding
permutation by block.
\end{thm}

\begin{proo}
To any element $\theta \otimes q_1\otimes \cdots \otimes q_b
\otimes \sigma \otimes p_1 \otimes \cdots \otimes p_a
\otimes\omega$ of $k[\mathbb{S}_m] \otimes \mathcal{Q}(\ol,\,
\ok)\otimes k[\Sy_N] \otimes \mathcal{P}(\oj,\, \oi) \otimes
k[\mathbb{S}_n] $, we associate a graph $\Psi(\theta \otimes
q_1\otimes \cdots \otimes q_b \otimes \sigma \otimes p_1 \otimes
\cdots \otimes p_a \otimes\omega)$ such that its vertices are
indexed by the $q_\beta$ and the $p_\alpha$, for $1\leq \beta \leq b$ and
$1\leq \alpha \leq a$. To do that, we
consider a geometric representation of the permutation $\sigma$.
We gather and index the outputs according to $\ok$ and the inputs
according to $\oj$. We index the vertices by the $q_\beta$ and the
$p_\alpha$. Then, for each $q_\beta$, we built $l_\beta$
outgoing edges whose roots are labelled by $1,\, \ldots ,\,
l_\beta$. We do the same with the $p_\alpha$ and the inputs of the
graph. Finally, we label the outputs of the graph according to
$\theta$ and the inputs according to $\omega$. For instance, the
element $(4123)\otimes q_1\otimes q_2\otimes (1324) \otimes p_1
\otimes p_2 \otimes (12435)$ gives the graph represented on
Figure~\ref{figraphecompo}.

\begin{figure}[h]
$$ \xymatrix{
\ar[dr]_(0.7){1} 1&2 \ar[d]_(0.6){2} &4\ar[dl]^(0.7){3} &3\ar[dr]_(0.7){1} & &5\ar[dl]^(0.7){2} \\
\ar@{--}[r]& *+[F-,]{p_1} \ar@{-}[dl]_(0.3){1}
\ar@{-}[drr]^(0.3){2} \ar@{--}[rrr]& & & *+[F-,]{p_2}
\ar@{-}[dr]^(0.3){2}   \ar@{-}[dll]_(0.3){1} |(0.75) \hole\newcommand{\YYY}{\vcenter{\xymatrix@M=0pt@R=6pt@C=6pt{
\ar@{-}[dr] &  &\ar@3{-}[dl]  \\
 &\ar@3{-}[d] &  \\  & &}}}
  \ar@{--}[r]& \\
 *=0{} \ar[dr]_(0.7){1} & &*=0{}\ar[dl]^(0.7){2} &*=0{}\ar[dr]_(0.7){1} & &*=0{} \ar[dl]^(0.7){2}\\
\ar@{--}[r] & *+[F-,]{q_1} \ar[d]_(0.3){1}  \ar@{--}[rrr]& & &
*+[F-,]{q_2} \ar[dl]_(0.3){1} \ar[d]_(0.3){2} \ar[dr]^(0.3){3}
\ar@{--}[r] & \\
&4 & &1 &2 &3 } $$ \caption{Image of $(4123)\otimes
q_1\otimes q_2\otimes (1324) \otimes p_1 \otimes p_2 \otimes
(12435)$ under $\Psi$ .} \label{figraphecompo}
\end{figure}
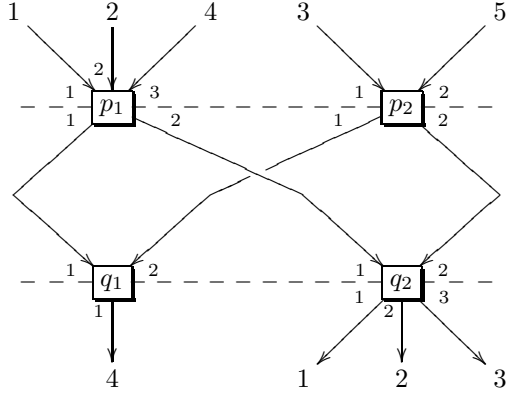

Let $\bar{\Psi}(\theta \otimes q_1\otimes \cdots \otimes q_b
\otimes \sigma \otimes p_1 \otimes \cdots \otimes p_a
\otimes\omega)$ be the equivalence class of $\Psi(\theta \otimes
q_1\otimes \cdots \otimes q_b \otimes \sigma \otimes p_1 \otimes
\cdots \otimes p_a \otimes\omega)$ for the relation $\approx$. Therefore, the
application $\bar{\Psi}$ naturally induces an application from the
quotient $k[\mathbb{S}_m]
 \otimes_{\mathbb{S}_\ol}
\mathcal{Q}(\ol,\, \ok)\otimes_{\mathbb{S}_\ok} k[\Sy_N]
\otimes_{\mathbb{S}_\oj} \mathcal{P}(\oj,\, \oi)
\otimes_{\mathbb{S}_\oi} k[\mathbb{S}_n]$. The equivalence
relation $\sim$ corresponds, via $\Psi$, to an (homeomorphic) rearrangement in
space of the graphs. Therefore, the non-planar graph created by
$\bar{\Psi}$ is invariant on the equivalent class of $\theta
\otimes q_1\otimes \cdots \otimes q_b \otimes \sigma \otimes p_1
\otimes \cdots \otimes p_a \otimes\omega$ under the relation
$\sim$. This finally defines an application $\widetilde{\Psi}$
which is an isomorphism of $\Sy$-bimodules.
\end{proo}

\begin{pro}
\label{AsCompo}
The composition product $\boxtimes$ is associative. More
precisely, let $\mathcal{R}$, $\Qo$ and $\Po$ be three
$\Sy$-bimodules. There is a natural isomorphism of $\Sy$-bimodules
$$   (\mathcal{R}\boxtimes
\mathcal{Q}) \boxtimes \mathcal{P} \cong \mathcal{R}\boxtimes
(\mathcal{Q} \boxtimes \mathcal{P}).$$
\end{pro}

\begin{proo}
Denote by $\G^3$ the set of 3-level graphs. The two
$\Sy$-bimodules $(\mathcal{R}\boxtimes \mathcal{Q}) \boxtimes
\mathcal{P}$ and $\mathcal{R}\boxtimes (\mathcal{Q} \boxtimes
\mathcal{P})$ are isomorphic to the $\Sy$-bimodule given by the
following formula
\begin{eqnarray*}
&&\Big( \bigoplus_{g\in\G^3} \bigotimes_{\nu \in \N_3}
\mathcal{R}(|Out(\nu)|,\, |In(\nu)|) \otimes \bigotimes_{\nu \in
\N_2} \Qo(|Out(\nu)|,\, |In(\nu)|) \otimes \\
&& \quad\quad \quad  \bigotimes_{\nu \in
\N_1} \Po(|Out(\nu)|,\, |In(\nu)|)\Big) \Bigg/ \approx \ .
\end{eqnarray*}
\end{proo}

\subsection{Horizontal and connected vertical composition products of $\Sy$-bimodules}

\begin{dei}[Concatenation product $\otimes$]
Let $\Po$ and $\Qo$ be two $\Sy$-bimodules. We define their
\emph{concatenation product} by the formula
$$ \Po\otimes \Qo (m,\,n) := \bigoplus_{m'+m''=m \atop n'+n''=n}
k[\Sy_{m'+m''}]\otimes_{\Sy_{m'}\times\Sy_{m''}} \Po(m',\,n')
\otimes_k \Qo(m'',\, n'') \otimes_{\Sy_{n'}\times
\Sy_{n''}}k[\Sy_{n'+n''}].$$
\end{dei}

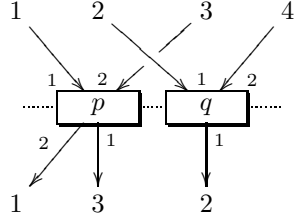
\begin{figure}[h]
$$ \xymatrix@C=10pt{ 1\ar[dr]_(0.6){1} & 2 \ar[dr]^(0.85){1} & 3
\ar[dl]_(0.85){2} | \hole &4 \ar[dl]^(0.6){2} \\
\ar@{..}[r]& *+[F-,]{\quad p\quad} \ar[d]^(0.35){1}
\ar[dl]_(0.5){2}
\ar@{..}[r]& *+[F-,]{\quad q\quad}\ar[d]^(0.35){1} \ar@{..}[r] & \\
1& 3& 2& }$$ \caption{Concatenation of two operations.}
\label{concatenation}
\end{figure}

The product $\otimes$ corresponds to the intuitive notion of
concatenation of operations (\emph{cf.}
Figure~\ref{concatenation}).

The concatenation product $\otimes$ is also called the
\emph{horizontal} product by contrast with the composition product
$\boxtimes$, which is called the \emph{vertical} product.

\begin{pro}
The concatenation product induces a symmetric monoidal category
structure on the category of $\Sy$-bimodules. The unit is given by
the $\Sy$-bimodule $k$ defined by
$$ \left\{
\begin{array}{ll}
k(0,\, 0)=k, & \\
k(m,\, n)=0 & \textrm{otherwise}.
\end{array}
\right.$$
\end{pro}

\begin{proo}
The unit relation comes from
$$ \Po\otimes k (m,\, n) = k[\Sy_m]\otimes_{\Sy_m} \Po(m,\, n) \otimes_k k
\otimes_{\Sy_n} k[\Sy_n]\cong \Po(m,\, n). $$ And the
associativity relation comes from
\begin{eqnarray*}
&&\big(\Po(m,\,n)\otimes \Qo(m',\, n')\big)\otimes R(m'',\, n'')
\cong
\Po(m,\,n)\otimes \big(\Qo(m',\, n')\otimes R(m'',\, n'')\big) \cong \\
&& k[\Sy_{m+m'+m''}]\otimes_{\Sy_m\times \Sy_{m'}\times \Sy_{m''}}
\Po(m,\,n)\otimes_k \Qo(m',\, n')\otimes_k R(m'',\, n'')\\
&& \otimes_{\Sy_n\times \Sy_{n'}\times \Sy_{n''}} k[\Sy_{n+n'+n''}].
\end{eqnarray*}
The symmetry isomorphism is given by the following formula
\begin{eqnarray*}
\Po(m,\, n) \otimes \Qo(m',\, n') &\to& (1,\,2)_{m,\, m'}
\big( \Po(m,\, n) \otimes \Qo(m',\, n') \big) (1,\, 2)_{n,\, n'}
\\&&  \cong\Qo(m',\, n') \otimes \Po(m,\, n).
\end{eqnarray*}
\end{proo}

\begin{rem}
The monoidal product $\otimes$ is bilinear. (It means that the
functors $\Po\otimes \bullet$ and $\bullet \otimes \Po$ are
additive for every $\Sy$-bimodule $\Po$).
\end{rem}

The notions of $\Sy$-bimodules represents the operations acting on
algebraic structures. The composition product $\boxtimes$ models
their compositions. For some gebras, one can restrict to
\emph{connected} compositions, based on connected graphs, without
loss of information (\emph{cf.} \ref{BiLie}).

\begin{dei}[Connected graph]
A \emph{connected graph} $g$ is a directed graph which is
connected as topological space. We denote the set of such graphs
by $\G_c$.
\end{dei}

\begin{dei}[Connected composition product $\boxtimes_c$]
Given two $\Sy$-bimodules $\Qo$ and $\Po$, we define their
\emph{connected composition product} by the following formula
$$\Qo\boxtimes_c \Po := \left( \bigoplus_{g\in\G^2_c}
\bigotimes_{\nu \in \N_2} \Qo(|Out(\nu)|,\, |In(\nu)|) \otimes_k
\bigotimes_{\nu \in \N_1} \Po(|Out(\nu)|,\, |In(\nu)|)\right)
\Bigg/ \approx \ .$$
\end{dei}

In this case, the algebraic writing is more complicated. It
involves a special class of permutations of the symmetric group
$\Sy_N$.

\begin{dei}[Connected permutations]
Let $N$ be an integer. Let $\ok=(k_1,\, \ldots,\, k_b)$ be a
$b$-tuple and $\oj=(j_1,\, \ldots ,\, j_a)$ be a $a$-tuple such
that $|\ok| = k_1+\cdots+k_b =|\oj|=j_1+\cdots+ j_a= N$.
A \emph{$(\ok,\, \oj)$-connected permutation} $\sigma$ of $\Sy_N$
is a permutation of $\Sy_N$ such that the graph of a geometric
representation of $\sigma$ is connected if one gathers the inputs
labelled by $j_1+\cdots+ j_i +1, \,
 \ldots,\, j_1+\cdots +j_i + j_{i+1}$, for $0\leq i \leq a-1$, and
 the outputs labelled by $k_1+\cdots+ k_i +1, \,
 \ldots,\, k_1+\cdots +k_i + k_{i+1}$, for $0\leq i \leq b-1$.

The set of $(\ok,\, \oj)$-connected permutations is denoted by
$\Sc_\kj \,$.
\end{dei}

\begin{ex1}
Consider the permutation $(1324)$ in $\Sy_4$ and the following
geometric representation
$$
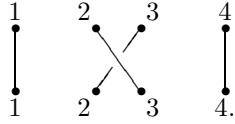
\xymatrix@C=15pt{1 \ar@{{*}-{*}}[d] &2 \ar@{{*}-{*}}[dr]&
3 \ar@{{*}-{*}}[dl] | \hole &4 \ar@{{*}-{*}}[d]\\
1 &2 & 3&4.}$$ Take  $\ok=(2,\, 2)$ and $\oj=(2,\, 2)$. If one
links the inputs $1$, $2$ and $3$, $4$ and the outputs $1$, $2$ and
$3$, $4$, it gives the following connected graph
$$
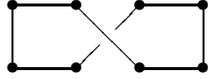
\xymatrix{ *=0{\bullet} \ar@{-}[r] \ar@{-}[d] &*=0{\bullet}\ar@{-}[dr]
&*=0{\bullet} \ar@{-}[dl] |\hole \ar@{-}[r]
&*=0{\bullet} \ar@{-}[d] \\
 *=0{\bullet} \ar@{-}[r] &*=0{\bullet} &*=0{\bullet}\ar@{-}[r] &*=0{\bullet} } $$
Therefore, the permutation $(1324)$ is $((2,\, 2),\, (2,\,
2))$-connected.
\end{ex1}

\begin{cex1}
Consider the same permutation $(1324)$ of $\Sy_4$ but let
$\ok=(1,\, 1,\, 2)$ and $\oj=(2,\, 1,\, 1)$. One obtains the
following non-connected graph
$$
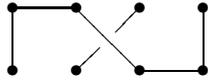
\xymatrix{ *=0{\bullet} \ar@{-}[r] \ar@{-}[d] &*=0{\bullet}\ar@{-}[dr] &*=0{\bullet} \ar@{-}[dl] |\hole
&*=0{\bullet} \ar@{-}[d] \\
 *=0{\bullet}  &*=0{\bullet} &*=0{\bullet}\ar@{-}[r] &*=0{\bullet} } $$
\end{cex1}

The next proposition will allow us the link the connected
composition product with the one defined in the theory of operads.

\begin{pro}
\label{permconnexesarbre} If $\ok=(N)$, all the permutations of
$\Sy_N$ are $((N),\, \oj)$-connected, for every $\oj$. Otherwise
stated, one has $\Sc_{(N),\, \oj}=\Sy_N$.

If $\ok$ is different from $(N)$, one has $\Sc_{\ok,\, (1,\,
\ldots,\, 1)}=\emptyset$.
\end{pro}

\begin{proo}
The proof is obvious.
\end{proo}

\begin{pro}
Let $\Po$ and $\Qo$ be two $\Sy$-bimodules. Their connected
composition product is isomorphic to the $\Sy$-bimodule given by
the formula $\mathcal{Q}\boxtimes_c \mathcal{P}(m,\, n) \cong$

$$\bigoplus_{N\in \mathbb{N}} \left( \bigoplus_{\ol,\, \ok,\, \oj,\, \oi} k[\mathbb{S}_m]
 \otimes_{\mathbb{S}_\ol}
\mathcal{Q}(\ol,\, \ok)\otimes_{\mathbb{S}_{\ok}} k[\Sc_\kj]
\otimes_{\mathbb{S}_\oj} \mathcal{P}(\oj,\, \oi)
\otimes_{\mathbb{S}_\oi} k[\mathbb{S}_n] \right) \Bigg/ \sim \ ,$$
where the direct sum runs over the $b$-tuples $\ol$, $\ok$ and the
$a$-tuples $\oj$, $\oi$ such that  $|\ol|=m$, $|\ok|=|\oj|=N$,
$|\oi|=n$.
\end{pro}

\begin{proo}
First, we have to prove that the right member is well defined.
For $\oj$ a $a$-tuple and $\ok$ a $b$-tuple such that $|\oj|=|\ok|=N$,
we denote $\oj':=(j_{\nu^{-1}(1)}$, $\ldots$,
$j_{\nu^{-1}(b)})$ and $\ok':=(k_{\tau^{-1}(1)},\, \ldots , \,
k_{\tau^{-1}(b)})$. They also verify $|\oj'|=N$ and $|\ok'|=N
$. Every $(\ok,\, \oj)$-connected permutation $\sigma$ gives by
composition on the left with a block permutation of the type
$\tau_\ok$ and by composition on the right with a bock permutation
of the type $\nu_\oj$ a $(\ok',\, \oj')$-connected permutation. The
geometric representations of $\sigma\in \Sy_N$ and
$\tau_\ok \circ \sigma \circ \nu_\oj$ are homeomorphic and one
links the same inputs and outputs. Therefore, if $\sigma\in\Sc_\kj$, one has
$\tau_\ok \, \sigma \, \nu_\oj \in \Sc_{\ok',\, \oj'}$.

The rest of the proof uses the same arguments as the proof of
Theorem~\ref{IsoCompositionProduct}.
\end{proo}

By opposition with the composition product $\boxtimes$, the connected composition
product $\boxtimes_c$ is a monoidal product in the category of $\Sy$-bimodules.
It remains to define the unit in this category. We denote
$$I:=\left\{
\begin{array}{l}
I(1,\, 1)=k,    \\
I(m,\,n)=0 \quad \textrm{otherwise}.
\end{array} \right.$$

\begin{pro}
The category $(\mathbb{S}\textrm{-biMod}, \, \boxtimes_c, \, I)$
is a monoidal category.
\end{pro}

\begin{proo}
To show the unit relation on $\Po\boxtimes_c I(m,\, n)$,
one has to study $\Sc_\kj$ in the case $\oj=(1,\, \ldots,\, 1)$. If
one does not gather all the outputs, this set is empty. Otherwise,
for $\ok=(n)$, one has $\Sc_{(n),\, (1,\, \ldots ,\, 1)}=\Sy_n$
(\emph{cf.} Proposition~\ref{permconnexesarbre}). This gives
$$ \mathcal{P}\boxtimes_c I (m,\, n) =
\mathcal{P}(m,\, n) \otimes_{\mathbb{S}_n} k\lbrack \mathbb{S}_n
\rbrack \otimes_{\mathbb{S}_1^{\times n}} k^{\otimes n} \cong
\mathcal{P}(m,\, n) \otimes k .id_n \otimes k^{\otimes n} \cong
\mathcal{P}(m,\, n).$$ (The case $I\boxtimes_c \Po\cong \Po$ is
symmetric.)

Once again, the associativity relation comes from the following
formula $\mathcal{R}\boxtimes_c
(\mathcal{Q} \boxtimes_c \mathcal{P})$ $\cong$ $(\mathcal{R}\boxtimes_c
\mathcal{Q}) \boxtimes_c \mathcal{P}$ $\cong$
\begin{eqnarray*}
&&\Big( \bigoplus_{g\in\G_3^c} \bigotimes_{\nu \in \N_3}
\mathcal{R}(|Out(\nu)|,\, |In(\nu)|) \otimes \bigotimes_{\nu \in
\N_2} \Qo(|Out(\nu)|,\, |In(\nu)|) \otimes\\
&&  \quad\quad\quad  \bigotimes_{\nu \in
\N_1} \Po(|Out(\nu)|,\, |In(\nu)|)\Big) \Bigg/\approx .
\end{eqnarray*}
\end{proo}

\begin{ex}
The monoidal categories $(k\textrm{-Mod}, \, \otimes_k, k)$ and
$(\Sy\textrm{-Mod}, \, \circ, I)$ are two full monoidal subcategories of
$(\mathbb{S}\textrm{-biMod}, \, \boxtimes_c, \, I)$.\\

The first category is the category of modules over $k$ endowed
with the classical tensor product. Every vector space $V$ can be seen as
the following $\Sy$-bimodule
$$\left\{ \begin{array}{l}
V(1,\, 1) :=V, \\
V(j,\, i) := 0 \quad \textrm{otherwise}.
\end{array} \right. $$

One has  $V\boxtimes_c W (1,\, 1)= V\otimes_k W$. Since there exists no
connected permutation associated to $a,\, b$-tuples of the form
$(1,\, \ldots ,\, 1)$ (\emph{cf.}
Proposition~\ref{permconnexesarbre}),
one has $V\boxtimes W (j,\, i)=
0$ for $(j,\, i) \ne (1,\, 1)$.
And the morphisms between two $\Sy$-bimodules only composed by an
$(\Sy_1,\,\Sy_1)$-module are exactly the morphisms of $k$-modules.\\

The category of $\Sy$-bimodules contains the category of
$\Sy$-modules related to the theory of operads (\emph{cf.}
\cite{GK, Loday3, May}). Given an $\Sy$-module $\Po$, we consider
the following $\Sy$-bimodule
$$\left\{ \begin{array}{ll}
\mathcal{P}(1,\, n) :=\mathcal{P}(n) & \textrm{for} \ n\in \mathbb{N}^*, \\
\mathcal{P}(j,\, i) := 0 & \textrm{if} \ j\ne 1.
\end{array} \right. $$

We have seen in Proposition~\ref{permconnexesarbre} that
$\Sc_{(N),\, (1,\, \ldots,\, 1)}=\Sy_N$. It gives here
$$ \mathcal{Q}\boxtimes_c \mathcal{P} (1,\, n)=
\bigoplus_{N\leq n} \left( \bigoplus_{i_1+\cdots +i_N =n}
\mathcal{Q}(1,\, N) \otimes_{\mathbb{S}_N} k[\mathbb{S}_N]
\otimes_{\mathbb{S}_1^{\times N}} \mathcal{P}(\overline{1},\, \oi)
\otimes_{\mathbb{S}_\oi} k[\mathbb{S}_n] \right) \Bigg/ \sim \ ,$$
where the equivalence relation $\sim$ is given by
$$q\otimes \sigma \nu
\otimes p_1 \otimes \cdots \otimes p_N \sim q\otimes \sigma
\otimes p_{\nu(1)}\otimes \cdots \otimes p_{\nu(N)}.$$

It is equivalent to take the coinvariants under the action of
$\Sy_N$ in the expression $\mathcal{Q}(1,\, N) \otimes_k
\mathcal{P}(1,\, i_1) \otimes_k \cdots \otimes_k \mathcal{P}(1,\,
i_N)$. We find here the monoidal product defined in the category
of $\Sy$-modules (\emph{cf.} \cite{GK, Loday3, May}), which is
well known under the following algebraic form (without the induced
representations)
$$  \mathcal{Q}\boxtimes_c \mathcal{P} (1,\, n)=
\bigoplus_{N\leq n} \left( \bigoplus_{i_1+\cdots +i_N =n}
\mathcal{Q}(1,\, N) \otimes \mathcal{P}(1,\, i_1) \otimes \cdots
\otimes \mathcal{P}(1,\, i_N) \right)_{\mathbb{S}_N} =
\mathcal{Q}\circ\mathcal{P}(n).$$

Notice that the restriction of the monoidal product $\boxtimes_c$
on $\Sy$-modules corresponds to the description of the composition
product $\circ$ defined by trees (\emph{cf.} \cite{GK, Loday3}).

When $j\ne 1$, the composition $\mathcal{Q}\boxtimes_c\mathcal{P}(j,\, i)$ represents
a concatenation of trees which is a non-connected graph. Therefore,
$\mathcal{Q}\boxtimes_c\mathcal{P}(j,\, i)=0$, for  $j\ne 1$.
\end{ex}

The default of the composition product $\boxtimes$ to be a monoidal product,
in the category of $\Sy$-bimodules, comes from the
fact that $\Po\boxtimes I$ corresponds to concatenations of elements of $\Po$ and
not only elements of $\Po$.

\begin{dei}[The $\Sy$-bimodules $T_\otimes(\Po)$ and $\Ss(\Po)$]
Since the monoidal product $\otimes$ is bilinear, the related free monoid
on an $\Sy$-bimodule $\Po$ is given by every possible concatenations of
elements of $\Po$
$$T_\otimes(\Po):=\bigoplus_{n\in \mathbb{N}} \Po^{\otimes n},$$
where $\Po^{\otimes 0}=k$ by convention.

The monoidal product is symmetric. Therefore, we can consider the
truncated symmetric algebra on $\Po$.
$$\Ss(\Po):=\bar{S}_\otimes(\Po)= \bigoplus_{n\in \mathbb{N}^*}
(\Po^{\otimes n})_{\Sy_n}.$$
\end{dei}

The functor $\Ss$ allows us to link the two composition products $\boxtimes$
and $\boxtimes_c$.

\begin{pro}
\label{lienentreproduits}
Let $\Po$ and $\Qo$ be two $\Sy$-bimodules. One has the following isomorphism
of $\Sy$-bimodules
$$\Ss(\Qo \boxtimes_c \Po) \cong \Qo \boxtimes \Po. $$
\end{pro}

\begin{proo}
The isomorphism lies on the fact that every graph of $\G^2$ is the concatenation
of connected graphs of $\G^2_c$. And the order between these connected graphs does
not matter.
\end{proo}

Remark that $\Po \boxtimes I = \Ss(\Po)$, which is different from $\Po$, in general.

\begin{dei}[Saturated $\Sy$-bimodules]
A \emph{saturated $\Sy$-bimodule} $\Po$ is an $\Sy$-bimodule such
that $\Po = \Ss(\Po)$. We denote this category sat-$\Sy$-biMod.
\end{dei}

\begin{ex1}
For every $\Sy$-bimodule $\Po$, $\Ss(\Po)$ is a saturated $\Sy$-bimodule.
\end{ex1}

Notice that $\Ss(I)(m,\, n)=0$ if $m\ne n$ and $\Ss(n,\,n)\cong k[\Sy_n]$.

\begin{pro}
The category $(\textrm{sat-}\Sy \textrm{-biMod},\, \boxtimes, \Ss(I))$ is a monoidal
category.
\end{pro}

\begin{proo}
Since $\Qo\boxtimes \Po \cong \Ss(\Qo\boxtimes \Po)$, the associativity relation
comes from Proposition~\ref{AsCompo}. If $\Po$ is a saturated $\Sy$-bimodule, one has
$\Po\boxtimes \Ss(I)=\Ss(\Po)=\Po$.
\end{proo}

We sum up these results in the diagram
$$(Vect,\, \otimes_k,\, k) \hookrightarrow (\Sy\textrm{-Mod}, \circ,\, I)
\hookrightarrow (\Sy\textrm{-biMod}, \boxtimes_c ,\, I) \xrightarrow{\Ss}
(\textrm{sat-}\Sy\textrm{-biMod}, \boxtimes,\, \Ss(I)),$$
where the first arrows are faithful functors of monoidal categories.

\begin{rem}
Since the concatenation product is symmetric and bilinear, its homological properties
are well known. In the rest of the text, we will mainly study the homological
properties of the composition products (\emph{cf.} Sections~\ref{dgframework} and
\ref{ComparisonLemmas}). Since the composition product $\boxtimes$ is obtained by concatenation
$\Ss$ from the connected composition product $\boxtimes_c$, we will only state the theorems for
the connected composition product $\boxtimes_c$ in the following of the text.
There exist analogue theorems for the composition product $\boxtimes$, which
are left to the reader.
\end{rem}

\subsection{Representation of the elements of $\Qo \boxtimes_c \Po$}

The product $\Qo \boxtimes_c\Po$ (and $\Qo \boxtimes \Po$) is isomorphic
to a quotient under the equivalence relation $\sim$
 that permutates the order of the elements of
$\Qo$ and $\Po$. For instance, to study the homological properties
of this product, we need to find natural representatives of the
classes of equivalence for the relation $\sim$.

In this section, all the $\Sy$-bimodules are reduced (\emph{cf.} \ref{reducedSbiMod}).

\begin{dei}[$\Sy_{i,\, \nu}^B$]
We denote by \emph{$\Sy_{i,\, \nu}^B$} the set of block permutations of type
$\underbrace{(i,\ldots,\,i)}_{\nu \ \textrm{times}}$ of $\Sy_{i\nu}$.
(\emph{cf.} Conventions).
\end{dei}

\begin{rem}
The set $\Sy_{i,\, \nu}^B$ is a subgroup of $\Sy_{i\nu}$.
\end{rem}

Let $\oi$ be a partition of the integer $n$. We consider the $n$-tuple
$\oi^{\ast}$ defined by $i^\ast_k:=|\{j \, / \, i_j=k\}|$. Denote by
$\Sy^B_{\oi}$ the subgroup image of $\prod_{k=1}^{n} \Sy_{k,\, i^\ast_k}^B $ in $\Sy_n$.

\begin{pro}
\label{isoquotient}
The product $\Qo \boxtimes_c \Po(m,\,n)$ is isomorphic to the $\Sy$-bimodule
defined by the following formula
$$\bigoplus_{\Theta} k[\Sy_m\big/{\Sy^B_{\ol}}]
\otimes_{\Sy_\ol} \mathcal{Q}(\ol,\,
\ok)\otimes_{\mathbb{S}_{\ok}} k[\Sy^c_\kj]
\otimes_{\mathbb{S}_{\oj}} \mathcal{P}(\oj,\, \oi)
\otimes_{\Sy_\oi} k[\Sy^B_{\oi}\big\backslash \Sy_n],$$
where the sum $\Theta$ runs over $\oi$ partition of the integer $n$,
$\ol$ partition of the integer $m$, $N$ in $\NN^*$ and $\oj$ $a$-tuple,
$\ok$ $b$-tuple such that $|\ok|=|\oj|=N$.
\end{pro}

\begin{proo}
After choosing the order of the elements of $\Po$ given by the
partition $\oi$ of the integer $n$, the remaining relations that
permutate the elements of $\Po$ involve only elements of
$\Sy^B_{\oi}$.
\end{proo}

We can do the same kind of work with the definition of $\boxtimes_c$ using
non-planar graphs. The goal here is to choose a planar representative of indexed
graphs.

\begin{dei}[Ordered partitions of $\lbrack n\rbrack$]
An \emph{ordered partition of $[n]$} is a sequence
$(\Pi_1,\ldots,\, \Pi_k)$ such that $\{\Pi_1,\ldots,\, \Pi_k\}$
is a partition of $[n]$ and such that
$\min(\Pi_1)<\min(\Pi_2)<\cdots <min(\Pi_k)$.
\end{dei}

Let $I$ be a set of $i$ elements. We denote by $\Po(j,\, i)(I)$ the module $\bigoplus_{f\, : \, [i] \to I}
\Po(j,\, i) \big/ \asymp$, where $f$ is a bijection from $[i]$ to $I$. Every
element of $\Po(j,\, i)(I)$ can be represented by a pair $(x,\, f)$. The
equivalence relation $\asymp$ is defined by $(x.\tau,\, f)\asymp (x,\, f\circ \tau)$
for $\tau$ in $\Sy_i$. Therefore,
$\Po(j,\, i)(I)$ is an $(\Sy_j,\, \Sy_I)$-module.
The module $\Po(j,\, i)(I)$ represents the vertices indexed by an element of
$\Po(j,\, i)$ with $i$ inputs labelled by elements of $I$.

\begin{pro}
\label{isopartition}
The product $\Qo \boxtimes_c \Po(m,\,n)$ is isomorphic, as $\Sy$-bimodule, to
\begin{eqnarray*}
\bigoplus_{\Theta'} && \big( (\Pi'_1)\Qo(|\Pi'_1|,\, k_1)\otimes_k
\cdots \otimes_k (\Pi'_b)\Qo(|\Pi'_b|,\, k_b)
\big) \otimes_{\Sy_{\ok}}  k[\Sy^c_\kj] \otimes_{\Sy_{\oj}}\\
&&\big( (\Po(j_1,\, |\Pi_1|)(\Pi_1)\otimes_k
\cdots \otimes_k \Po(j_a,\, |\Pi_a|)(\Pi_a)
\big),
\end{eqnarray*}
where the sum $\Theta'$ runs over pairs
$\big((\Pi_1',\ldots,\, \Pi_b'),\, (\Pi_1,\ldots,\, \Pi_a)
\big)$ of ordered partitions of $([m],\,[n])$,
$N$ in $\NN^*$ and $\oj$ $a$-tuple,
$\ok$ $b$-tuple such that $|\ok|=|\oj|=N$.
\end{pro}

As a consequence, we will often write the elements of $\Qo\boxtimes_c \Po$ like this
$$ q_1^{\Pi_1'}\otimes \cdots \otimes q_b^{\Pi_b'} \otimes
\sigma \otimes p_1^{\Pi_1} \otimes \cdots \otimes p_a^{\Pi_a} =
(q_1^{\Pi_1'} ,\, \ldots ,\,  q_b^{\Pi_b'}) \, \sigma
(p_1^{\Pi_1},\, \ldots ,\, p_a^{\Pi_a}), $$
or even without the $\Pi$'s.

\begin{rem}
The first operation $p_1$ has one input indexed by $1$ and the operation
$q_1$ has one output indexed by $1$. We will use this property
to build a contracting homotopy
for the augmented bar and cobar construction (\emph{cf.} Section~\ref{acyclicityaugm}).
\end{rem}

The same propositions can be proved for the composition product $\boxtimes$. One has to change
the set of connected permutations $\Sy^c_\kj$ by the symmetric group $\Sy_N$.

\section{Definitions of properad and prop}

In this section, we give the definitions of \emph{properad} and \emph{prop}.
These notions model the operations acting
on gebras. We give the construction of the free properad and prop. Therefore, we can define the notions of
\emph{quadratic} properad and prop and give examples.

\subsection{Definition of properad}

\begin{dei}[Properad]
A \emph{properad} is a monoid $(\Po,\, \mu,\ \eta)$ in the monoidal category
$(\Sy\textrm{-biMod},\, \boxtimes_c ,\, I)$. Equivalently, one defines a properad by
\begin{itemize}
\item An associative composition $\Po \boxtimes_c \Po \xrightarrow{\mu} \Po$,

 \item And a unit $I \xrightarrow{\eta} \Po$.
\end{itemize}
\end{dei}

\begin{ex}
$ $

\begin{itemize}
 \item The inclusions of monoidal categories
$$ (Vect,\, \otimes_k,\, k) \hookrightarrow (\Sy\textrm{-Mod}, \circ,\, I)
\hookrightarrow (\Sy\textrm{-biMod}, \boxtimes_c ,\, I) $$ show
that an associative algebra is an operad and that an operad is a
properad.

\item The fundamental example of properad is given by the
$\Sy$-bimodule $End(V)$ of multilinear applications of a vector
space $V$ (\emph{cf.} \ref{endv}).

\item The main examples of properads studied here are the free and quadratic
properads (\emph{cf.} \ref{free}, \ref{quadratic}).
\end{itemize}
\end{ex}

\begin{dei}[Augmented properad]
A properad $(\Po,\, \mu,\, \eta,\, \epsilon)$ is called \emph{augmented}
if the \emph{counit} (or augmentation map) $\epsilon \, : \, \Po \to I$ is a morphism
of properads.
\end{dei}

\begin{ex}
The free properad $\F(V)$ (\emph{cf.} \ref{free}) and the quadratic properads (\emph{cf.} \ref{quadratic})
are augmented properads.
\end{ex}

When $\Po$ is an augmented properad, we have $\Po=I \oplus \oPo$, where $\oPo$ is
the kernel of $\epsilon$ and is called the \emph{ideal of augmentation}.

\begin{dei}[Weight graded $\Sy$-bimodule]
A \emph{weight graded $\Sy$-bimodule} $M$ is a direct sum, indexed by
$\rho \in \NN$, of $\Sy$-bimodules $M=\bigoplus_{\rho \in
\mathbb{N}} M^{(\rho)}$.
\end{dei}

The morphisms of weight graded $\Sy$-bimodules are the morphisms
of $\Sy$-bimodules that preserve this decomposition. The set of
weight graded $\Sy$-bimodules with its morphisms forms a category denoted
$\textrm{gr-}\Sy\textrm{-biMod}$.\\

One can generalize the various products of $\Sy$-bimodules to the
weight graded framework. For instance, one has that $(\Qo
\boxtimes_c \Po)^{(\rho)}$ is given by the sum on 2-level
connected graphs indexed by elements of $\Qo$ and $\Po$ such that
the total sum of the weight of these elements is equal to $\rho$.

\begin{dei}[Weight graded properad]
A \emph{weight graded properad} is a monoid in the monoidal category
of weight graded $\Sy$-bimodules with the composition
product $\boxtimes_c$.
\end{dei}

A weight graded properad $\Po$ such that $\Po^{(0)}=I$ is called a
\emph{connected properad}.

\subsection{Definition of prop}

\begin{dei}[Prop]
A \emph{prop} $(\Po,\, \tilde{\mu},\, \tilde{\eta})$ is a monoid in the monoidal
category $(\textrm{sat-}\Sy\textrm{-}$ $\textrm{biMod},\, \boxtimes,\,
\Ss(I))$. In other words,
\begin{itemize}
\item The $\Sy$-bimodule $\Po$ is stable under concatenation
$\Po\otimes\Po \hookrightarrow \Po$,

\item The composition $\Po \boxtimes \Po \xrightarrow{\mu} \Po$ is
associative,

 \item The morphism $\Ss(I) \xrightarrow{\eta} \Po$ is
a unit.
\end{itemize}
\end{dei}

\begin{rem}
This definition corresponds to a $k$-linear version of the notion of PROP
(\emph{cf.} S. MacLane \cite{MacLane2} and F.W. Lawvere \cite{Lawvere}).  A \emph{prop} is
the support of the operations acting on algebraic structures.
That's-why we denote it with small letters. The name ``properad" is a contraction of ``prop" and ``operad".
\end{rem}

If one restricts the compositions of a prop to connected
graph, it gives a properad. Therefore, one can consider the
forgetful functor $U_c$ from Props to Properads
$$U_c \ : \ Props \to Properads .$$

\begin{pro}
The functor $\Ss$ induces a functor from Properads to Props. It is
the left adjoint to the forgetful functor $U_c$.
$$\xymatrix{\textrm{Props \ } \ar@<0.5ex>@^{>}[r]^(0.45){U_c}   &
\ar@<0.5ex>@^{>}[l]^(0.55){\Ss}{\ \textrm{Properads}}.} $$
\end{pro}

\begin{proo}
Let $(\Po,\, \mu ,\, \eta)$ be a properad. We have seen that the
$\Sy$-bimodule $\Ss(\Po)$ is saturated. The composition of the
prop $\Ss(\Po)$ is given by
$$\Ss(\Po)\boxtimes \Ss(\Po) \cong \Ss(\Po\boxtimes_c \Po) \xrightarrow{\Ss(\mu)} \Ss(\Po).$$
And the unit comes from
$$\Ss(I) \xrightarrow{\Ss(\eta)} \Ss(\Po).$$
The adjunction relation is left to the reader.
\end{proo}

\subsection{Definition of coproperad}

Dually, we define the notion of coproperad.

\begin{dei}[Coproperad]
A \emph{coproperad} is a comonoid $(\Co,\, \Delta,\ \varepsilon)$
in the monoidal category $(\Sy\textrm{-biMod},\, \boxtimes_c ,\,
I)$. A structure of coproperad on $\Co$ is given by
\begin{itemize}
\item A coassociative coproduct $\Co \xrightarrow{\Delta} \Co \boxtimes_c \Co$,
\item And a counit $\Co \xrightarrow{\varepsilon} I$.
\end{itemize}
\end{dei}

\begin{ex}
The main examples of coproperads studied in this article are the cofree connected coproperad
$\F^c(V)$ (\emph{cf.} \ref{connectedcofreecoproperad}) and the Koszul
dual $\Po^{\ac}$ of a properad $\Po$ (\emph{cf.} \ref{Koszuldualofproperad}).
\end{ex}

\subsection{The example of $End(V)$}
\label{endv}

\begin{dei}[$End(V)$]
We denote by $End(V)$ the set $\{ \textrm{linear applications} \ :
\ V^{\otimes n} \to V^{\otimes m}\}_{n,\, m}$. The symmetric
groups $\Sy_n$ and $\Sy_m$ act on $V^{\otimes n}$ and $V^{\otimes
m}$ by permutation of the variables. Therefore, $End(V)$ is an
$\Sy$-bimodule. The composition $End(V)\boxtimes End(V)
\xrightarrow{\chi} End(V)$ is given by the composition of linear
applications. And the inclusion $\eta \, : \, \Ss(I)
\hookrightarrow End(V)$ corresponds to linear applications
$V^{\otimes n} \to V^{\otimes n}$ obtained by permutation of the
variables.
\end{dei}

For instance, let  $p_\alpha \in Hom_{k\textrm{-Mod}}(V^{\otimes
i_\alpha},\, V^{\otimes j_\alpha})$ and $q_\beta\in
Hom_{k\textrm{-Mod}}(V^{\otimes k_\beta},\, V^{\otimes l_\beta})$,
the morphism $\chi(\theta\otimes q_1 \otimes \cdots \otimes q_b
\otimes \sigma \otimes p_1\otimes\cdots \otimes p_a \otimes
\omega)\in Hom_{k\textrm{-Mod}}(V^{\otimes n},\, V^{\otimes m})$
corresponds to the following composition
$$ \xymatrix@C=35pt{V^{\otimes n} \ar[r]^-{\omega}&V^{\otimes n}
\ar[r]^-{p_1\otimes \cdots \otimes p_a }&V^{\otimes N}
\ar[r]^-{\sigma}&V^{\otimes N}\ar[r]^-{q_1\otimes \cdots \otimes
q_b}& V^{\otimes m} \ar[r]^-{\theta}  &V^{\otimes m},}$$ where the
morphism $p_1\otimes\cdots \otimes p_a\, : \, V^{\otimes n}=
 V^{\otimes i_1} \otimes \cdots \otimes V^{\otimes i_b} \to
 V^{\otimes j_1} \otimes \cdots \otimes V^{\otimes j_b}=
 V^{\otimes N}$ is the concatenation of morphisms $p_\alpha$.\\

The $\Sy$-bimodule $End(V)$ is saturated and $(End(V),\, \chi,\, \eta)$ is a prop.
Its restriction to connected compositions $End(V)\boxtimes_c
End(V) \xrightarrow{\chi_c} End(V)$ induces a structure of
properad.

\subsection{$\Po$-module}

\label{P-module}
We recall here the basic definitions related to the notion of module
over a monoid that will be applied to properads and props throughout the text. For more
details, we refer the reader to the book of S. MacLane \cite{MacLane1}.\\

A structure of \emph{module (on the right)} $L$ over a monoid
$\Po$ in a monoidal category is given by a morphism $L \boxtimes_c
\Po \xrightarrow{\rho} \Po$ such that the following diagrams commute
$$
\begin{array}{cc}
\xymatrix{L\boxtimes_c \Po \boxtimes_c \Po \ar[r]^(0.6){L \boxtimes_c
\mu} \ar[d]^{\rho \boxtimes_c \Po} & L\boxtimes_c\Po \ar[d]^{\rho} \\
L \boxtimes_c \Po \ar[r]^{\rho} & L} &
\xymatrix{L\boxtimes_c I \ar[r]^{L\boxtimes_c \eta} \ar[rd]
& L \boxtimes_c \Po \ar[d]^{\rho}\\ & L.}
\end{array}$$
In this case, the object $L$ is called a \emph{$\Po$-module (on the right)}.\\

The free module on the right on an object $M$ is given by $L:=M\boxtimes_c \Po$.
When $\Po$ is a prop, the free $\Po$-module on a saturated $\Sy$-module $M$
is $M\boxtimes \Po$. In the rest of the text, we will call the \emph{free $\Po$-module}
on an $\Sy$-bimodule $M$, the following $\Po$-module $\Ss(M)\boxtimes \Po =M \boxtimes \Po$.

\begin{dei}[Relative composition product]
Let $(L,\, \rho)$ be a right $\Po$-module and $(R,\, \lambda)$
be a left $\Po$-module. Their \emph{relative composition product}
$L{\boxtimes_c}_{\Po} R$ is given by the cokernel of the following maps
$$\xymatrix{L \boxtimes_c \Po \boxtimes_c R \,  \ar@<0.5ex>[r]^(0.6){\rho\boxtimes_c R}
\ar@<-0.5ex>[r]_(0.6){L \boxtimes_c \lambda} & L\boxtimes_c R.}$$
\end{dei}

When $(\Po,\, \mu, \, \eta,\, \epsilon)$ is an augmented monoid,
one can define a (left) action of $\Po$ on $I$ by the
formula
$$ \lambda_\epsilon\ : \
\Po \boxtimes_c I \xrightarrow{\epsilon \boxtimes_c I} I \boxtimes_c
I \cong I.$$

\begin{dei}[Indecomposable quotient]
Let $(\Po,\, \mu, \, \eta,\, \epsilon)$ be an augmented monoid
 in a monoidal category and let
$(L,\, \rho)$ be a right $\Po$-module. The \emph{indecomposable quotient}
of $L$ is the relative composition product of $L$ and $I$, that's-to-say
the cokernel of the maps
$$\xymatrix{L \boxtimes_c \Po \boxtimes_c I \ar@<0.5ex>[r]^{\rho\boxtimes_c I}
\ar@<-0.5ex>[r]_{L \boxtimes_c \lambda_\epsilon} & L\boxtimes_c I \cong L.}$$
\end{dei}

When $(\Po,\, \mu, \, \eta,\, \epsilon)$ is an augmented monoid,
one has that the indecomposable quotient of the free module
$L=M\boxtimes_c \Po$ is isomorphic to $M$.

\subsection{$\Po$-gebra}

The notion of \emph{$\Po$-gebra} is the natural generalization to
the notion of an algebra over an operad.

\begin{dei}[$\Po$-gebra]
Let $(\Po,\, \mu,\, \eta)$ be a properad (or a prop). A structure
of \emph{$\Po$-gebra} on a $k$-module $V$ is given by a morphism
of properads (or props) : $\Po \to End(V)$.
\end{dei}

\begin{rem}
A $\Po$-gebra is an ``algebra" over a prop $\Po$. The term
``algebra" is to be taken here in the largest sense. For instance, a
$\Po$-gebra can be equipped with coproducts (operations with multiple outputs).
It is the case for
bialgebras, and Lie bialgebras.
\end{rem}

Consider the $\Sy$-module $T(V)$ defined by $T(V)(n):=V^{\otimes
n}$, where the (left) action is given by the permutation of
variables. One can see that a morphism of (left) $\Sy$-modules
$\mu_V \, : \, \Po \boxtimes V \to T(V)$ induces a unique morphism
of $\Sy$-modules $\widetilde{\mu_V} \, : \, \Po \boxtimes T(V) \to
T(V)$. For instance, the element
$$ 
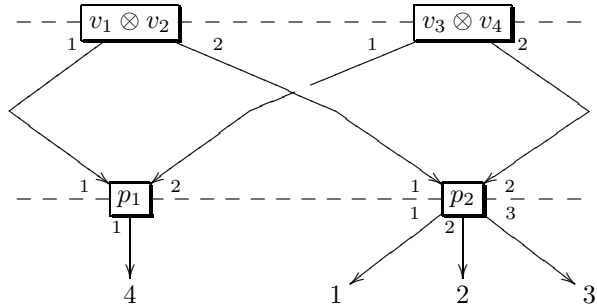
\xymatrix{\ar@{--}[r]& *+[F-,]{v_1\otimes v_2} \ar@{-}[dl]_(0.4){1}
\ar@{-}[drr]^(0.4){2} \ar@{--}[rrr]& & & *+[F-,]{v_3\otimes v_4}
\ar@{-}[dr]^(0.4){2}   \ar@{-}[dll]_(0.4){1} |(0.75) \hole  \ar@{--}[r]& \\
 *=0{} \ar[dr]_(0.7){1} & &*=0{}\ar[dl]^(0.7){2} &*=0{}\ar[dr]_(0.7){1}
 & &*=0{} \ar[dl]^(0.7){2}\\
\ar@{--}[r] & *+[F-,]{p_1} \ar[d]_(0.3){1}  \ar@{--}[rrr]& & &
*+[F-,]{p_2} \ar[dl]_(0.3){1} \ar[d]_(0.3){2} \ar[dr]^(0.3){3}
\ar@{--}[r] & \\
&4 & &1 &2 &3 } $$ can be written
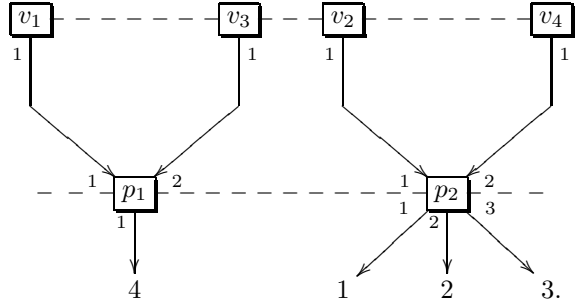
$$ \xymatrix{*+[F-,]{v_1} \ar@{-}[d]_(0.4){1}  \ar@{--}[rr] &
& *+[F-,]{v_3} \ar@{-}[d]^(0.4){1}] \ar@{--}[r] & *+[F-,]{v_ 2}
\ar@{-}[d]_(0.4){1} \ar@{--}[rr] &
& *+[F-,]{v_4} \ar@{-}[d]^(0.4){1} \\
 *=0{} \ar[dr]_(0.7){1} & &*=0{}\ar[dl]^(0.7){2} &*=0{}\ar[dr]_(0.7){1}
 & &*=0{} \ar[dl]^(0.7){2}\\
\ar@{--}[r] & *+[F-,]{p_1} \ar[d]_(0.3){1}  \ar@{--}[rrr]& & &
*+[F-,]{p_2} \ar[dl]_(0.3){1} \ar[d]_(0.3){2} \ar[dr]^(0.3){3}
\ar@{--}[r] & \\
&4 & &1 &2 &3. } $$

\begin{pro}
A $\Po$-gebra structure on a $k$-module $V$ is equivalent to a
morphism $\mu_V\, :\, \Po\boxtimes V \to T(V)$ such that the
following diagram commutes
$$\xymatrix{ \Po \boxtimes \Po \boxtimes V  \ar[r]^-{\Po
\boxtimes \mu_V} \ar[d]^-{\mu \boxtimes V} &
\Po \boxtimes T(V)\ar[d]^-{\widetilde{\mu_V}} \\
\Po \boxtimes V \ar[r]^-{\mu_V} & T(V). } $$
\end{pro}

\begin{rem}
The same arguments hold for a gebra over a properad. When $\Po$ is
an operad, we find the classical notion of an algebra over an
operad.
\end{rem}

\subsection{Free properad and free prop}
\label{free}

In this section, we complete the following diagram of forgetful and left adjunct functors
$$\xymatrix{ \textrm{Props} \ar@<0.5ex>@^{>}[d]^{U}  \ar@<0.5ex>@^{>}[r]^{U_c}
& \ar@<0.5ex>@^{>}[l]^{\Ss} \textrm{Properads} \ar@<0.5ex>@^{>}[d]^{U} \\
 \ar@<0.5ex>@^{>}[u]^{\widetilde{{\F}}} \textrm{sat-}\Sy\textrm{-biMod}
 \ar@<0.5ex>@^{>}[r]^(0.55){U} &
 \ar@<0.5ex>@^{>}[u]^{\F}  \ar@<0.5ex>@^{>}[l]^(0.45){\Ss} \Sy\textrm{-biMod},} $$
with the description of the functors $\F$ and $\widetilde{\F}$. To do this, we will
use the construction of the free monoid described in \cite{BVfreemonoid}.\\

We recall the construction of the free properad given in \cite{BVfreemonoid}
(section 5).

\begin{thm}
\label{proplibre}
The free properad on an $\Sy$-bimodule $V$ is given by the sum
on connected graphs (without level) $\G$ with the vertices indexed by elements
of $V$
$$\F(V)=\left( \bigoplus_{g\in \mathcal{G}_c} \bigotimes_{\nu \in
\mathcal{N}} V(|Out(\nu)|,\, |In(\nu)|)\right) \Bigg/ \approx \
.$$
The composition $\mu$ comes from the composition of directed graphs.
\end{thm}

\begin{rem}
If $V$ is an $\Sy$-module, the connected graphs involved here are trees.
Therefore, we find here the same construction of the free operad given in \cite{GK}.

By the same arguments, we have that the free prop $\widetilde{\F}(V)$
on a saturated $\Sy$-bimodule $V$ is given by the sum on graphs (without level) indexed by elements
of $V$. We find here the
same construction of the free prop as B. Enriquez and P.
Etingof in \cite{EE}.
\end{rem}

\begin{dei}[Free prop]
In this paper, we call \emph{free prop} on an $\Sy$-module $V$, the
image of $V$ by the following composition of left adjunct functors
$$\widetilde{\F}(V):=\Ss \circ \F(V) .$$
\end{dei}

Recall from \cite{BVfreemonoid}, the notion of \emph{split analytic functor}.

\begin{dei}[Split analytic functors]
\label{splitanalyticfunctor}
A functor $f$ is a \emph{split analytic functor} if it
is a direct sum $\bigoplus_{n \in \NN} f_{(n)}$ of homogenous
polynomial functors of degree $n$.
\end{dei}

One can decompose the set of directed (connected) graphs with the number of vertices.
We denote by $\G_{c,\, (n)}$ the set of connected graphs with $n$ vertices. Hence, one has
\begin{eqnarray*}
\F(V) &=& \left( \bigoplus_{g\in \mathcal{G}_c}
\bigotimes_{\nu \in \mathcal{N}} V(|Out(\nu)|,\, |In(\nu)|)
\right) \Bigg/ \approx\\
&=& \bigoplus_{n\in \mathbb{N}} \underbrace{
\left( \bigoplus_{g\in \mathcal{G}_{c,\, (n)}} \bigotimes_{i=1}^n
V(|Out(\nu_i)|,\, |In(\nu_i)|)\right)
\Bigg/ \approx}_{\F_{(n)}(V)}.
\end{eqnarray*}

\begin{pro}
\label{Fanalytique}
The functor $\F \, : \, V \mapsto \F(V)$
is a split analytic functor, where the part of weight $n$ is denoted
$\F_{(n)}(V)$. This graduation is stable under the composition $\mu$ of the
free properad. Therefore, the free properad is a weight graded properad.

We define the counit $\epsilon$ by the following projection $\F(V) \to I=\F_{(0)}(V)$.
With this augmentation map, the free properad is an augmented properad.
\end{pro}

\begin{lem}
\label{bigraduation} Let $f=\bigoplus_{n=0}^\infty f_{(n)}$ be a
split analytic functor in the category of $\Sy$-bimodules. Let
$M=\bigoplus_{\rho\in \mathbb{N}} M^{(\rho)}$ be a weight graded
$\Sy$-bimodule. The $\Sy$-bimodule $f(M)$ is bigraded with one
graduation induced by the analytic functor and the other one given
by the total weight.
\end{lem}

\begin{proo}
Denote $f_{(n)}=f_n\circ \Delta_n$ where $f_n$ is an $n$-linear
morphism. The first graduation can be written
$f(M)_{(n)}:=f_n(M,\ldots,\, M)$. The second one corresponds to
$$f(M)^{(\rho)}:=  \sum_{n\ge 1 \atop i_1+\cdots+i_n=\rho} f_n(M^{(i_1)},\ldots,\, M^{(i_n)}).$$
\end{proo}

\begin{cor}
Every free properad on a weight graded $\Sy$-bimodule is bigraded.
\end{cor}

\subsection{Cofree connected coproperad}
\label{connectedcofreecoproperad}

We define the structure of cofree connected coproperad on the
same $\Sy$-bimodule as the free properad. Denote
$$\F^c(V):= \F(V)=\left( \bigoplus_{g\in \mathcal{G}^c}
\bigotimes_{\nu \in \mathcal{N}} V(|Out(\nu)|,\, |In(\nu)|)\right)
\Bigg/ \approx.$$

The projection on the summand generated by the trivial graph
$\ \vcenter{\xymatrix@M=0pt@R=10pt@C=6pt{ {}\ar@{-}[d] \\
{}}} \ $ gives the counit $\varepsilon \, : \, \F^c(V) \to I$.

The coproduct $\Delta$ lies on the set of cuttings of graphs into
two parts. The image of an element $g(V)$ represented by a graph
$g$ indexed by operations of $V$ under the coproduct $\Delta$ is
given by the formula
$$\Delta(g(V)):=\sum_{\left (g_1(V),\, g_2(V) \right)} g_1(V) \boxtimes_c g_2(V),$$
where the sum is over pairs $(g_1(V),\, g_2(V))$ of family of elements such that
$\mu(g_1(V) \boxtimes_c g_2(V))=g(V)$.

\begin{pro}
For every $\Sy$-bimodule $V$, $(\F^c(V),\, \Delta,\, \varepsilon)$
is a weight graded connected coproperad. This properad is cofree
is the category of connected coproperads.
\end{pro}

\begin{deo}
The counit relation comes from
\begin{eqnarray*}
(\varepsilon \boxtimes_c id)\circ \Delta \left( g(V) \right) &=&
(\varepsilon \boxtimes_c id)\circ \Delta \left(\sum_{\left
(g_1(V),\, g_2(V) \right)} g_1(V) \boxtimes_c g_2(V)\right) \\
&=& \sum_{\left (g_1(V),\, g_2(V) \right)} \varepsilon(g_1(V))
\boxtimes_c g_2(V) \\
&=& I \boxtimes_c g(V)\\
 &=& g(V).
\end{eqnarray*}
The coassociativity relation comes from
\begin{eqnarray*}
&& (\Delta\boxtimes_c id)\circ \Delta (g(V)) =(id \boxtimes_c
\Delta)\circ \Delta (g(V)) = \\
&& \sum_{(g_1(V),\ g_2(V),\, g_3(V))} g_1(V)\boxtimes_c
g_2(V)\boxtimes_c g_3(V),
\end{eqnarray*}
where the sum is over the triples $(g_1(V),\ g_2(V),\, g_3(V))$
such that $\mu(g_1(V) \boxtimes_c g_2(V) \boxtimes_c g_3(V))=g(V)$.

Let $C$ be a connected coproperad and let $f\, : \, C \to V$ be a
morphism of $\Sy$-bimodules. The analytic decomposition of
$\F^c(V)$, according to the number of vertices of graphs, and the
fact that $C$ is connected allows to define by induction a
morphism of connected coproperads $\bar{f} \, : \, C \to \F^c(V)$.
This morphism is determined by $f$ and is the unique morphism of
connected coproperads such that the following diagram commutes
$$\xymatrix{V & \ar[l] \F^c(V) \\
& \ar[ul]^{f} \quad C.  \ar[u]_{\bar{f}}}$$ $\cqfd$
\end{deo}

\begin{rem}
This construction generalizes the construction given by T. Fox in \cite{F}
for cofree connected coalgebras.
\end{rem}

\subsection{Quadratic properads and quadratic props}
\label{quadratic}

We define the notions of quadratic properad and quadratic prop. We relate these two notions and give
examples. \\

In this paper, we are interested in the study of properads and props defined by generators and relations.
Let $V$ be an $\Sy$-bimodule and $R$ be a sub-$\Sy$-bimodule of $\F(V)$. Consider the
ideal of $\F(V)$ generated by $R$. (For the notion of ideal in a monoidal category, we refer the reader to
\cite{BVthese}). The quotient properad $\F(V)/(R)$ is generated by $V$ and the
relations $R$. \\

The Koszul duality for props deals with a subclass of properads and props
composed by \emph{quadratic properads} and \emph{quadratic props}.

\begin{dei}[Quadratic properads]
A \emph{quadratic properad} is a properad defined by generators $V$ and relations
$R$ such that $R$ is a sub-$\Sy$-bimodule of $\F_{(2)}(V)$ (the part of weight $2$
of $\F(V)$).
\end{dei}

It means that $R$ is composed by linear combinations of
connected graphs indexed by $2$ elements of $V$.

\begin{pro}
Every properad $\Po=\F(V)/(R)$ defined by generators $V$ and relations
$R$ such that $R$ is a sub-$\Sy$-bimodule of
$\F_{(n)}(V)$ is a weight graded properad.
\end{pro}

\begin{proo}
Since the ideal $(R)$ is generated by an homogenous weight graded
$\Sy$-bimodule, the quotient properad $\F(V)/(R)$ is weight graded. The elements of
weight $n$ of $\Po$ corresponds to the image of $\F_{(n)}(V)$ under the projection
$\F(V)\twoheadrightarrow \Po=\F(V)/(R)$.
\end{proo}

\begin{cor}
Every quadratic properad is weight graded.
\end{cor}

Once again, the same results hold for props.

\begin{dei}[Quadratic props]
A \emph{quadratic prop} is a prop defined by a generating $\Sy$-bimodule $V$ and relations
$R$ such that $R$ is a sub-$\Sy$-bimodule of $\widetilde{\F}_{(2)}(V)$ (the part of weight $2$
of $\widetilde{\F}(V)$).
\end{dei}

Quadratic properads and quadratic props are related by the following proposition.

\begin{pro}
\label{quadraticconnectedprop}
Let $V$ be an $\Sy$-bimodule and $R$ a sub-$\Sy$-bimodule of $\F_{(2)}(V)
\subset \widetilde{\F}_{(2)}(V)$. The quotient prop $\widetilde{\F}(V)/(R)$ is isomorphic to
$\Ss(\F(V)/(R))$, where
$\F(V)/(R)$ is the quadratic properad
generated by $V$ and $R$.
\end{pro}

\begin{proo}
Consider the morphism of props
$$\widetilde{\F}(V)=\widetilde{\F}\circ \Ss (V) =\Ss\circ \F(V) \xrightarrow{\Phi}
\Ss(\F(V)/(R)).$$
Since its kernel is the ideal generated by $R$ in $\F(V)$, $\Phi$ induces
an isomorphism $\widetilde{\F}(V)/(R) \xrightarrow{\bar{\Phi}} \Ss(\F(V)/(R))$.
\end{proo}

\begin{rem}
This proposition implies that, when the relations of a quadratic prop are defined
by linear combinations of connected graphs, it is enough to study the related quadratic
properad. In this case, the quadratic prop is obtained by concatenation
from the quadratic properad. Since the concatenation product
is well known in homology, it is enough to study the properad associated to some
algebraic structures (gebras)
to understand their homological properties (deformation, gebra up to
homotopy). (\emph{cf.} \ref{KDprops}).
\end{rem}

\begin{ex}
Once again, the first examples come from the full subcategories $k$-Mod and $\Sy$-Mod.
Quadratic algebras, like $S(V)$ and $\Lambda(V)$ (the original
 examples of J.-L. Koszul), are quadratic properads. And quadratic operads,
like $\A$, $\C$ and $\Li$ (coding associative, commutative and Lie
algebras), are quadratic properads.
\end{ex}

\begin{ex}
Let us give new examples treated by this theory. All the following properads
are quadratic.
\begin{itemize}
\item \label{BiLie}
The properad $\BLi$
coding Lie bialgebras. This properad is generated by the following $\Sy$-bimodule
$$V=\left\{ \begin{array}{l}
V(1,\, 2)=\lambda.k\otimes sgn_{\Sy_2}, \\
V(2,\, 1)=\Delta.sgn_{\Sy_2}\otimes k, \\
V(m,\, n) = 0 \quad \textrm{otherwise.}
\end{array}
\right.= \vcenter{\xymatrix@M=0pt@R=6pt@C=4pt{{\scriptstyle 1} & &
{\scriptstyle 2}
\\\ar@{-}[dr] &
&\ar@{-}[dl]  \\  &\ar@{-}[d] & \\
& & \\ & & }}\otimes sgn_{\Sy_2} \oplus
\vcenter{\xymatrix@M=0pt@R=6pt@C=4pt{  & \ar@{-}[d]&  \\
&\ar@{-}[dl] \ar@{-}[dr]& \\ & & \\{\scriptstyle 1}&
&{\scriptstyle 2} }} \otimes sgn_{\Sy_2},$$
where $sgn_{\Sy_2}$ is the signature representation of $\Sy_2$. Its
quadratic relations are given by
\begin{eqnarray*}
R&=&\left\{
\begin{array}{l}
\lambda (\lambda,\, 1) \big( (123)+(231)+(312)  \big) \\
\oplus \ \big((123)+(231)+(312) \big)(\Delta,\, 1)\Delta  \\
\oplus \ \Delta \otimes \lambda - (\lambda,\, 1)\otimes(213)\otimes(1,\, \Delta) \\
\quad -(1,\, \lambda)(\Delta,\, 1) - (1,\,\lambda
)(\Delta,\,1)-(1,\, \lambda )\otimes(132)\otimes(\Delta,1)
\end{array}
\right. \\
&=&\left\{ \begin{array}{l}
\vcenter{\xymatrix@M=0pt@R=6pt@C=4pt{{\scriptstyle 1} & &
{\scriptstyle 2}& & {\scriptstyle 3}\\
\ar@{-}[dr] &
&\ar@{-}[dl] & & \ar@{-}[dl]  \\
& \ar@{-}[dr] & &\ar@{-}[dl]  & \\
& &\ar@{-}[d] & & \\
& & \\ & & }} + \vcenter{\xymatrix@M=0pt@R=6pt@C=4pt{{\scriptstyle
2} & &
{\scriptstyle 3}& & {\scriptstyle 1}\\
\ar@{-}[dr] &
&\ar@{-}[dl] & & \ar@{-}[dl]  \\
& \ar@{-}[dr] & &\ar@{-}[dl]  & \\
& &\ar@{-}[d] & & \\
& & \\ & & }} +\vcenter{\xymatrix@M=0pt@R=6pt@C=4pt{{\scriptstyle
3} & &
{\scriptstyle 1}& & {\scriptstyle 2}\\
\ar@{-}[dr] &
&\ar@{-}[dl] & & \ar@{-}[dl]  \\
& \ar@{-}[dr] & &\ar@{-}[dl]  & \\
& &\ar@{-}[d] & & \\
& & \\ & & }} \\
\oplus \vcenter{\xymatrix@M=0pt@R=6pt@C=4pt{ & & & & \\
 & &\ar@{-}[d] & & \\
& & \ar@{-}[dl]  \ar@{-}[dr] & & \\
& \ar@{-}[dl]  \ar@{-}[dr] & &\ar@{-}[dr] & \\
& & & &\\ {\scriptstyle 1} & & {\scriptstyle 2}& &{\scriptstyle 3}
}} +
\vcenter{\xymatrix@M=0pt@R=6pt@C=4pt{ & & & & \\
 & &\ar@{-}[d] & & \\
& & \ar@{-}[dl]  \ar@{-}[dr] & & \\
& \ar@{-}[dl]  \ar@{-}[dr] & &\ar@{-}[dr] & \\
& & & &\\ {\scriptstyle 2} & & {\scriptstyle 3}& &{\scriptstyle 1}
}} + \vcenter{\xymatrix@M=0pt@R=6pt@C=4pt{ & & & & \\
 & &\ar@{-}[d] & & \\
& & \ar@{-}[dl]  \ar@{-}[dr] & & \\
& \ar@{-}[dl]  \ar@{-}[dr] & &\ar@{-}[dr] & \\
& & & &\\ {\scriptstyle 3} & & {\scriptstyle 1}& &{\scriptstyle
2}}} \\
\ \\
\oplus \vcenter{\xymatrix@M=0pt@R=6pt@C=4pt{{\scriptstyle 1} & &
{\scriptstyle 2}
\\\ar@{-}[dr] &  &\ar@{-}[dl]  \\  &\ar@{-}[d] & \\
& \ar@{-}[dl]\ar@{-}[dr]& \\ & & \\ {\scriptstyle 1} & &
{\scriptstyle 2}}} -
 \vcenter{\xymatrix@M=0pt@R=6pt@C=2pt{ & {\scriptstyle 1} & &{\scriptstyle 2}
 \\ & \ar@{-}[d] &  &\ar@{-}[d]  \\  & \ar@{-}[dl]
 \ar@{-}[dr] &  & \ar@{-}[dl]\\
\ar@{-}[d]&  & \ar@{-}[d] & \\ & & & \\ {\scriptstyle 1} & &
{\scriptstyle 2} &}} +
 \vcenter{\xymatrix@M=0pt@R=6pt@C=2pt{ & {\scriptstyle 2} & & {\scriptstyle 1}
  \\ & \ar@{-}[d] &  &\ar@{-}[d]  \\  & \ar@{-}[dl]
 \ar@{-}[dr] &  & \ar@{-}[dl]\\
\ar@{-}[d]&  & \ar@{-}[d] & \\ & & & \\ {\scriptstyle 1} & &
{\scriptstyle 2} &}} -
\vcenter{\xymatrix@M=0pt@R=6pt@C=2pt{{\scriptstyle 1} & & {\scriptstyle 2} & \\
\ar@{-}[d] & &\ar@{-}[d] &  \\ \ar@{-}[dr] &
& \ar@{-}[dl] \ar@{-}[dr]& \\
& \ar@{-}[d]&  & \ar@{-}[d] \\ & & & \\ &{\scriptstyle 1} &
&{\scriptstyle 2}}}+
\vcenter{\xymatrix@M=0pt@R=6pt@C=2pt{{\scriptstyle 2} & &
{\scriptstyle 1} &
\\ \ar@{-}[d] & &\ar@{-}[d] &  \\ \ar@{-}[dr] &
& \ar@{-}[dl] \ar@{-}[dr]& \\
& \ar@{-}[d]&  & \ar@{-}[d] \\ & & & \\ & {\scriptstyle 1} & &
{\scriptstyle 2}}}\ .
\end{array}
\right. \end{eqnarray*}

The $\BLi$-gebras are exactly the Lie bialgebras introduced by V. Drinfeld
  (\emph{cf.} \cite{Drinfeld}).

\item One variation of Lie bialgebras has been introduced by M. Chas (\emph{cf.}
\cite{C}) in the framework of String Topology. An \emph{involutive Lie bialgebra}
is a Lie bialgebra such that $\lambda \circ \Delta=0$. Therefore, the associated
properad is the same than $\BLi$, where the ``loop''
$\vcenter{\xymatrix@M=0pt@R=5pt@C=5pt{& \ar@{-}[d] &   \\
&\ar@{-}[dl]
\ar@{-}[dr] & \\
\ar@{-}[dr]& & \ar@{-}[dl]\\ &\ar@{-}[d] & \\ & & }}$ is added to the space of
relations. We denote it by $\BLi^0$, since every composition based on level graphs of
genus $>0$ is null.

\item \label{InfBi} The properad $\IBi$ coding infinitesimal Hopf
algebras. This properad is generated by the $\Sy$-bimodule
$$V=\left\{ \begin{array}{l}
V(1,\, 2)=m.k\otimes k[\Sy_2], \\
V(2,\, 1)=\Delta.k[\Sy_2]\otimes k, \\
V(m,\, n) = 0 \quad \textrm{otherwise.}
\end{array}
\right.  = \vcenter{\xymatrix@M=0pt@R=6pt@C=6pt{\ar@{-}[dr] &
&\ar@{-}[dl]  \\  &\ar@{-}[d] & \\
& & \\ & & }} \oplus
\vcenter{\xymatrix@M=0pt@R=6pt@C=6pt{  & \ar@{-}[d]&  \\
&\ar@{-}[dl] \ar@{-}[dr]& \\ & &\\&&}}\ ,$$ And its relations are given by
\begin{eqnarray*}
R&=&\left\{ \begin{array}{l}
m(m,\, 1)-m(1,\, m) \\
\oplus \ (\Delta,\,  1)\Delta - (1,\Delta)\Delta  \\
\oplus \ \Delta \otimes m - (m,\, 1)(1,\, \Delta) - (1,\,
m)(\Delta,\, 1).
\end{array}
\right. \\
&=& \left\{ \begin{array}{l} \vcenter{\xymatrix@M=0pt@R=6pt@C=6pt{
& & & & \\ \ar@{-}[ddrr] &
&\ar@{-}[dl] & & \ar@{-}[ddll]  \\
&  & &  & \\
& &\ar@{-}[d] & & \\
& & }} - \vcenter{\xymatrix@M=0pt@R=6pt@C=6pt{& & & & \\
\ar@{-}[dr] &
&\ar@{-}[dr] & & \ar@{-}[dl]  \\
& \ar@{-}[dr] & &\ar@{-}[dl]  & \\
& &\ar@{-}[d] & & \\
& & }} \\
\oplus \vcenter{\xymatrix@M=0pt@R=6pt@C=6pt{ & & & & \\
 & &\ar@{-}[d] & & \\
& & \ar@{-}[dl]  \ar@{-}[dr] & & \\
& \ar@{-}[dl]  \ar@{-}[dr] & &\ar@{-}[dr] & \\
& & & &}} -
\vcenter{\xymatrix@M=0pt@R=6pt@C=6pt{ & & & & \\
 & &\ar@{-}[d] & & \\
& & \ar@{-}[dl]  \ar@{-}[dr] & & \\
& \ar@{-}[dl]   & &\ar@{-}[dl] \ar@{-}[dr] & \\
& & & & }}  \\
\ \\
\oplus \vcenter{\xymatrix@M=0pt@R=6pt@C=6pt{\ar@{-}[dr] &  &\ar@{-}[dl]  \\  &\ar@{-}[d] & \\
& \ar@{-}[dl]\ar@{-}[dr]& \\ & & }} -
 \vcenter{\xymatrix@M=0pt@R=6pt@C=6pt{  & \ar@{-}[d] &  &\ar@{-}[d]  \\  & \ar@{-}[dl]
 \ar@{-}[dr] &  & \ar@{-}[dl]\\
\ar@{-}[d]&  & \ar@{-}[d] & \\ & & & }} -
\vcenter{\xymatrix@M=0pt@R=6pt@C=6pt{ \ar@{-}[d] & &\ar@{-}[d] &
\\ \ar@{-}[dr] &
& \ar@{-}[dl] \ar@{-}[dr]& \\
& \ar@{-}[d]&  & \ar@{-}[d] \\ & & & }}\ .
\end{array} \right. \end{eqnarray*}

These $\IBi$-gebras are associative and non-commutative analogs of Lie bialgebras. They correspond
to the infinitesimal Hopf algebras defined by M. Aguiar
(\emph{cf.} \cite{Ag1}, \cite{Ag2} and \cite{Ag3}).

\item The properad $\frac{1}{2}\Bi$, degenerated analogue of the properad of bialgebras.
The generating $\Sy$-bimodule $V$ is the same as in the case of $\IBi$,
composed by a product and a coproduct. The relations $R$ are the following ones
\begin{eqnarray*}
R&=&\left\{ \begin{array}{l}
m(m,\, 1)-m(1,\, m) \\
\oplus \ (\Delta,\,  1)\Delta - (1,\Delta)\Delta  \\
\oplus \ \Delta \otimes m \qquad \qquad \qquad (\kappa).
\end{array}
\right. \\
&=& \left\{ \begin{array}{l} \vcenter{\xymatrix@M=0pt@R=6pt@C=6pt{
& & & & \\ \ar@{-}[ddrr] &
&\ar@{-}[dl] & & \ar@{-}[ddll]  \\
&  & &  & \\
& &\ar@{-}[d] & & \\
& & }} - \vcenter{\xymatrix@M=0pt@R=6pt@C=6pt{& & & & \\
\ar@{-}[dr] &
&\ar@{-}[dr] & & \ar@{-}[dl]  \\
& \ar@{-}[dr] & &\ar@{-}[dl]  & \\
& &\ar@{-}[d] & & \\
& & }} \\
\oplus \vcenter{\xymatrix@M=0pt@R=6pt@C=6pt{ & & & & \\
 & &\ar@{-}[d] & & \\
& & \ar@{-}[dl]  \ar@{-}[dr] & & \\
& \ar@{-}[dl]  \ar@{-}[dr] & &\ar@{-}[dr] & \\
& & & &}} -
\vcenter{\xymatrix@M=0pt@R=6pt@C=6pt{ & & & & \\
 & &\ar@{-}[d] & & \\
& & \ar@{-}[dl]  \ar@{-}[dr] & & \\
& \ar@{-}[dl]   & &\ar@{-}[dl] \ar@{-}[dr] & \\
& & & & }}  \\
\ \\
\oplus \vcenter{\xymatrix@M=0pt@R=6pt@C=6pt{\ar@{-}[dr] &  &\ar@{-}[dl]  \\  &\ar@{-}[d] & \\
& \ar@{-}[dl]\ar@{-}[dr]& \\ & & }} \qquad \qquad \qquad
(\kappa).\end{array} \right.
\end{eqnarray*}
This example has been introduced by M. Markl in \cite{Markl3}
to describe a minimal model for the prop of bialgebras.
\end{itemize}
\end{ex}

\begin{cex}
$ $

\begin{itemize}
\item The previous properad is a quadratic version of the properad
$\Bi$ coding bialgebras. This last one is defined by generators
and relations but is not quadratic. The definition of $\Bi$ is the
same than $\frac{1}{2}\Bi$ where the relation $(\kappa)$ is
replaced by the following one
\begin{eqnarray*}
(\kappa')\, &:& \, \Delta \otimes m - (m,\,m)\otimes (1324)\otimes
(\Delta,\,\Delta)\\ &=&
\vcenter{\xymatrix@M=0pt@R=6pt@C=6pt{\ar@{-}[dr] &  &\ar@{-}[dl]  \\
&\ar@{-}[d] & \\
& \ar@{-}[dl]\ar@{-}[dr]& \\ & & }} -
\vcenter{\xymatrix@R=6pt@C=2pt{ &\ar@{-}[d] & &\ar@{-}[d] & \\
*=0{}& *=0{} \ar@{-}[dl] \ar@{-}[ddrr] |\hole & *=0{}&*=0{}
\ar@{-}[ddll]\ar@{-}[dr] &*=0{} \\
*=0{}\ar@{-}[dr] & *=0{}&*=0{} &*=0{} &*=0{}\ar@{-}[dl] \\
*=0{}&*=0{} \ar@{-}[d]&*=0{} &*=0{} \ar@{-}[d]&*=0{} \\
& & & & \\}}
\end{eqnarray*}
Since the relation $(\kappa')$ includes a composition involving $4$
generating operations,
the properad $\Bi$ is neither quadratic nor weight graded.

\item Consider the prop coding unitary infinitesimal
Hopf algebras defined by J.-L. Loday in
\cite{Loday4}. Its definition is the same than the prop
$\IBi$ where the third relation is
$$\Delta
\otimes m - (m,\, 1)(1,\, \Delta) - (1,\, m)(\Delta,\, 1)- (1,\,
1)= \vcenter{\xymatrix@M=0pt@R=6pt@C=6pt{\ar@{-}[dr] &  &\ar@{-}[dl]  \\  &\ar@{-}[d] & \\
& \ar@{-}[dl]\ar@{-}[dr]& \\ & & }} -
 \vcenter{\xymatrix@M=0pt@R=6pt@C=6pt{  & \ar@{-}[d] &  &\ar@{-}[d]  \\  & \ar@{-}[dl]
 \ar@{-}[dr] &  & \ar@{-}[dl]\\
\ar@{-}[d]&  & \ar@{-}[d] & \\ & & & }} -
\vcenter{\xymatrix@M=0pt@R=6pt@C=6pt{ \ar@{-}[d] & &\ar@{-}[d] &
\\ \ar@{-}[dr] &
& \ar@{-}[dl] \ar@{-}[dr]& \\ & \ar@{-}[d]&  & \ar@{-}[d] \\ & & &
}} -
\vcenter{\xymatrix@M=0pt@R=6pt@C=6pt{ \ar@{-}[d] \\
\ar@{-}[d] \\ \ar@{-}[d] \\{} }}
\ \  \vcenter{\xymatrix@M=0pt@R=6pt@C=6pt{ \ar@{-}[d] \\
 \ar@{-}[d] \\ \ar@{-}[d] \\ {} }} \ .$$
Since this relation is neither connected nor homogenous of weight $2$, one
has that this prop is not quadratic and it does not come from a properad.
\end{itemize}
\end{cex}

\section{Differential graded framework}
\label{dgframework}

In this section, we generalize the previous notions in the differential graded framework. For instance,
we define the compositions products of differential graded $\Sy$-bimodule. We give the definitions of
\emph{quasi-free $\Po$-modules} and \emph{quasi-free properads}.

\subsection{The category of differential graded $\Sy$-bimodules}

We now work in the category dg-Mod of differential graded $k$-modules. Recall
that this category can be endowed with a tensor product $V\otimes_k W$ defined
by the following formula
$$(V\otimes_k W)_d:=\bigoplus_{i+j=d}V_i\otimes_k W_j, $$
where the boundary map on $V\otimes_k W$ is given by
$$\delta(v\otimes w):=\delta(v)\otimes w  + (-1)^{|v|}v\otimes \delta(w),$$
for $v\otimes w$ an elementary tensor and where $|v|$ denotes the
homological degree of $v$. The symmetry isomorphism involve
Koszul-Quillen sign rules
$$\tau_{v,\, w}\ : \ v\otimes w \mapsto (-1)^{|v||w|}w\otimes v.$$
More generally, the commutation of two elements (morphisms,
differentials, objects, etc ...) of degree $d$ and $e$ induces a
sign $(-1)^{de}$.

\begin{dei}[dg-$\Sy$-bimodule]
A \emph{dg-$\Sy$-bimodule} $\Po$ is a collection $(\Po(m,\,
n))_{m,\, n \in \mathbb{N}}$ of differential graded modules with an action
of the symmetric group $\Sy_m$ on the left and an action of $\Sy_n$
on the right. These groups acts by morphisms of dg-modules. We denote by
$\Po_d(m,\, n)$ the sub-$(\Sy_m,\, \Sy_n)$-module composed by elements
of degree $d$ of the chain complex $\Po(m,\, n)$.
\end{dei}

\subsection{Composition products of dg-$\Sy$-bimodules}
\label{compositionproductofdg}

The composition product $\boxtimes_c$ and $\boxtimes$ of $\Sy$-bimodules
can be generalized to the differential graded framework. The definition remains
the same except that the equivalence relation $\sim$ is based on symmetry
isomorphisms and therefore involves Koszul-Quillen sign rules.
For instance, one has
\begin{eqnarray*}
\theta \otimes q_1\otimes \cdots \otimes q_i \otimes q_{i+1}
\otimes \cdots \otimes q_b \otimes \sigma \otimes p_1 \otimes
\cdots \otimes p_a \otimes \omega
\sim \\
(-1)^{|q_i||q_{i+1}|} \theta' \otimes q_1\otimes \cdots \otimes
q_{i+1} \otimes q_i \otimes \cdots \otimes q_b \otimes \sigma'
\otimes p_1 \otimes \cdots \otimes p_a\otimes \omega,
\end{eqnarray*}
where $\sigma'=(1\ldots i+1\ i\ldots\, a)_\ok\ \sigma$ and
$\theta'=\theta (1\ldots i+1\ i\ldots\, a)^{-1}_{\ol}$.

\begin{rem}
To lighten the notations, we will often omit the induced representations
$\theta$ and $\omega$ in the following of the text.
\end{rem}

We define the image of an element $q_1\otimes \cdots \otimes q_b \otimes \sigma
\otimes p_1 \otimes \cdots \otimes p_a$ of $\Qo(\ol,\,
\ok)\otimes_{\Sy_\ok} k[\Sy_\kj^c] \otimes_{\Sy_\ol} \Po(\oj,\,
\oi)$ under the boundary map $\delta$ by the following formula

\begin{eqnarray*}
&& \delta(q_1\otimes \cdots \otimes q_b \otimes \sigma \otimes p_1
\otimes \cdots \otimes p_a)= \\
&&\sum_{\beta=1}^b (-1)^{|q_1|+\cdots +|q_{\beta-1}|}q_1\otimes
\cdots \otimes \delta(q_\beta)\otimes \cdots \otimes q_b \otimes
\sigma \otimes p_1 \otimes \cdots \otimes p_a +\\
&&\sum_{\alpha=1}^a (-1)^{|q_1|+\cdots
+|q_b|+|p_1|+\cdots+|p_{\alpha-1}|}  q_1 \otimes \cdots \otimes
q_b \otimes \sigma \otimes p_1 \otimes \cdots \otimes
\delta(p_\alpha) \otimes \cdots \otimes p_a.
\end{eqnarray*}

\begin{lem}
\label{differentielle}
The differential $\delta$ is constant on equivalence classes under
the relation $\sim$.
\end{lem}

\begin{proo}
The proof is straightforward and left to the reader.
\end{proo}

One can  also generalize the concatenation product $\otimes$ to the category of
dg-$\Sy$-bimodules.

\begin{dei}[Saturated differential graded $\Sy$-bimodule]
A \emph{saturated differential graded $\Sy$-bimodule} is a dg-$\Sy$-bimodule
$\Po$ such that $\Ss(\Po)\cong\Po$. We denote this category $\textrm{sat-dg-}
\Sy\textrm{-biMod}$.
\end{dei}



One can extend Proposition~\ref{isoquotient} and
Proposition~\ref{isopartition} in the differential
graded framework.

\begin{pro}
\label{representation}
Let $\Qo$ and $\Po$ be two differential graded $\Sy$-bimodules.
Their connected composition product $\Qo\boxtimes_c \Po(m,\,n)$ is isomorphic
to the dg-$\Sy$-bimodule
$$\bigoplus_{\Theta} k[\Sy_m\big/{\Sy^B_{\ol}}]
\otimes_{\Sy_\ol} \mathcal{Q}(\ol,\,
\ok)\otimes_{\mathbb{S}_{\ok}} k[\Sy_\kj^c]
\otimes_{\mathbb{S}_{\oj}} \mathcal{P}(\oj,\, \oi)
\otimes_{\Sy_\oi} k[\Sy^B_{\oi}\big\backslash \Sy_n]$$
and to the
dg-$\Sy$-bimodule
\begin{eqnarray*}
&& \bigoplus_{\Theta'} \big( (\Pi'_1)\Qo(|\Pi'_1|,\, k_1)\otimes_k
\cdots \otimes_k (\Pi'_b)\Qo(|\Pi'_b|,\, k_b)
\big) \otimes_{\Sy_{\ok}}  k[\Sy_\kj^c] \otimes_{\Sy_{\oj}}\\
&&\quad \quad \big( (\Po(j_1,\, |\Pi_1|)(\Pi_1)\otimes_k
\cdots \otimes_k \Po(j_a,\, |\Pi_a|)(\Pi_a)
\big),
\end{eqnarray*}
\end{pro}

Once again, we have the same propositions for the composition product $\boxtimes$.
We are going to study the homological properties of the
composition product of two dg-$\Sy$-bimodules. \label{structurediffproduit}\\

Remark that the chain complex $(\Qo\boxtimes_c \Po,\, \delta)$
corresponds to the total complex of a bicomplex. We define the
bidegree of an element
$$(q_1,\, \ldots,\, q_b)\,
\sigma\, (p_1,\ldots,\, p_a) \quad  \textrm{of} \quad
k[\Sy_m\big/{\Sy^B_{\ol}}] \otimes_{\Sy_\ol}\Qo(\ol,\,
\ok)\otimes_{\Sy_{\ok}} k[\Sy_\kj^c] \otimes_{\Sy_{\oj}} \Po (\oj,\,
\oi)\otimes_{\Sy_\oi} k[\Sy^B_{\oi}\big\backslash \Sy_n]$$ by
$(|q_1|+\cdots+|q_b|,\, |p_1|+\cdots+|p_a|)$. The horizontal
differential $\delta_h \, :\, \Qo\boxtimes_c \Po \to \Qo\boxtimes_c
\Po$ is induced by the one of $\Qo$ and the vertical differential
$\delta_v \, :\, \Qo\boxtimes_c \Po \to \Qo\boxtimes_c \Po$ is induced
by the one of $\Po$. One has explicitly,

\begin{eqnarray*}
&& \delta_h\big((q_1,\, \ldots,\, q_b)\, \sigma\, (p_1,\ldots,\,
p_a) \big)= \\
&&\sum_{\beta=1}^b (-1)^{|q_1|+\cdots+|q_{\beta-1}|}(q_1,\ldots,\,
\delta(q_\beta),\ldots,\, q_b)\, \sigma\, (p_1,\ldots,\, p_a)\\
\textrm{and} && \delta_v\big((q_1,\, \ldots,\, q_b)\, \sigma\,
(p_1,\ldots,\, p_a) \big)=\\
&& \sum_{\alpha=1}^a
(-1)^{|q_1|+\cdots+|q_b|+|p_1|+\cdots+|p_{\alpha-1}|}(q_1,\ldots,\,
q_b)\, \sigma\, (p_1,\ldots,\,\delta(p_\alpha),\ldots,\, p_a).
\end{eqnarray*}

One can see that $\delta=\delta_h+\delta_v$ and that
$\delta_h\circ \delta_v =-\delta_v\circ \delta_h$.

Like every bicomplex, this one gives rise to two spectral
sequences $I^r(\Qo\boxtimes \Po)$ et $II^r(\Qo\boxtimes \Po)$ that
converge to the total homology $H_*(\Qo\boxtimes \Po, \, \delta)$.
Recall that these two spectral sequences verify

\begin{eqnarray*}
 && I^1(\Qo\boxtimes_c \Po)=H_*(\Qo\boxtimes_c \Po, \, \delta), \qquad
I^2(\Qo\boxtimes_c \Po)=H_*(H_*(\Qo\boxtimes_c \Po, \,
\delta_v),\,\delta_h) \\
&and& II^1(\Qo\boxtimes_c \Po)=H_*(\Qo\boxtimes_c \Po, \, \delta_h),
\qquad II^2(\Qo\boxtimes_c \Po)=H_*(H_*(\Qo\boxtimes_c \Po, \,
\delta_h),\, \delta_v).
\end{eqnarray*}

Proposition~\ref{representation} gives that the composition
product $\Qo \boxtimes_c \Po$ can be decomposed as the direct sum of
chain complexes
$$ \bigoplus_{\Theta}
k[\Sy_m\big/{\Sy^B_{\ol}}] \otimes_{\Sy_\ol}\Qo(\ol,\,
\ok)\otimes_{\Sy_{\ok}} k[\Sy^c_\kj] \otimes_{\Sy_{\oj}} \Po (\oj,\,
\oi)\otimes_{\Sy_\oi} k[\Sy^B_{\oi}\big\backslash \Sy_n].$$

This decomposition is compatible with the structure of bicomplex.
Therefore, the spectral sequences can be decomposed in the same
way

\begin{eqnarray*}
&& I^r(\Qo\boxtimes_c \Po)=\bigoplus_{\Theta}
I^r\left(k[\Sy_m\big/{\Sy^B_{\ol}}] \otimes_{\Sy_\ol}\Qo(\ol,\,
\ok)\otimes_{\Sy_\ok} k[\Sy_\kj^c] \otimes_{\Sy_\oj} \Po(\oj,\, \oi)
\otimes_{\Sy_\oi} k[\Sy^B_{\oi}\big\backslash \Sy_n]
\right) \\
&\textrm{and}&II^r(\Qo\boxtimes_c \Po)=\bigoplus_{\Theta}
II^r\left(k[\Sy_m\big/{\Sy^B_{\ol}}] \otimes_{\Sy_\ol}\Qo(\ol,\,
\ok)\otimes_{\Sy_\ok} k[\Sy_\kj^c] \otimes_{\Sy_\oj} \Po(\oj,\, \oi)
\otimes_{\Sy_\oi} k[\Sy^B_{\oi}\big\backslash \Sy_n]\right)
\end{eqnarray*}

From these two formulas, one has the following proposition.

\begin{pro}
\label{I1} When $k$ is a field of characteristic $0$, one has
\begin{eqnarray*}
&&I^1\left(k[\Sy_m\big/{\Sy^B_{\ol}}] \otimes_{\Sy_\ol}\Qo(\ol,\,
\ok)\otimes_{\Sy_\ok} k[\Sy_\kj^c] \otimes_{\Sy_\oj} \Po(\oj,\,
\oi)\otimes_{\Sy_\oi} k[\Sy^B_{\oi}\big\backslash \Sy_n]
\right)\\
&=& k[\Sy_m\big/{\Sy^B_{\ol}}] \otimes_{\Sy_\ol}\Qo(\ol,\,
\ok)\otimes_{\Sy_\ok}k[\Sy_\kj^c]\otimes_{\Sy_\oj} H_*\left(\Po(\oj,\,
\oi)\right)\otimes_{\Sy_\oi}
k[\Sy^B_{\oi}\big\backslash \Sy_n] \ \textrm{and} \\
&&II^1\left(k[\Sy_m\big/{\Sy^B_{\ol}}] \otimes_{\Sy_\ol}\Qo(\ol,\,
\ok)\otimes_{\Sy_\ok} k[\Sy_\kj^c] \otimes_{\Sy_\oj} \Po(\oj,\, \oi)
\otimes_{\Sy_\oi} k[\Sy^B_{\oi}\big\backslash
\Sy_n]\right)\\
&=& k[\Sy_m\big/{\Sy^B_{\ol}}]
\otimes_{\Sy_\ol}H_*\left(\Qo(\ol,\,
\ok)\right)\otimes_{\Sy_\ok}k[\Sy_\kj^c]\otimes_{\Sy_\oj} \Po(\oj,\,
\oi)\otimes_{\Sy_\oi} k[\Sy^B_{\oi}\big\backslash \Sy_n].
\end{eqnarray*}
\end{pro}

\begin{proo}
By Mashke's Theorem, we have the rings of the form $k[\Sy_\oj]$
and $k[\Sy_\oi]$ are semi-simple. Therefore, the modules
$k[\Sy_m\big/{\Sy^B_{\ol}}] \otimes_{\Sy_\ol}\Qo(\ol,\,
\ok)\otimes_{\Sy_\ok} k[\Sy_\kj^c]$ are $k[\Sy_\oj]$ projective
modules (on the right) and the modules
$k[\Sy^B_{\oi}\big\backslash \Sy_n]$ are $k[\Sy_\oi]$ projective
modules (on the left) .
\end{proo}

\begin{pro}
\label{I2II2} When $k$ is a field of characteristic $0$, one has
\begin{eqnarray*}
&& I^2\left(k[\Sy_m\big/{\Sy^B_{\ol}}] \otimes_{\Sy_\ol}
\Qo(\ol,\, \ok)\otimes_{\Sy_\ok} k[\Sy_\kj^c] \otimes_{\Sy_\oj}
\Po(\oj,\, \oi)\otimes_{\Sy_\oi} k[\Sy^B_{\oi}\big\backslash
\Sy_n] \right)\\ &=& II^2\left(k[\Sy_m\big/{\Sy^B_{\ol}}]
\otimes_{\Sy_\ol}\Qo(\ol,\, \ok)\otimes_{\Sy_\ok} k[\Sy_\kj^c]
\otimes_{\Sy_\oj} \Po(\oj,\, \oi)\otimes_{\Sy_\oi}
k[\Sy^B_{\oi}\big\backslash \Sy_n] \right)\\ &=&
k[\Sy_m\big/{\Sy^B_{\ol}}] \otimes_{\Sy_\ol}H_*\Qo(\ol,\,
\ok)\otimes_{\Sy_\ok}k[\Sy_\kj^c]\otimes_{\Sy_\oj} H_*\Po(\oj,\,
\oi)\otimes_{\Sy_\oi} k[\Sy^B_{\oi}\big\backslash \Sy_n].
\end{eqnarray*}
\end{pro}

\begin{proo}
We use the same arguments to show that
\begin{eqnarray*}
&& I^2\left(k[\Sy_m\big/{\Sy^B_{\ol}}] \otimes_{\Sy_\ol}
\Qo(\ol,\, \ok)\otimes_{\Sy_\ok} k[\Sy_\kj^c] \otimes_{\Sy_\oj}
\Po(\oj,\, \oi) \otimes_{\Sy_\oi} k[\Sy^B_{\oi}\big\backslash
\Sy_n] \right) \\
&=& k[\Sy_m\big/{\Sy^B_{\ol}}]
\otimes_{\Sy_\ol}H_*\left(\Qo(\ol,\,
\ok)\right)\otimes_{\Sy_\ok}k[\Sy_\kj^c]\otimes_{\Sy_\oj}
H_*\left(\Po(\oj,\, \oi)\right)\otimes_{\Sy_\oi}
k[\Sy^B_{\oi}\big\backslash \Sy_n] .
\end{eqnarray*}
We conclude by K\"unneth's formula
\begin{eqnarray*}
H_*\left(\Qo(\ol,\, \ok)\right)&=& H_*\left(\Qo(l_1,\, k_1)\otimes
\cdots \otimes \Qo(l_a,\, k_a) \right) \cong \\ && H_*\Qo(l_1,\,
k_1)\otimes \cdots \otimes H_*\Qo(l_a,\, k_a) = H_*\Qo(\ol,\,
\ok).
\end{eqnarray*}
\end{proo}

\begin{rem}
Let $\Phi \, : \, M \to M'$ and $\Psi \, : \, N \to N'$ be two
morphisms of dg-$\Sy$-bimodules. The morphism $\Phi\boxtimes_c \Psi
\, :\, M\boxtimes_c N \to M' \boxtimes_c N'$ is a morphism of
bicomplexes. It induces morphisms of spectral sequences
$I^r(\Phi\boxtimes_c \Psi)\, :\, I^r(M\boxtimes_c N) \to
I^r(M'\boxtimes_c N')$ and $ II^r(\Phi\boxtimes_c \Psi)\, :\,
II^r(M\boxtimes_c N) \to II^r(M'\boxtimes_c N')$.
\end{rem}

More generally, we have the following proposition.

\begin{pro}
\label{homologieproduit} When $k$ is a field of characteristic
$0$, one has
$$H_*(\Qo \boxtimes_c \Po)(m,\, n)= \bigoplus_{\Theta}
k[\Sy_m\big/{\Sy^B_{\ol}}] \otimes_{\Sy_\ol}H_*\Qo(\ol,\,
\ok)\otimes_{\Sy_\ok} k[\Sy_\kj^c] \otimes_{\Sy_\oj} H_*\Po(\oj,\,
\oi) \otimes_{\Sy_\oi} k[\Sy^B_{\oi}\big\backslash \Sy_n].$$
That's-to-say, $H_*(\Qo \boxtimes_c \Po) =(H_*\Qo) \boxtimes_c
(H_*\Po)$.
\end{pro}

\begin{proo}
It is a direct corollary of Mashke's Theorem (for coinvariants)
and K\"un\-neth's Theorem (for tensor products).
\end{proo}

Once again, the similar results apply to the connected
composition products $\boxtimes$.

\begin{rem}
These propositions generalize the results of B. Fresse for the
composition product $\circ$ of operads to the composition products
$\boxtimes_c$ and $\boxtimes$ of properads and props (\emph{cf.}
\cite{Fresse} 2.3.).
\end{rem}

\subsection{Differential properads and differential props}

We give here a differential variation of the notions of properads and
props.

\begin{dei}[Differential properad]
A \emph{differential properad $(\Po,\, \mu,\, \eta)$} is a monoid in
the monoidal category $(\textrm{dg-}\Sy\textrm{-biMod},\,
\boxtimes_c,\, I) $ of differential graded
$\Sy$-bimodules.
\end{dei}

The composition morphism $\mu$ of a differential properad is a
morphism of chain complexes of degree $0$ compatible with the
action of the symmetric groups $\Sy_m$ and $\Sy_n$. One has the
derivation type relation :
\begin{eqnarray*}
&& \delta \left( \mu((p_1,\ldots,\, p_b)\, \sigma \,
(p'_1,\ldots,\, p'_a) ) \right) = \\
&& \sum_{\beta=1}^b (-1)^{|p_1|+\cdots +|p_{\beta-1}|} \mu((
p_1,\ldots ,\, \delta(p_\beta),\ldots ,\,
p_b) \, \sigma \, (p'_1,\ldots ,\,p'_a)) +\\
&& \sum_{\alpha=1}^a (-1)^{|p_1|+\cdots
+|p_b|+|p'_1|+\cdots+|p'_{\alpha-1}|}  \mu((p_1,\ldots,\, p_b) \,
\sigma \, (p'_1,\ldots,\, \delta(p'_\alpha),\ldots,\,p'_a)).
\end{eqnarray*}

Proposition~\ref{homologieproduit} gives the following property.

\begin{pro}
When $k$ is a field of characteristic $0$, the homology groups
$(H_*(\Po),\ H_*(\mu),\, H_*(\eta))$ of a differential properad
$(\Po,\, \mu,\, \eta)$ form a graded properad.
\end{pro}

\begin{dei}[Differential weight graded $\Sy$-bimodule]
A \emph{differential weight graded $\Sy$-bimodule} $M$ is a direct
sum over $\rho \in \mathbb{N}$ of differential $\Sy$-bimodules
$M=\bigoplus_{\rho\in \mathbb{N}} M^{(\rho)}$.

The morphisms of differential weight graded $\Sy$-bimodules are
the morphisms of differential $\Sy$-bimodules that preserve this
decomposition. This forms a category denoted gr-dg-$\Sy$-biMod.
\end{dei}

Remark that the weight is a graduation separate from the
homological degree that gives extra information.

\begin{dei}[Differential weight graded properad]
A \emph{differential weight graded properad} is a monoid in the
monoidal category of differential weight graded
$\Sy$-bimodules.
\end{dei}

A differential weight graded properad $\Po$ such that $\Po^{(0)}=I$ is
called \emph{connected}.

\begin{rem}
Dually, we can define the notion of \emph{differential coproperad} and
\emph{weight graded differential coproperad}. A differential coproperad is
a coproperad where the coproduct is a morphism of dg-$\Sy$-bimodules
verifying a coderivation type relation.
\end{rem}

Similarly one has the notions of \emph{differential weight graded
prop} and \emph{coprop}.

\subsection{Free and quasi-free $\Po$-modules}

We define a variation of the notion of free $\Po$-module, namely the notion of
\emph{quasi-free $\Po$-module}, and we prove some of its basic properties. \\

Let $\Po$ be a (weight graded) dg-properad. Recall from \ref{P-module}, that the
free differential $\Po$-module on the right over a dg-$\Sy$-bimodule $M$ is given by
$L=M\boxtimes_c \Po$, where the differential $\delta$ comes from the one of $M$ and the one of
$\Po$ as explained before. When $(\Po,\, \mu, \, \eta,\, \epsilon)$ is an augmented dg-properad,
one has that the indecomposable quotient of the free dg-$\Po$-module
$L=M\boxtimes_c \Po$ is isomorphic to $M$, as dg-$\Sy$-bimodule.\\

In the following part of the text, we will encounter differential
$\Po$-modules which are not strictly free. They are built on the
dg-$\Sy$-bimodule $M\boxtimes_c \Po$ but the boundary map is the sum of the
canonical differential $\delta$ with an homogenous morphism of degree $-1$.

\begin{dei}[Quasi-free $\Po$-module]
Let $\Po$ be differential properad. A \emph{(right) quasi-free
$\Po$-module} $L$ is a dg-$\Sy$-bimodule of the form $M\boxtimes_c
\Po$, where the differential $\delta_\theta$ is the sum
$\delta+d_\theta$ of the canonical differential $\delta$ on the
product $M\boxtimes_c \Po$ with an homogenous morphism $d_\theta \,
: \, M\boxtimes_c \Po \to M\boxtimes_c \Po$ of degree $-1$
 such that
\begin{eqnarray*}
&& d_\theta\left((m_1,\ldots,\, m_b)\, \sigma \, (p_1,\ldots,\,
p_a) \right) =\\
&& \sum_{\beta=1}^b (-1)^{|m_1|+\cdots+|m_{\beta-1}|}
(m_1,\ldots,\, d_\theta(m_\beta),\ldots,\, m_b)\, \sigma \,
(p_1,\ldots,\, p_a).
\end{eqnarray*}

In the weight graded framework, we ask $d_\theta$ to be an
homogenous map for the weight such that
$${d_\theta} \ : \ M^{(\rho)} \to \big( \underbrace{M}_{(<\rho)}
\boxtimes_c \underbrace{\Po}_{(>1)} \big)^{(\rho)}.$$
In this case, the boundary map $\delta_\theta$ preserves the global weight of $L$.
\end{dei}

Therefore, the boundary map $\delta_\theta \, : \, M\boxtimes_c \Po \to M\boxtimes_c \Po$
verifies the following relation
\begin{eqnarray*}
&& \delta_\theta\left((m_1,\ldots,\, m_b)\, \sigma \,
(p_1,\ldots,\,
p_a) \right) =\\
&& \sum_{\beta=1}^b (-1)^{|m_1|+\cdots+|m_{\beta-1}|}
(m_1,\ldots,\, \delta_\theta(m_\beta),\ldots,\, m_b)\, \sigma \,
(p_1,\ldots,\, p_a) +\\
&& \sum_{\alpha=1}^a
(-1)^{|m_1|+\cdots+|m_b|+|p_1|+\cdots+|p_{\alpha-1}|}
(m_1,\ldots,\, m_b)\, \sigma \, (p_1,\ldots,\, \delta(p_\alpha)
,\ldots,\, p_a).
\end{eqnarray*}

Since $\delta^2=0$, the identity ${\delta_\theta}^2=0$ is
equivalent to $\delta d_\theta + d_\theta \delta +
{d_\theta}^2=0$. \\

If one writes $d_\theta(m_\beta)=(m'_1,\ldots,\, m'_{b'})\,
\sigma'\,(p'_1,\ldots,\, p'_{a'})\in M\boxtimes_c \Po$, then one has
\begin{eqnarray*}
&& (m_1,\ldots,\, d_\theta(m_\beta),\ldots,\, m_b)\, \sigma \,
(p_1,\ldots,\, p_a) \\
&& (m_1,\ldots,\, m'_1,\ldots,\, m'_{b'},\ldots ,\ m_b)\, \sigma' \,
\mu\left( (p'_1,\ldots,\, p'_{a'}) \, \sigma'' \, (p_1,\ldots,\,
p_a) \right) \in M\boxtimes_c \Po.
\end{eqnarray*}

Notice that the morphism $d_\theta$ is determined by its restriction on
$\xymatrix{M \ar[r]^-{M\boxtimes_c \eta}&
M\boxtimes_c \Po}$.

\begin{pro}
The morphism $\xymatrix{M \ar[r]^-{M\boxtimes_c \eta}& M\boxtimes_c
\Po}$ is a morphism of dg-$\Sy$-bimodules if and only if $d_\theta=0$.
\end{pro}

When $\Po$ is a augmented dg-properad, the condition for the morphism $M\boxtimes_c
\varepsilon$ to be a morphism of
differential graded modules is weaker.

\begin{pro}
Let $(\Po,\, \mu,\, \eta,\, \varepsilon)$ be an augmented dg-properad.
The morphism  $\xymatrix{M\boxtimes_c\Po
\ar[r]^-{M\boxtimes_c \varepsilon}& M}$ is a morphism of
dg-$\Sy$-bimodules is and only if the following composition is null
$$\xymatrix{ M\boxtimes_c \Po \ar[r]^-{d_\theta}& M \boxtimes_c
\Po \ar[r]^-{M\boxtimes_c \varepsilon}& M }=0.$$
\end{pro}

\begin{cor}
In this case, the indecomposable quotient of $L=M\boxtimes_c \Po$
is isomorphic to $M$ as a dg-$\Sy$-bimodule.
\end{cor}

\begin{dei}[Morphism of weight graded quasi-free $\Po$-module]
Let $L=M\boxtimes_c \Po$ and $L'=M'\boxtimes_c \Po$ be two weight graded
quasi-free $\Po$-modules. A morphism $\Phi \, : \, L \to L'$ of
weight graded differential $\Po$-modules is a morphism of
\emph{weight graded quasi-free $\Po$-modules} if $\Phi$ preserves the graduation
on the left
$$\Phi \ :\ \underbrace{M}_{(\rho)} \boxtimes_c \Po \to
\underbrace{M'}_{(\rho)} \boxtimes_c \Po.$$
\end{dei}

In the following of the text, the examples of weight graded
quasi-free $\Po$-modules that we will have to treat will have a particular form.

\begin{dei}[Analytic quasi-free $\Po$-module]
\label{analyticquasifreemodule}
Let $\Po$ be a weight graded differential properad. An
\emph{analytic quasi-free $\Po$-module} is quasi-free $\Po$-module $L=M\boxtimes_c \Po$
such that the weight graded differential $\Sy$-bimodule $M$ is given by the image $M=\Upsilon(V)$
of a weight graded differential $\Sy$-bimodule $V$ under a split analytic functor $\Upsilon$
(\emph{cf.} \ref{splitanalyticfunctor}).
Denote $\Upsilon=\bigoplus_{n\in \mathbb{N}} \Upsilon_{(n)}$, where
$\Upsilon_{(n)}$ is an homogenous polynomial functor of degree $n$ and such that
$\Upsilon_{(0)}=I$. Moreover, $V$ must verify $V^{(0)}=0$ so that
$\Upsilon(V)^{(0)}=I$. As we have seen in Proposition~\ref{bigraduation}, $\Upsilon(V)$
is bigraded. The morphism $d_\theta$ is supposed to verify
$$d_\theta \ : \ \Upsilon_{(n)}(V)^{(\rho)} \to \underbrace{\Upsilon(V)}_{n-1}\boxtimes_c
\underbrace{\Po}_{1},$$ where the graduation $(\rho)$ is the global graduation coming from the one
of $V$ and $(n)$ is the degree of the analytic functor $\Upsilon$.
\end{dei}

\begin{ex}
The examples of analytic quasi-free $\Po$-modules fall into two parts.
\begin{itemize}
\item When $V=\oPo$, we have the example of the (differential) bar construction
$\Upsilon(V)=\bar{\B}(\Po)=\F^c(\Sigma\oPo)$ (\emph{cf.} section~\ref{barcobar}).
The augmented bar construction
$\bar{\B}(\Po)\boxtimes_c \Po$ is an analytic quasi-free $\Po$-modules.

\item When $\Po$ is quadratic properad generated by a
dg-$\Sy$-bimodule $V$, the Koszul dual coproperad $\Po^{\ac}$
is a split analytic functor in $V$ (\emph{cf.} section~\ref{KoszulDuality}).
Therefore, the Koszul complex $\Po^{\ac}\boxtimes_c \Po$ is an analytic
quasi-free $\Po$-module.
\end{itemize}
\end{ex}

We have immediately the following proposition.

\begin{pro}
Let $\Po$ be an augmented weight graded dg-properad and $L=\Upsilon(V)\boxtimes_c \Po$
be an analytic quasi-free $\Po$-module. We have
$$\xymatrix{ M\boxtimes_c \Po \ar[r]^-{d_\theta}& M \boxtimes_c \Po
\ar[r]^-{M\boxtimes_c \varepsilon}& M }=0.$$
Therefore, we can identify the indecomposable quotient of $L$ to $\Upsilon(V)$.
\end{pro}

\begin{rem}
Dually, we have the notions of differential cofree $\Co$-comodule on a dg-$\Sy$-bimodule $M$,
given by $M\boxtimes_c \Co$, quasi-cofree $\Co$-comodule and analytic quasi-cofree $\Co$-comodule.
The augmented cobar construction $\bar{\B}^c(\Co)\boxtimes_c \Co$ is an example of analytic
quasi-free $\Co$-comodule.
\end{rem}

\subsection{Free and quasi-free properads}

We define here the notion of free differential properad. As in the
case of $\Po$-modules, we will have to deal with differential
properads that are not strictly free. Therefore, we introduce the
notion of \emph{quasi-free properad} which is a free properad where the
differential is given by the sum of the canonical one with a
\emph{derivation}.

Let $(V,\, \delta)$ be a dg-$\Sy$-bimodule. We have seen in
\ref{free} that the free properad on $V$ is given by the weight graded
$\Sy$-bimodule $\F(V)$ which is the direct sum on connected graphs (without
level) where the vertices are indexed by elements of $V$,
$$ \F(V)=\left( \bigoplus_{g\in\mathcal{G}_c} \bigotimes_{\nu\in
\mathcal{N}}V(|Out(\nu)|,\, |In(\nu)|)\right) \Bigg/ \approx \ .$$

This construction corresponds to the quotient of the simplicial
$\Sy$-bimodule $\mathcal{FS}(V)=\bigoplus_{n\in\mathbb{N}}
(V_+)^{\boxtimes_c n}$ under the relation generated by $I\boxtimes_c
V \equiv V \boxtimes_c I$. In the differential graded framework,
one has to take the signs into account. In this case, we quotient
$\FS(V)$ by the relation $$(\nu_1,\, \nu,\, \nu_2)\, \sigma \,
(\nu'_1,\, 1,\, \nu'_2)\equiv (-1)^{|\nu_2|+|\nu'_1|}(\nu_1,\, 1,\, 1,\, \nu_2)\, \sigma' \,
(\nu'_1,\, \nu,\, \nu'_2)$$, for instance, where $\nu$ belongs to
$V(2,\, 1)$.

\begin{pro}
The canonical differential $\delta$ defined on $\FS(V)$ by the one
of $V$ is constant on the equivalence classes for the relation
$\equiv$.
\end{pro}

\begin{proo}
The verification of signs is straightforward and is the same as in
 Lemma~\ref{differentielle}.
\end{proo}

\begin{thm}
The differential properad $(\F(V),\, \delta)$ is free on the
dg-$\Sy$-bimodule $(V,\, \delta)$.
\end{thm}

\begin{rem}
Since the analytic decomposition of the functor $\F(V)$ given in Proposition~\ref{Fanalytique}
is stable
under the differential $\delta$, we have that $\F(V)$ is a weight graded differential properad.

The projection $\F(V) \to I$ is a morphism of dg-properads.
\end{rem}

The functor $\F$ shares interesting homological properties.

\begin{pro}
\label{freefuncteurexact}
Over a field $k$ of characteristic $0$, the functor $\F \,  : \, \textrm{dg-}\Sy\textrm{-biMod}
\to \textrm{dg-properad}$ is exact. We have $H_*(\F(V))=\F(H_*(V))$.
\end{pro}

\begin{proo}
The proof of this proposition is the same as the proof of Theorem $4$ of the article of
M. Markl and A. A. Voronov \cite{MV} and lies, once again, on Mashke's theorem.
\end{proo}

\begin{dei}[Derivation]
Let $\Po$ be a dg-properad. An homogenous morphism of $\Sy$-bimodules
$d \, :\, \Po \to \Po$ of (homological) degree $|d|$ is called a \emph{derivation} if
it satisfies the following identity
\begin{eqnarray*}
&& \delta \left( \mu((p_1,\ldots,\, p_b)\, \sigma \,
(p'_1,\ldots,\, p'_a) ) \right) = \\
&& \sum_{\beta=1}^b (-1)^{(|p_1|+\cdots +|p_{\beta-1}|)|d|} \mu((
p_1,\ldots ,\, \delta(p_\beta),\ldots ,\,
p_b) \, \sigma \, (p'_1,\ldots ,\,p'_a)) +\\
&& \sum_{\alpha=1}^a (-1)^{(|p_1|+\cdots
+|p_b|+|p'_1|+\cdots+|p'_{\alpha-1}|)|d|}  \mu((p_1,\ldots,\, p_b)
\, \sigma \, (p'_1,\ldots,\, \delta(p'_\alpha),\ldots,\,p'_a)).
\end{eqnarray*}
\end{dei}

For example, the differential $\delta$ of a differential properad is a derivation of degree $-1$
such that $\delta^2=0$. Notice that the sum $\delta+d$ of the differential $\delta$ with
a derivation $d$ of degree $-1$ is a derivation of the same degree. In this case, since
$\delta^2=0$ the equation $(\delta^+d)^2$ is equivalent to $\delta d +d \delta +d^2=0$.

\begin{rem}
Dually, one can define the notion of \emph{coderivation} on a differential coproperad $\Co$.
\end{rem}

The next proposition gives the form of the derivations on the free differential properad $\F(V)$.

\begin{pro}
\label{derivation}
Let $\F(V)$ be the free dg-properad on the
dg-$\Sy$-bimodule $V$. For any homogenous morphism $\theta \, :
\, V \to \F(V)$, there exists a unique derivation $d_\theta \, :\,
\F(V) \to \F(V)$ with the same degree such that its restriction to
$V$ is equal to $\theta$. This correspondence is one-to-one.

Moreover, if  $\theta\, : \, V \to\F_{(r)}(V)\subset \F(V)$, we have
$d_\theta\left( \F_{(s)}(V)\right)\subset \F_{(r+s-1)}(V)$.
\end{pro}

\begin{proo}
Let $g$ be a graph with $n$ vertices. We define $d_\theta$ on a
representative of an equivalence class for the relation $\equiv$
(vertical move of a vertex up to sign) associated to the graph
$g$. This corresponds to a choice of an order for writing the
elements of $V$ that index the vertices of $g$. Write this
representative $\bigotimes_{i=1}^n \nu_i$, with $\nu_i\in V$.
Hence, we define $d_\theta$ by the following formula
$$d_\theta(\bigotimes_{i=1}^n \nu_i)=\sum_{i=1}^n (-1)^{|\nu_1|
+\cdots+|\nu_{i-1}|}\left( \bigotimes_{j=1}^{i-1}\nu_j
\right)\otimes \theta(\nu_i) \otimes \left( \bigotimes_{j=i+1}^{n}
\nu_j \right).$$
One can see that $d_\theta$ is well defined on the equivalence classes for the
relation $\equiv$. Applying $d_\theta$ is given by applying $\theta$ to
every element of $V$ indexing a vertex of $g$. (Since the application $\theta$
is an application of $\Sy$-bimodules, the morphism $d_\theta$ is well defined on
the equivalence classes for the relation $\approx$).

The product $\mu$ of the free properad $\F(V)$ is surjective. Hence, it shows that
every derivation is of this form.

If the image under $\theta$ of the elements of $V$ is a sum of elements of $\F(V)$
that can be written by graphs with $r$ vertices ($\theta$ creates $r-1$ vertices), the form
of $d_\theta$ shows that $d_\theta\left( \F_{(s)}(V)\right)\subset \F_{(r+s-1)}(V)$.
\end{proo}

Dually, we have the same proposition for the coderivations of the cofree
connected coproperad.

\begin{pro}
\label{coderivation}
Let $\F^c(V)$ be the cofree connected dg-coproperad on the dg-$\Sy$-bimodule $V$.
For any homogenous morphism $\theta \, : \, \F^c(V) \to V$, there exists a unique coderivation
$d_\theta \, :\, \F^c(V) \to \F^c(V)$ such that the composition
$$\xymatrix{\F^c(V) \ar[r]^{d_\theta} & \F^c(V) \ar[r]^{proj} & V}$$ is equal to $\theta$.
This correspondence is one-to-one.

Moreover, if $\theta \, : \, \F^c(V) \to V$ is null on the components $\F^c_{(r)}(V)\subset \F^c(V)$,
for $s\neq r$, we have $d_\theta\left( \F_{(s+r-1)}(V)\right)\subset \F_{(s)}(V)$, for all $s>0$.
\end{pro}

Graphically, applying the coderivation $d_\theta$ on an element of $\F^c(V)$ represented
by a graph $g$ corresponds to apply $\theta$ on every admissible sub-graph of $g$.\\

We can now state the definition of a \emph{quasi-free properad}.

\begin{dei}(Quasi-free properad)
A \emph{quasi-free properad} $(\Po,\, \delta_\theta)$ is a dg-properad of the form
$\Po=\F(V)$, with $(V,\, \delta)$ a dg-$\Sy$-bimodule and where the differential
$\delta_\theta \, :\, \F(V) \to \F(V)$ is the sum of the canonical differential
$\delta$ with a derivation $d_\theta$ of degree $-1$, $\delta_\theta=\delta +d_\theta$.
\end{dei}

The previous proposition shows that every derivation $d_\theta$ on $\F(V)$
is determined by its restriction $\theta$ on $V$. Therefore, the inclusion
$V \to \F(V)$ is a morphism of dg-$\Sy$-bimodules if and only if
$\theta$ is null. On the contrary, the projection
$\F(V) \to V$ is a morphism of dg-$\Sy$-bimodules under a weaker condition.

\begin{pro}
The projection $\F(V) \to V$ is a morphism of dg-$\Sy$-bimodules
if and only if $\theta(V) \subset \bigoplus_{r\ge
2}\F_{(r)}(V)$.
\end{pro}

In this case, we say that the differential is \emph{decomposable}.

Dually, we have the notion of \emph{quasi-cofree coproperad}.

\begin{dei}[Quasi-cofree coproperad]
A \emph{quasi-cofree coproperad} $(\Co,\, \delta_\theta)$ is a dg-coproperad
of the form $\Co=\F^c(V)$, with $(V,\, \delta)$ a dg-$\Sy$-bimodule and
where the differential $\delta_\theta \, :\, \F^c(V) \to \F^c(V)$ is the sum
of the canonical differential $\delta$ with a coderivation
$d_\theta$ of degree $-1$, $\delta_\theta=\delta +d_\theta$.
\end{dei}

The previous proposition on coderivations of $\F^c(V)$ shows that
the projection $\F^c(V)\to V$ is a morphism of dg-$\Sy$-bimodules
if and only if $\theta$ is null and that the inclusion $V \to \F^c
(V)$ is a morphism of dg-$\Sy$-bimodules if and only if $\theta$
is null on $\F^c_{(1)}(V)=V$.

\begin{ex}
The fundamental examples of quasi-free properad and quasi-cofree
coproperad are given by the cobar construction $\F\left(
\Sigma^{-1} \overline{\mathcal{C}} \right)$ and the bar
construction $\F^c\left( \Sigma \oPo \right)$. These two
constructions are studied in the next section.
\end{ex}

\section{Bar and cobar construction}
\label{barcobar}

We generalize the bar and cobar constructions of associative algebras (\emph{cf.} S\'eminaire
H. Cartan \cite{Cartan}) and for operads (\emph{cf.} E. Getzler, J. Jones \cite{GJ}) to properads.
First, we study the properties of the partial composition products and coproducts.
We define the reduced bar and cobar constructions and then the bar and cobar constructions
 with coefficients. These constructions share interesting homological properties. For instance,
 they induce resolutions for differential properads and props (\emph{cf.} Theorem~\ref{criteredeKoszul} and
Section~\ref{KoszulDuality}). To show these resolutions, we prove that
the augmented bar and cobar constructions are acyclic, at the end of this section.

\subsection{Partial composition product and coproduct}

We define the notions of \emph{suspension} and \emph{desuspension} of a dg-$\Sy$-bimodules.
Then, we give the definitions and first properties of the \emph{partial composition products} and
\emph{coproducts}.

\subsubsection{\bf Suspension and desuspension of a dg-$\Sy$-bimodule}

Let $\Sigma$ be the graded $\Sy$-bimodule defined by
$$\left\{ \begin{array}{ll}
\Sigma(0,\,0) = k.\Sigma &\textrm{where $\Sigma$ is an element of degree +1,} \\
\Sigma(m,\, n)=0 & \textrm{otherewise.}
\end{array}\right.$$

\begin{dei}[Suspension $\Sigma V$]
The \emph{suspension} of the dg-$\Sy$-bimodule $V$ is the dg-$\Sy$-bimodule
$\Sigma V := \Sigma \otimes V$.
\end{dei}

Taking the suspension of a dg-$\Sy$-bimodule $V$ corresponds to
taking the tensor product of every element of $V$ with $\Sigma$.
To an element $v\in V_{d-1}$, one associates $s\otimes v \in
(\Sigma V)_d$, which is often denoted $\Sigma v$. Therefore,
$(\Sigma V)_d$ is naturally isomorphic to $V_{d-1}$. Following the
same rules as in the previous section, the differential on $\Sigma
V$ is given by the formula $\delta (\Sigma v)= -\Sigma \delta
(v)$, for every $v$ in $V$. The suspension $\Sigma V$ corresponds
to the introduction of an extra element of degree $+1$. This
involves Koszul-Quillen
sign rules when permutating objects.\\

In the same way, we define the dg-$\Sy$-bimodule $\Sigma^{-1}$ by
$$\left\{ \begin{array}{ll}
\Sigma^{-1}(0,\,0) = k.\Sigma^{-1} & \textrm{where $\Sigma^{-1}$ is an element of degree $-1$,} \\
\Sigma^{-1}(m,\, n)=0 & \textrm{otherwise.}
\end{array}\right.$$

\begin{dei}[Desuspension $\Sigma^{-1} V$]
The \emph{desuspension} of the dg-$\Sy$-bimodule $V$ is the dg-$\Sy$-bimodule
$\Sigma^{-1} V := \Sigma^{-1} \otimes V$.
\end{dei}

\subsubsection{\bf Partial composition product}

Let $(\Po,\, \mu,\, \eta,\, \varepsilon)$ be an augmented properad. Since
$2$-level connected graphs are endowed with a bigraduation given by the number
of vertices on each level, we can decompose the composition product $\mu$ on connected
graphs in the
following way $\mu=\bigoplus_{r,\, s \, \in \mathbb{N}} \mu_{(r,\,s)}$, where
$$\mu_{(r,\,s)} \ : \ (I\oplus \underbrace{\oPo}_{r}) \boxtimes_c
 (I \oplus \underbrace{\oPo}_{s}) \to \oPo.$$
The product $\mu_{(r,\, s)}$ gives the composition of $r$ non-trivial operations
with $s$ non-trivial operations.

\begin{dei}[Partial composition product]
We call \emph{partial composition product} the restriction $\mu_{(1,\, 1)}$
of the composition product $\mu$ to $ (I\oplus
\underbrace{\oPo}_{1})\boxtimes_c (I\oplus
\underbrace{\oPo}_{1})$.
\end{dei}

Notice that the partial composition product is the composition product
of two non-trivial elements of $\Po$ where at least one output of the first one
is linked to one input of the second one.

\begin{rem}
When $\Po$ is an operad, this partial composition product is
exactly the partial product denoted $\circ_i$ by V. Ginzburg and
M.M. Kapranov in \cite{GK}. In the framework of
$\frac{1}{2}$-PROPs introduced by M. Markl and A.A. Voronov in
\cite{MV}, these authors only consider partial composition
products where the upper element has only one output or where the
lower element has only one input. The partial composition product
defined by W.L. Gan in \cite{Gan} on dioperads only involves
graphs where the two elements of $\oPo$ are linked by a unique
edge. In these cases, the related monoidal product is generated by
the partial composition product (every composition can be written
with a finite number of partial products). One can see that the
composition product of a properad has the same property.
\end{rem}

\begin{lem}
\label{thetainduitproduit}
Let $(\Po,\, \mu,\, \eta,\, \varepsilon)$ be a augmented dg-properad. The partial
composition product $\mu_{(1,\, 1)}$ induces an homogenous morphism $\theta$ of degree $-1$
$$\theta \ : \ \F_{(2)}^c(\Sigma \oPo)=(I\oplus
\underbrace{\Sigma \oPo}_{1})\boxtimes_c (I\oplus
\underbrace{\Sigma \oPo}_{1}) \to \Sigma \oPo.$$
\end{lem}

Recall that when $\Po$ is a dg-$\Sy$-bimodule, the free properad
$\F(\Sigma \oPo)$ on $\Sigma \oPo$ is bigraded. The first graduation comes
from the number of operations of $\Sigma \oPo$ used to represent an element of
$\F(\Sigma \oPo)$. We denote this graduation $\F_{(n)}(\Sigma \oPo)$.
And the second graduation , the homological degree, is given
by the sum of the homological degrees of the operations of $\Sigma \oPo$.

\begin{proo}
Let $(1,\ldots,\, 1,\, \Sigma q,\, 1,\ldots ,\,
1)\sigma(1,\ldots,\, 1,\, \Sigma p,\, 1,\ldots ,\, 1)$
represent a homogenous element of homological degree $|q|+|p|+2$ of
$\F_{(2)}^c(\Sigma \oPo)$. We define $\theta$ by the following formula
\begin{eqnarray*}
&& \theta\left((1,\ldots,\, 1,\, \Sigma q,\, 1,\ldots ,\,
1)\sigma(1,\ldots,\, 1,\, \Sigma p,\, 1,\ldots ,\, 1)\right):=\\
&& (-1)^{|q|}\Sigma\mu_{(1,\, 1)}(1,\ldots,\, 1,\,  q,\, 1,\ldots
,\, 1)\sigma(1,\ldots,\, 1,\, p,\, 1,\ldots ,\, 1).
\end{eqnarray*}
Since $\Po$ is a differential properad, this last element has degree
$|q|+|p|+1$.
\end{proo}

With Proposition~\ref{coderivation}, we can associate a coderivation
$d_\theta \, : \,\F^c(\Sigma \oPo) \to \F^c(\Sigma \oPo) $ to the homogenous morphism
any morphism $\theta$.

\begin{pro}
\label{propriétés-codérivation}
The coderivation $d_\theta$ has the following properties :
\begin{enumerate}
\item The equation $\delta d_\theta +d_\theta \delta=0$  is true if and only if
the partial composition product $\mu_{(1,\, 1)}$ is a morphism of dg-$\Sy$-bimodules.
\item We have ${d_\theta}^2=0$.
\end{enumerate}
\end{pro}

\begin{proo}
$\ $
\begin{enumerate}
\item Following Proposition~\ref{coderivation}, applying $d_\theta$ to an element
of $\F^c(\Sigma\oPo)$ represented by a graph $g$ corresponds
to apply $\theta$ to every possible sub-graph of $g$. Here, we apply $\theta$ to every
couple of vertices linked by at least on edge and such that there exists no vertex between.
For a couple of such vertices indexed by $\Sigma q$ and $\Sigma p$, we have
\begin{eqnarray*}
&& \delta \circ d_\theta \left( (1,\ldots,\, 1,\, \Sigma q,\,
1,\ldots ,\, 1)\sigma(1,\ldots,\, 1,\, \Sigma p,\, 1,\ldots ,\, 1)
\right)  = \\
&&(-1)^{|q|}\delta \left( \Sigma \mu_{(1,\ 1)}((1,\ldots,\, 1,\,
q,\, 1,\ldots ,\, 1)\sigma(1,\ldots,\, 1,\, p,\, 1,\ldots ,\, 1))
\right)  = \\
&& (-1)^{|q|+1}\Sigma \, \delta \left( \mu_{(1,\ 1)}((1,\ldots,\,
1,\, q,\, 1,\ldots ,\, 1)\sigma(1,\ldots,\, 1,\, p,\, 1,\ldots ,\,
1)) \right)
\end{eqnarray*}
and
\begin{eqnarray*}
&& d_\theta \circ \delta \left( (1,\ldots,\, 1,\, \Sigma q,\,
1,\ldots ,\, 1)\sigma(1,\ldots,\, 1,\, \Sigma p,\, 1,\ldots ,\, 1)
\right)  = \\
&& d_\theta \big( -(1,\ldots,\, 1,\, \Sigma \delta(q),\, 1,\ldots
,\, 1)\sigma(1,\ldots,\, 1,\, \Sigma p,\, 1,\ldots ,\, 1) +
 \\
 &&(-1)^{|q|} (1,\ldots,\, 1,\, \Sigma q,\, 1,\ldots ,\,
1)\sigma(1,\ldots,\, 1,\, \Sigma \delta(p),\, 1,\ldots ,\, 1)
\big)= \\
&& (-1)^{|q|}\Sigma \mu_{(1,\, 1)}((1,\ldots,\, 1,\, \delta(q),\,
1,\ldots ,\, 1)\sigma(1,\ldots,\, 1,\, p,\, 1,\ldots ,\, 1))+ \\
&& \Sigma \mu_{(1,\, 1)} ((1,\ldots,\, 1,\, q,\, 1,\ldots ,\,
1)\sigma(1,\ldots,\, 1,\, \delta(p),\, 1,\ldots ,\, 1)).
\end{eqnarray*}
Therefore, $\delta d_\theta +d_\theta \delta =0$ implies that
\begin{eqnarray*}
&&\delta \left( \mu_{(1,\, 1)} ((1,\ldots,\, 1,\,  q,\, 1,\ldots
,\, 1)\sigma(1,\ldots,\, 1,\,  p,\, 1,\ldots ,\, 1))   \right) = \\
&&\mu_{(1,\, 1)}( (1,\ldots,\, 1,\,  \delta(q),\, 1,\ldots ,\,
1)\sigma(1,\ldots,\, 1,\, p,\, 1,\ldots ,\, 1)) + \\
&& (-1)^{|q|}\mu_{(1,\, 1)}( (1,\ldots,\, 1,\,  q,\, 1,\ldots ,\,
1)\sigma(1,\ldots,\, 1,\, \delta(p),\, 1,\ldots ,\, 1)),
\end{eqnarray*}
which means that the partial composition product $\mu_{(1,\,1)}$
is a morphism of dg-$\Sy$-bimodule. In the other way, applying
$\delta d_\theta +d_\theta \delta$ corresponds to do a sum of
expressions of the previous type. \\

\item To show that $d_\theta^2=0$, we have to consider the pairs of
couples of vertices of a graph $g$. There are two possible cases :
the pairs of couples of vertices such that the four vertices are different (a),
the pairs of couples of vertices with one common vertex (b).\\

\begin{enumerate}
\item To an element of the form
$$X\otimes \Sigma q_1 \otimes
\Sigma p_1 \otimes Y \otimes \Sigma q_2 \otimes \Sigma p_2 \otimes
Z, $$
we apply $\theta$ twice beginning with one different couple each time. It gives
\begin{eqnarray*}
&&\big( (-1)^{|X|+|q_1|}(-1)^{|X|+|q_1|+|p_1|+1+|Y|+|q_2|} +
(-1)^{|X|+|q_1|+|p_1|+|Y|+|q_2|} (-1)^{|X|+|q_1|} \big) \\
&&X\otimes \Sigma \mu_{(1,\, 1)}(q_1\otimes p_1) \otimes Y \otimes
\Sigma \mu_{(1,\, 1)}(q_2\otimes p_2)\otimes Z =0.
\end{eqnarray*}

\item In this case, there are two possible configurations
\begin{itemize}
\item The common vertex is between the two other vertices (following
the flow of the graph) (\emph{cf.} Figure~\ref{entresommets} ).

\begin{figure}[h]
$$ \xymatrix{ & & \\
 & *+[F-,]{\Sigma p}  \ar@{.}[u] \ar@{.}[ul]\ar@{.}[ur] \ar@/_/@{-}[d]\ar@/^/@{-}[d]& \\
 & *+[F-,]{\Sigma q} \ar@{-}[d]& \\
 & *+[F-,]{\Sigma r} \ar@{.}[ul]\ar@{.}[ur]\ar@{.}[dl]\ar@{.}[dr]& \\
 & &} $$ \caption{}  \label{entresommets}
\end{figure}

Therefore, to an element of the form
$$X\otimes \Sigma r \otimes \Sigma q \otimes \Sigma p \otimes Y,$$
if we apply $\theta$ twice beginning by a different couple each time. We have
\begin{eqnarray*}
&& (-1)^{|X|+|r|}(-1)^{|X|+|r|+|q|} X\otimes \Sigma \mu_{(1,\,
1)}(\mu_{(1,\, 1)}(r\otimes q)\otimes p)\otimes Y + \\
&& (-1)^{|X|+|r|+1+|q|}(-1)^{|X|+|r|} X\otimes \Sigma
\mu_{(1,\,1)} (r\otimes \mu_{(1,\,1)}(q\otimes p))\otimes Y=0,
\end{eqnarray*}
by associativity of the partial product $\mu_{(1,\,1)}$ (which comes from
the associativity of $\mu$).
\item Otherwise, the common vertex is above (\emph{cf.}
Figure~\ref{audessussommets}) or under the two others.

\begin{figure}[h]
$$ \xymatrix{& & & \\
& &*+[F-,]{\Sigma p} \ar@{.}[u] \ar@{.}[ul] \ar@{.}[ur] \ar@{.}[d] & \\
& *+[F-,]{\Sigma r}
\ar@{-}[ur]\ar@{.}[u]\ar@{.}[ul]\ar@{.}[d]\ar@{.}[dl]  & *=0{}
\ar@{.}[d]&
*+[F-,]{\Sigma q} \ar@{-}@/_/[ul]\ar@{-}@/^/[ul]\ar@{.}[u] \ar@{.}[d] \\
 & & & } $$ \caption{}  \label{audessussommets}
\end{figure}
We same kind of calculus gives this time
\begin{eqnarray*}
&&(-1)^{|X|+|r|+1+|q|}(-1)^{|X|+|r|}X\otimes \Sigma
\mu_{(1,\,1)}(r\otimes
\mu_{(1,\,1)}(q\otimes p))\otimes Y+ \\
&&(-1)^{(|r|)+1(|q|+1)+|X|+|q|+1+|r|}(-1)^{|X|+|q|}X\otimes
\Sigma\mu_{(1,\, 1)}(q\otimes \mu_{(1,\, 1)} (r\otimes p))\otimes
Y=0.
\end{eqnarray*}
And the associativity of $\mu_{(1,\, 1)}$ gives
$$\mu_{(1,\, 1)}(q\otimes \mu_{(1,\, 1)} (r\otimes p))= (-1)^{|r||q|}
\mu_{(1,\,1)}(r\otimes \mu_{(1,\,1)}(q\otimes p)). $$
\end{itemize}
\end{enumerate}
We conclude by remarking that the image under the morphism $d_\theta^2$ is a sum
of expressions of such type.
\end{enumerate}
\end{proo}

\begin{rem}
This proposition lies on the natural Koszul-Quillen sign rules.
In the article of \cite{GK}, these
signs appear via the operad $Det$.
\end{rem}

\begin{dei}[Pair of adjacent vertices]

Any pair of vertices of a graph linked by at least on edge and without a vertex between
is called a \emph{pair of adjacent vertices}. It corresponds to sub-graphs of the form
$\F_{(2)}(V)=\F^c_{(2)}(V)$ and are the composable pairs under
the coderivation $d_\theta$.
\end{dei}

\begin{rem}
The image under the coderivation $d_\theta$ is given by the composition of all
pairs of adjacent vertices. In order to define the homology of graphs, M. Kontsevich
has introduced in \cite{Ko} the notion of \emph{edge contraction}.
This notion corresponds to the composition of a pair of vertices linked by only one edge.
In the case of operads (trees), $\frac{1}{2}$-PROPs and dioperads (graphs of genus $0$),
adjacent vertices are linked by only one edge. Therefore, this notion applies there.
But in the framework of properads, this notion is too restrictive and its natural generalization
is given by the coderivation $d_\theta$.
\end{rem}

\subsubsection{\bf Partial coproduct}

We dualize here the results of the previous section.

\begin{dei}[Partial coproduct]

Let $(\mathcal{C},\, \Delta,\, \varepsilon,\, \eta)$ be a coaugmented coproperad.
The following composition is called \emph{partial coproduct}
$$\overline{\mathcal{C}}  \xrightarrow{\Delta } \mathcal{C}\boxtimes_c
 \mathcal{C} \twoheadrightarrow
(I\oplus \underbrace{\oC}_{1}) \boxtimes_c
(I\oplus \underbrace{\oC}_{1}).$$
It corresponds to the part of the image of $\Delta$ represented by graphs
with two vertices. We denote it $\Delta_{(1,\,1)}$.
\end{dei}

\begin{lem}
\label{thetainduitcoproduit}
Let $(\mathcal{C},\, \Delta ,\,\varepsilon,\, \eta)$ be
a differential coaugmented coproperad. The partial coproduct induces
an homogenous morphism of degree $-1$
$$\theta' \ : \ \Sigma^{-1} \oC \to \F_{(2)}(\Sigma^{-1} \oC)=
(I\oplus \underbrace{\oC}_1) \boxtimes_c
(I\oplus \underbrace{\oC}_1).$$
\end{lem}

\begin{proo}
Let $c$ be an element of $\oC$. Using the same kind of notations as Sweedler
for Hopf algebras, we have
$$\Delta_{(1,\,1)}(c)=\sum_{(c',\, c'')}
(1,\ldots,\, 1,\,  c',\, 1,\ldots ,\, 1)\sigma(1,\ldots,\, 1,\,
 c'',\, 1,\ldots ,\, 1).$$
Hence, we define  $\theta'(\Sigma^{-1}c)$ by the formula
$$\theta'(\Sigma^{-1}c):=-\sum_{(c',\, c'')} (-1)^{|c'|}
(1,\ldots,\, 1,\,  \Sigma^{-1}\, c',\, 1,\ldots ,\,
1)\sigma(1,\ldots,\, 1,\, \Sigma^{-1}\,c'',\, 1,\ldots ,\, 1).$$
\end{proo}

\begin{rem}
The sign $-$ appearing in the definition of $\theta'$ is not essential here. His role
will be obvious in the proof of the bar-cobar resolution (\emph{cf.}
Theorem~\ref{barcobarresolution}).
\end{rem}

With Proposition~\ref{derivation}, we can associate a derivation
$d_{\theta'} \, :\, \F(\Sigma^{-1}\oC)\to \F(\Sigma^{-1} \oC)$ of degree $-1$
to the homogenous morphism $\theta'$ which verifies the following properties.

\begin{pro}
$ \ $
\begin{enumerate}
\item The equation $\delta d_{\theta'} +d_{\theta'} \delta =0$ is
true if and only if the partial coproduct $\Delta_{(1,\, 1)}$ is a
morphism of dg-$\Sy$-bimodules. \item We have ${d_{\theta'}}^2=0$.
\end{enumerate}
\end{pro}

\begin{proo}
The proof of Proposition~\ref{derivation} shows that the image of
an element of $\F(\Sigma^{-1} \oC)$ represented by a graph $g$
under the derivation $d_{\theta'}$ is a sum where the application
$\theta'$ is applied on each vertex of $g$. Therefore, the
arguments to show this proposition are the same as for the
previous proposition. For instance, the second point comes from
the coassociativity of $\Delta_{(1,\,1)}$, which is induced by the
one of $\Delta$, and from the rules of signs.
\end{proo}

\begin{rem}
The derivation $d_{\theta'}$ corresponds to take the partial
coproduct of every element indexing a vertex of a graph.
Therefore, this derivation is the natural generalization of the
notion of \emph{vertex expansion} introduced by M. Kontsevich in
\cite{Ko} in the framework of graph cohomology. In the case of
operads, the derivation $d_{\theta'}$ is the same morphism as the
one defined by the ``vertex expansion" of trees.
\end{rem}

\subsection{Definitions of the bar and cobar constructions}
\label{Defsbarcobar}

With the propositions given above, we can now define the bar and cobar constructions for
properads.

\subsubsection{\bf Reduced bar and cobar constructions}

\begin{dei}[Reduced bar construction]
To any differential augmented properad $(\Po,\, \mu,\,
\eta,\, \varepsilon)$, we can associate the quasi-cofree coproperad
$\bar{\B}(\Po) := \F^c(\Sigma \oPo)$ with the differential defined by
$\delta_\theta=\delta +d_\theta$, where $\theta$ is the morphism
induced by the partial product of $\Po$ (\emph{cf.} Lemma~\ref{thetainduitproduit}).
This quasi-cofree coproperad is called the \emph{reduced bar construction of $\Po$}.
\end{dei}

We denote $\bar{\B}_{(s)}(\Po) := \F^c_{(s)}(\Sigma \oPo)$. Therefore, $d_\theta$ defines
the following chain complex
$$\xymatrix{\cdots \ar[r]^(0.4){d_\theta} & \bar{\B}_{(s)}(\Po)
\ar[r]^(0.5){d_\theta} & \bar{\B}_{(s-1)}(\Po)
\ar[r]^(0.6){d_\theta} & \cdots \ar[r]^(0.4){d_\theta} &
\bar{\B}_{(1)}(\Po) \ar[r]^(0.5){d_\theta}&
\bar{\B}_{(0)}(\Po).}$$

Proposition~\ref{propriétés-codérivation} gives that $\delta d_\theta+d_\theta \delta =0$ and
that ${d_\theta}^2=0$. Hence, we have ${\delta_\theta}^2=0$
and we see that the reduced bar construction is the total complex of a bicomplex.

Dually, we have the following definition :

\begin{dei}[Reduced cobar construction]
To any differential coaugmented coproperad $(\mathcal{C}$, $\Delta$, $\varepsilon$, $\eta)$
we can associate the quasi-free properad  $\bar{\B}^c(\mathcal{C}) :=
\F(\Sigma^{-1} \oC)$ with the differential defined by $\delta_{\theta'}=\delta +d_{\theta'}$,
where $\theta'$ is the morphism induced by the partial coproduct of $\Co$
(\emph{cf.} Lemma~\ref{thetainduitcoproduit}). This quasi-free properad is called the
\emph{reduced cobar construction of $\Co$}.
\end{dei}

Once again, the derivation $d_{\theta'}$ defines a cochain complex
$$\xymatrix{\bar{\B}^c_{(0)}(\mathcal{C}) \ar[r]^{d_{\theta'}} &
\bar{\B}^c_{(1)}(\mathcal{C}) \ar[r]^(0.6){d_{\theta'}} & \cdots
\ar[r]^(0.45){d_{\theta'}}&
\bar{\B}^c_{(s)}(\mathcal{C})\ar[r]^(0.45){d_{\theta'}}&
\bar{\B}^c_{(s+1)}(\mathcal{C}) \ar[r]^(0.6){d_{\theta'}}& \cdots
.}$$

The reduced cobar construction is the total complex of a bicomplex.

\begin{rem}
When $\Po$ (or $\Co$) is an algebra or an operad (a coalgebra or a cooperad), we find
the classical definitions of bar and cobar construction
(\emph{cf.}
S\'eminaire H. Cartan \cite{Cartan} and \cite{GK}).
\end{rem}

\subsubsection{\bf Bar construction with coefficients}
\label{barconstructrioncoef}

Let $(\Po, \, \mu ,\, \eta,\, \varepsilon)$ be an augmented differential properad.
On $\bar{\B}(\Po)$ we can define two homogenous morphisms $\theta_r\, :\,
\bar{\B}(\Po)\to\bar{\B}(\Po)\boxtimes_c \Po $ and $\theta_l\, :\,
\bar{\B}(\Po)\to \Po \boxtimes_c \bar{\B}(\Po)$ of degree $-1$
by
$$\theta_r \ : \ \bar{\B}(\Po)=\F^c(\Sigma
\oPo)\xrightarrow{\Delta}  \F^c(\Sigma \oPo)\boxtimes_c
\F^c(\Sigma \oPo) \twoheadrightarrow \F^c(\Sigma \oPo)\boxtimes_c
(I \oplus \underbrace{\oPo}_1),$$ and by
$$\theta_l \ : \ \bar{\B}(\Po)=\F^c(\Sigma
\oPo)\xrightarrow{\Delta}  \F^c(\Sigma \oPo)\boxtimes_c
\F^c(\Sigma \oPo) \twoheadrightarrow (I \oplus
\underbrace{\oPo}_1)\boxtimes_c \F^c(\Sigma \oPo).$$

Remark that these two morphisms correspond to extract one non-trivial operation from the top or the
bottom of graphs that represent elements of $\F^c(\Sigma \oPo)$.

\begin{lem}
\label{lemme-bar-construction-augmentée}
The morphism $\theta_r$ induces an homogenous morphism $d_{\theta_r}$
of degree $-1$
$$d_{\theta_r}\ : \ \bar{\B}(\Po)\boxtimes_c \Po \to
\bar{\B}(\Po)\boxtimes_c \Po .$$
And the dg-$\Sy$-bimodule $\bar{\B}(\Po)\boxtimes_c \Po$ endowed with the differential
defined by the sum of the canonical differential $\delta$ with the coderivation $d_\theta$
of $\bar{\B}(\Po)$ and the morphism $d_{\theta_r}$ is an analytic quasi-free
$\Po$-module (on the right).
\end{lem}

\begin{proo}
Since $(\delta+d_\theta)^2=0$,
the equation $(\delta +d_\theta +d_{\theta_r})^2=0$ is equivalent to
$$\left\{  \begin{array}{l}
(\delta+d_\theta)d_{\theta_r}+d_{\theta_r}(\delta+d_\theta)=0 \quad \textrm{and}\\
{d_{\theta_r}}^2=0,
\end{array} \right.$$
These two equations are proved in the same arguments as in the proof of
Proposition~\ref{propriétés-codérivation}.
The functor $\bar{\B}(\Po)$ is analytic in $\oPo$, therefore the
quasi-free $\Po$-module $\bar{\B}(\Po)\boxtimes_c \Po$ is analytic.
\end{proo}

The morphism $d_{\theta_r}$ corresponds to take one by one
operations indexing top vertices
of the graph of an element of $\bar{\B}(\Po)$
and compose them with operations of $\Po$ that are on the first level of
$\bar{\B}(\Po) \boxtimes_c\Po$.

\begin{dei}[Augmented bar construction]
The analytic quasi-free $\Po$-module $\bar{\B}(\Po)\boxtimes_c\Po$
is called the \emph{augmented bar construction (on the right)}.
\end{dei}

One has the same lemma with $\theta_l$ and
$d_{\theta_l}$ and define the \emph{augmented bar construction on the left}
$\Po \boxtimes_c \bar{\B}(\Po)$.\\

In the case of algebras and operads, the augmented bar construction is acyclic.
In the following part of the text, we generalize this result to the augmented
bar construction of properads.\\

Denote $\B(\Po,\, \Po,\, \Po):=\Po\boxtimes_c \bar{\B}(\Po)
\boxtimes_c \Po$. The morphisms $\theta_r$ and $\theta_l$ induce two homogenous morphisms
$$d_{\theta_r} ,\, d_{\theta_l}\ : \\B(\Po,\, \Po,\, \Po) \to \B(\Po,\, \Po,\, \Po), $$
of degree -1. We define the differential $d$ on  $\B(\Po,\, \Po,\, \Po)$ by the sum of
the following morphisms :

\vspace{6pt}

\begin{tabular}{l}
$\bullet \quad$ the canonical differential $\delta$ induced by the one of $\Po$, \\
$\bullet\quad$ the coderivation $d_\theta$ defined on $\bar{\B}(\Po)$,\\
$\bullet\quad$ the homogenous morphism $d_{\theta_r}$ of degree $-1$, \\
$\bullet\quad$ the homogenous morphism $d_{\theta_r}$ of degree $-1$.
\end{tabular}

\begin{lem}
The morphism $d$ verifies the equation $d^2=0$.
\end{lem}

\begin{proo}
Once again, the proof lies on rules of sign.
\end{proo}

\begin{dei}[Bar construction with coefficients]
Let $(\Po,\, \mu,\ \eta,\, \varepsilon)$ be augmented dg-properad,
$L$ be a differential right $\Po$-module and $R$
be a differential left $\Po$-module.
The \emph{bar construction of $\Po$ with coefficients in $L$ and $R$} is the
dg-$\Sy$-bimodule $\B(L,\, \Po,\,
R):= L{\boxtimes_c}_\Po \B(\Po,\, \Po,\, \Po) {\boxtimes_c}_\Po R$, where
the differential $d$ is induced by the one of $\B(\Po,\, \Po,\, \Po)$, defined
previously, and by the ones of $L$ and $R$ denoted $\delta_R$ and $\delta_L$.
\end{dei}

\begin{pro}
The bar construction $\B(L,\, \Po,\, R)$ with coefficients in $L$ and $R$ is
isomorphic, as dg-$\Sy$-bimodule, to $L\boxtimes_c \bar{\B}(\Po) \boxtimes_c R$,
where the differential $d$ is defined by the sum of the following morphisms :

\vspace{6pt}

\begin{tabular}{l}
$\bullet \quad$  the canonical differential $\delta$
induced by the one of $\Po$ on $\bar{\B}(\Po)$,\\
$\bullet \quad$ the canonical differential $\delta_L$ induced by the one of $L$,\\
$\bullet \quad$ the canonical differential $\delta_R$ induced by the one of $R$,\\
$\bullet\quad$ the coderivation $d_\theta$ defined on $\bar{\B}(\Po)$,\\
$\bullet\quad$ a homogenous morphism $d_{\theta_R}$ of degree $-1$,\\
$\bullet\quad$ a homogenous morphism $d_{\theta_L}$ of degree $-1$.
\end{tabular}
\end{pro}

\begin{proo}
One has to understand the effect of the morphism $d_{\theta_r}$ on the relative
composition product. The morphism $d_{\theta_r}$ corresponds to the morphism
$d_{\theta_R}$ defined by the following composition
$$\xymatrix@C=60pt{ L\boxtimes_c \bar{\B}(\Po) \boxtimes_c R
\ar[r]^{\widetilde{\theta_r}}& L \boxtimes_c \bar{\B}(\Po)
\boxtimes_c \Po \boxtimes_c R \ar[r]^(0.55){L\boxtimes_c
\bar{\B}(\Po)\boxtimes_c r} & L\boxtimes_c \bar{\B}(\Po)
\boxtimes_c R,}$$
where the morphism $\widetilde{\theta_r}$ is induced by $\theta_r$.
\end{proo}

From this description of the bar construction with coefficients, we have immediately
the following corollary.

\begin{cor}
$\ $
\begin{itemize}
\item The reduced bar construction is isomorphic to the bar construction with
scalar coefficients $L=R=I$. We have $\bar{\B}(\Po)=\B(I,\, \Po,\, I)$.

\item The chain complex $\B(L,\, \Po,\, \Po)=L\boxtimes_c
\bar{\B}(\Po)\boxtimes_c \Po$ (resp. $\B(\Po,\, \Po,\, R)=\Po
\boxtimes_c \bar{\B} (\Po) \boxtimes_c R$) is an analytic quasi-free
$\Po$-module on the right (resp. on the left).
\end{itemize}
\end{cor}

In the same way as Proposition~\ref{freefuncteurexact} and Theorem~\ref{lemmecomparaisonproperades},
we have that the bar construction is an exact functor.

\begin{pro}
\label{BarFoncteurExact}
The bar construction $\B(L,\, \Po,\, R)$ preserves quasi-isomorphisms.
\end{pro}

\begin{proo}
Let $\beta \, : \, \Po \to \Po'$ be a quasi-isomorphism of augmented dg-properads. Let
$\alpha \, : \, L \to L'$ be a quasi-isomorphism of right dg-$\Po$,$\Po'$-modules and let
$\gamma \, : \, R \to R'$ be a quasi-isomorphism of left dg-$\Po$,$\Po'$-modules. Denote by $\Phi$ the
induced morphism $\B(L,\, \Po,\, R) \to \B(L',\, \Po',\, R')$. To show that $\Phi$ is a
quasi-isomorphism, we introduce the following filtration $F_r:=\oplus_{s\leq r} \B_{(s)}(L,\, \Po,\, R)$, where
$s$ is the number of vertices of the underlying graph indexed by elements of $\Sigma \oPo$.
One can see that this filtration is stable under the differential $d$ of the bar construction. The canonical
differentials $\delta$, $\delta_L$ and $\delta_R$ go from $\B_{(s)}(L,\, \Po,\, R)$ to $\B_{(s)}(L,\, \Po,\, R)$
and the maps $d_\theta$, $d_{\theta_R}$ and $d_{\theta_R}$ go from $\B_{(s)}(L,\, \Po,\, R)$ to $\B_{(s-1)}(L,\, \Po,\, R)$.
Therefore, this filtration induces a spectral sequence $E^*_{s,\, t}$ such that $E^0_{s,\, t}=\B_{(s)}(L,\, \Po,\, R)_{s+t}$
and such that $d^0=\delta+\delta_L+\delta_R$. By the same ideas as in the proof of Proposition~\ref{freefuncteurexact},
we have that $E^0_{s,\, t}(\Phi)$ is a quasi-isomorphism. Since the filtration $F_r$ is bounded below and exhaustive, the
result is given by the classical convergence theorem of spectral sequences
(\emph{cf.} J.A. Weibel \cite{Weibel} $5.5.1$).

\end{proo}

\subsubsection{\bf Cobar construction with coefficients}

We dualize the previous arguments to define the cobar construction with coefficients. \\

Let $(\mathcal{C},\,  \Delta,\ \varepsilon,\, \eta)$ be
a differential coaugmented coproperad and let $(L,\ l)$
and $(R,\, r)$ be two differential $\mathcal{C}$-comodules on the right
and on the left.

The two $\Co$-comodules induce the following morphisms :
$$\theta'_R\ : \ R \xrightarrow{r} \mathcal{C}\boxtimes_c R
\twoheadrightarrow (I\oplus \underbrace{\oC}_1)\boxtimes_c R,$$
$$\theta'_L \ : \ L \xrightarrow{r}  L\boxtimes_c \mathcal{C}
\twoheadrightarrow L\boxtimes_c (I\oplus \underbrace{\oC}_1).$$

\begin{lem}
The applications $\theta'_R$ and $\theta'_L$ both induce an homogenous
morphism of degree $-1$ on $L\boxtimes_c \bar{\B}^c(\mathcal{C})\boxtimes_c R$
of the form
\begin{eqnarray*}
d_{\theta'_R} \ &:& \ \xymatrix{ \ L\boxtimes_c
\bar{\B}^c(\mathcal{C})\boxtimes_c
R\ar[r]^(0.4){\widetilde{\theta'_R}} & \ {\begin{array}{l} \\ L
\boxtimes_c \bar{\B}^c(\mathcal{C}) \boxtimes_c (I\oplus
\underbrace{\oC}_1) \boxtimes_c R \end{array}} \ar[r]& }\\
&&\xymatrix@C=80pt{{\begin{array}{l} \\ L \boxtimes_c
\bar{\B}^c(\mathcal{C}) \boxtimes_c (I\oplus
\underbrace{\Sigma^{-1} \oC}_1) \boxtimes_c R \end{array}}
\ar[r]^(0.6){L \boxtimes_c \mu_{\bar{\B}^c(\mathcal{C})}
\boxtimes_c R} & L\boxtimes_c \bar{\B}^c(\mathcal{C})\boxtimes_c
R},
\end{eqnarray*}

\begin{eqnarray*}
\theta'_L \ &:& \ \xymatrix{d_{\theta'_L} \ : \ L\boxtimes_c
\bar{\B}^c(\mathcal{C})\boxtimes_c
R\ar[r]^(0.4){\widetilde{\theta'_L}} & {\begin{array}{l} \\ L
\boxtimes_c (I\oplus \underbrace{\oC}_1) \boxtimes_c
\bar{\B}^c(\mathcal{C}) \boxtimes_c R \end{array}}
\ar[r]& }\\
&&\xymatrix@C=80pt{{\begin{array}{l} \\ L \boxtimes_c (I\oplus
\underbrace{\Sigma^{-1} \oC}_1) \boxtimes_c
\bar{\B}^c(\mathcal{C}) \boxtimes_c R \end{array}}\ar[r]^(0.6){L
\boxtimes_c \mu_{\bar{\B}^c(\mathcal{C})} \boxtimes_c R} &
L\boxtimes_c \bar{\B}^c(\mathcal{C})\boxtimes_c R}.
\end{eqnarray*}
\end{lem}

\begin{pro}
The differential $d$ on the $\Sy$-bimodule $L\boxtimes_c
\bar{\B}^c(\mathcal{C})\boxtimes_c R$ is the sum of the following
morphisms  :

\vspace{6pt}

\begin{tabular}{l}
$\bullet \quad$ the canonical differential $\delta$ induced by the one of $\Co$ on $\bar{\B}^c(\mathcal{C})$,\\
$\bullet \quad$ the canonical differential $\delta_L$ induced by the one of $L$,\\
$\bullet \quad$ the canonical differential $\delta_R$ induced by the one of $R$,\\
$\bullet\quad$ the derivation $d_\theta$ defined on $\bar{\B}^c(\mathcal{C})$,\\
$\bullet\quad$ the homogenous morphism $d_{\theta'_R}$ of degree $-1$,\\
$\bullet\quad$ the homogenous morphism $d_{\theta'_L}$ of degree $-1$.
\end{tabular}
\end{pro}

\begin{dei}[Cobar construction with coefficients]

The dg-$\Sy$-bimodule $L\boxtimes_c
\bar{\B}^c(\mathcal{C})\boxtimes_c R$ endowed with the differential
$d$ is called the \emph{cobar construction with coefficients in
$L$ and $R$} and is denoted $\B^c(L,\, \mathcal{C},\ R)$.
\end{dei}

\begin{cor}
$ \ $
\begin{itemize}
\item The reduced cobar construction $\bar{\B}^c(\mathcal{C})$
is isomorphic to the cobar construction with coefficients in the
scalar $\Co$-comodule $I$.

\item The chain complex $\B^c(L,\, \mathcal{C},\,
\mathcal{C})=L\boxtimes_c \bar{\B}^c(\mathcal{C})\boxtimes_c
\mathcal{C}$ (resp. $\B(\mathcal{C},\, \mathcal{C},\,
R)=\mathcal{C} \boxtimes_c \bar{\B}^c (\mathcal{C}) \boxtimes_c
R$) is an analytic quasi-cofree $\mathcal{C}$-comodule on the right
 (resp. on the left)
\end{itemize}
\end{cor}

\subsection{Acyclicity of the augmented bar and cobar constructions}
\label{acyclicityaugm}

When $\Po$ is an algebra $A$, that's-to-say when the
dg-$\Sy$-bimodule $\Po$ is concentrated in $\Po(1,\,1)=A$, one has
that the augmented bar construction is acyclic. To prove it, one
introduces directly a contracting homotopy (\emph{cf.} S\'eminaire
H. Cartan \cite{Cartan}). In the case of operads, B. Fresse has
proved that the augmented bar construction on left $\Po\circ
\bar{\B}(\Po)$ is acyclic using a contracting homotopy. This
homotopy lies on the fact the monoidal product $\circ$ is linear
on the left. The acyclicity of the augmented bar construction on
the right comes from this result and the use of comparison lemmas
for quasi-free modules (\emph{cf.} \cite{Fresse} section $4.6$).

In the framework of properads, since the composition product
$\boxtimes_c$ is not linear on the left nor on the right, we can
not use the same method. We begin by defining a filtration on the
chain complex $\left(\Po \boxtimes_c \bar{\B}(\Po),\, d \right)$.
This filtration induces a convergent spectral sequence. Then, we
introduce a contracting homotopy on the chain complexes $\left(
E^0_{p,\, *},\, d^0 \right)$, for $p>0$, to show that they are
acyclic.

\subsubsection{\bf Acyclicity of the augmented bar construction}

Let $(\Po,\, \mu,\,\eta,\, \varepsilon)$ be a augmented dg-properad.
The dg-$\Sy$-bimodule $\Po\boxtimes_c\bar{\B}(\Po)$ is the image of the functor
$$\oPo \mapsto (I\oplus \oPo)\boxtimes_c\F^c(\Sigma \oPo).$$
Since this functor is a sum a polynomial functors, it is an split analytic functor
(\emph{cf.} \ref{splitanalyticfunctor}). For an element of $(\Po,\, \mu,\,
\eta,\, \varepsilon)$, the graduation is given by the number of
vertices indexed by elements of $\oPo$. Once again, we denote this (total)
graduation by
$$\Po\boxtimes_c\bar{\B}(\Po)=\bigoplus_{s\in\mathbb{N}} \left(
\Po\boxtimes_c\bar{\B}(\Po)\right)_{(s)}.$$

We consider the filtration defined by
$$ F_i=F_i\left(\Po\boxtimes\bar{\B}(\Po)\right):=\bigoplus_{s\leq i} \left(
\Po\boxtimes_c\bar{\B}(\Po) \right)_{(s)}.$$

Therefore, $F_i$ is composed by elements
$\Po\boxtimes\bar{\B}(\Po)$ represented by graphs indexed at most
$i$ elements of $\oPo$.

\begin{lem}
The filtration $F_i$ of the dg-$\Sy$-bimodule
$\Po\boxtimes_c\bar{\B}(\Po)$ is stable under the differential $d$
of the augmented bar construction.
\end{lem}

\begin{proo}
Lemma~\ref{lemme-bar-construction-augmentée} give the form of the
differential $d$ of the augmented bar construction. It is the sum
of the three following terms :
\begin{enumerate}
\item The differential $\delta$ coming from the one of $\Po$. This one does not change the number
of elements of $\oPo$. Hence, one has $\delta(F_i)\subset F_i$.
\item The coderivation $d_\theta$ of the reduced bar construction $\bar{\B}(\Po)$.
This morphism corresponds to the composition of pairs of adjacent vertices. Therefore,
it verifies $d_\theta(F_i)\subset
F_{i-1}$.
\item The morphism $d_{\theta_l}$. Applying this morphism corresponds to
taking an operation of $\Sigma \oPo$ of $\bar{\B}(\Po)$ from the bottom
and compose it with elements of $\Po$ on the last level of the graph,
$$ \Po\boxtimes_c\bar{\B}(\Po) \to  \Po \boxtimes_c(I\oplus \oPo)
\boxtimes_c\bar{\B}(\Po) \to  \Po\boxtimes_c\bar{\B}(\Po).$$
The total number of elements of $\oPo$ is decreasing and one has
$d_{\theta_l}(F_i)\subset F_i$.
\end{enumerate}
\end{proo}

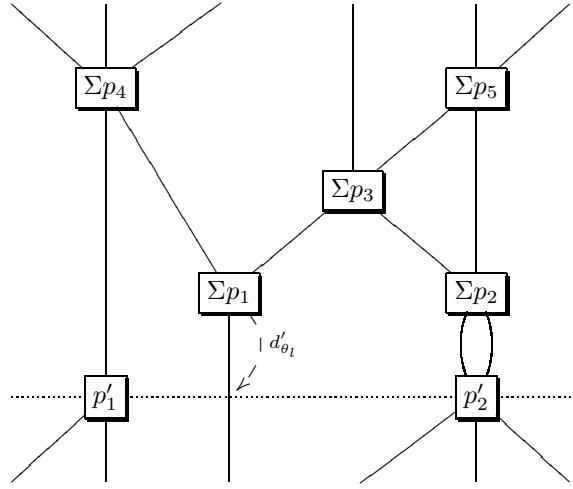
\begin{figure}[h]
$$ \xymatrix{ \ar@{-}[dr]&\ar@{-}[d] & &\ar@{-}[dd] &\ar@{-}[d] &\ar@{-}[dl] \\
 &*+[F-,]{\Sigma p_4} \ar@{-}[ddd] \ar@{-}[ddr] \ar@{-}[ur]& & & *+[F-,]{\Sigma p_5} \ar@{-}[dl] \ar@{-}[dd]&\\
 & & & *+[F-,]{\Sigma p_3}\ar@{-}[dl]\ar@{-}[dr] & &\\
 & & *+[F-,]{\Sigma p_1}\ar@{-}[dd] \ar@/^1pc/@{-->}[d]^{d'_{\theta_l}}& & *+[F-,]{\Sigma p_2} \ar@{-}@/_/[d]\ar@{-}@/^/[d]&\\
 &*+[F-,]{p'_1}\ar@{-}[d]\ar@{-}[dl]  \ar@{.}[l] \ar@{.}[rrr]& & & *+[F-,]{p'_2}\ar@{-}[d]\ar@{-}[dl]\ar@{-}[dr] \ar@{.}[r]&\\
& & & & & } $$ \caption{Form of the morphism $d'_{\theta_l}$}  \label{d'}
\end{figure}

This lemma shows that the filtration $F_i$ induces a spectral
sequence, denote $E^*_{p,\, q}$. The first term is equal to
$$E^0_{p,\, q}=F_p\left( (\Po\boxtimes_c\bar{\B}(\Po))_{p+q}\right)/
F_{p-1}\left( (\Po\boxtimes_c\bar{\B}(\Po))_{p+q}\right),$$ where
$p+q$ is the homological degree. Hence, the module $E^0_{p,\, q}$
is given by the graphs with exactly $p$ vertices indexed by
elements of $\oPo$. We have $E^0_{p,\, q}=\left(
 \left(\Po\boxtimes_c\bar{\B}(\Po)\right)_{(p)}\right)_{p+q}$ and
 the differential $d^0$ is a sum of the morphisms $d^0=\delta
 +d'_{\theta_l}$. The morphism $\delta$ is induced by the one of
 $\Po$ and the morphism $d'_{\theta_l}$ is equal to the morphism
 $d_{\theta_l}$ when it does not change the total number of
 vertices indexed by $\oPo$. Therefore, the morphism
 $d'_{\theta_l}$ takes one operation of $\Sigma \oPo$
in $\bar{\B}(\Po)$ from the bottom and inserts it on the line of
elements of $\Po$ (in the last level) without composing it with
operations of $\oPo$. Otherwise, the morphism $d_{\theta_l}$
implies a non-trivial composition with elements of $\Po$, in this
case $d'_{\theta_l}$ is null. (\emph{cf.} figure~\ref{d'}).\\

The spectral sequence $E^*_{p,\, q}$ is concentrated in the
half-plane $p\ge0$. We compute the homology of the chain complexes
$\left( E^0_{p,\, *},\, d^0 \right)$ to show that the spectral
sequence collapses at rank $E^1_{p,\, q}$.

\begin{lem}
\label{degeneree} At rank $E^1_{p,\, q}$, we have
$$ E^1_{p,\, q}=\left\{ \begin{array}{l}
I \quad \textrm{if} \quad p=q=0, \\
0 \quad \textrm{otherwise}.
\end{array} \right.$$
\end{lem}

\begin{proo}
When $p=0$, we have
$$ E^0_{0,\, q}=\left\{ \begin{array}{l}
I \quad \textrm{if} \quad q=0, \\
0 \quad \textrm{otherwise}.
\end{array} \right.$$
and the differential $d^0$ is null on the modules $E^0_{0,\, q}$.
Hence, the homology of these modules is equal to
$$ E^1_{0,\, q}=\left\{ \begin{array}{l}
I \quad \textrm{if} \quad q=0, \\
0 \quad \textrm{otherwise}.
\end{array} \right.$$\\

When $p>0$, we introduce a contracting homotopy $h$ for each chain
complexes $\left(E^0_{p,*},\ d^0 \right)$.

Thanks to Proposition~\ref{representation}, we can
represent an element of $\Po \boxtimes \bar{\B}(\Po)$ by $(p_1,\,
\ldots ,\, p_r)$ $\sigma$ $(b_1,\, \ldots ,\, b_s)$. When $p_1\in I$,
we define
$$h\left((p_1,\, \ldots ,\, p_r)\sigma (b_1,\, \ldots ,\, b_s)
\right)=0,$$ otherwise $h$ is given by the formula
\begin{eqnarray*}
&& h\left( (p_1,\, \ldots ,\, p_r)\sigma (b_1,\, \ldots ,\,
b_s)\right)=\\
&& (-1)^{(|p_1|+1)(|p_2|+\cdots+|p_r|)} (1,\,\ldots,\
1,\,  p_2,\, \ldots,\, p_r)\sigma'(b',\,b_{i+1},\ \ldots ,\, b_s),
\end{eqnarray*}
where $b'=\mu_{\F(\Sigma \oPo)} \left(\Sigma p_1\otimes
(b_1,\,\ldots,\,b_i)\right)$ with $b_1,\, \ldots,\, b_i$ the
elements of $\bar{\B}(\Po)$ linked via at least one edge to the
vertex indexed by $p_1$ in the graphic representative of $(p_1,\,
\ldots ,\, p_r)\sigma (b_1,\, \ldots ,\, b_s)$. This morphism $h$
is homogenous of degree $+1$ ($\oPo \to \Sigma \oPo$) and does not
change the number of elements of $\oPo$. Therefore, we have
$h\,:\, E^0_{p,\, q} \to E^0_{p,\, q+1}$. Intuitively, the
application $h$ takes the ``first" non-trivial operation between
the ones of $\Po$ on the last level, suspends it, and lifts it
between the ones of $\bar{\B}(\Po)$. Let verify that
$hd^0+d^0h=i_d$.
\begin{enumerate}
\item The rules of signs and suspensions show that $h\delta+\delta
h =0$.
 \item We have $h d'_{\theta_l}+d'_{\theta_l} h=i_d$. We first
 compute $d'_{\theta_l}$ :
\begin{eqnarray*}
&& d'_{\theta_l}\left( (p_1,\, \ldots ,\, p_r)\sigma (b_1,\,
\ldots
,\, b_s)\right)= \\
&&\qquad  \sum
(-1)^{|p'|(|p_{j+1}|+\cdots+|p_r|+|b_1|+\cdots+|b_k|)}\\
&& \qquad (p_1,\ldots,\, p',\, p_{j+1} \ldots,\
p_r)\widetilde{\sigma}(b_1,\ldots,\, b_k,\, b'_{k+1},\ldots ,\,
b_s),
\end{eqnarray*}
where $b_{k+1}=\mu_{\F(\Sigma \oPo)}(\Sigma p'\otimes b'_{k+1})$.
Then, we have
\begin{eqnarray*}
&&hd'_{\theta_l}\left( (p_1,\, \ldots ,\, p_r)\sigma (b_1,\,
\ldots ,\, b_s)
\right) =\\
&&\qquad  \sum (-1)^{|p'|(|p_{j+1}|+\cdots+|p_r|+|b_1|+\cdots+|b_k|)+(|p_1|+1)(|p_2|+\cdots+|p'|+\cdots+|p_r|)} \\
&&\qquad (1,\dots,\, 1,\ p_2,\ldots,\, p',\, p_{j+1} \ldots,\
p_r)\widetilde{\sigma}'(b',\, b_{i+1},\ldots,\, b_k,\,
b'_{k+1},\ldots ,\, b_s).
\end{eqnarray*}
The same calculus for $d'_{\theta_l} h$ gives
\begin{eqnarray*}
&& d'_{\theta_l} h\left( (p_1,\, \ldots ,\, p_r)\sigma (b_1,\,
\ldots
,\, b_s)\right)=\\
&& \qquad (p_1,\, \ldots ,\, p_r)\sigma (b_1,\, \ldots ,\, b_s) - \\
&& \qquad \sum (-1)^{|p'|(|p_{j+1}|+\cdots+|p_r|+|b_1|+\cdots+|b_k|)+(|p_1|+1)(|p_2|+\cdots+|p'|+\cdots+|p_r|)} \\
&&\qquad (1,\dots,\, 1,\ p_2,\ldots,\, p',\, p_{j+1} \ldots,\
p_r)\widetilde{\sigma}'(b',\, b_{i+1},\ldots,\, b_k,\,
b'_{k+1},\ldots ,\, b_s),
\end{eqnarray*}
the sign $-1$ comes from the commutation of $p'$ with $\Sigma
p_1$.
\end{enumerate}
Since $h$ is an homotopy, the chain complexes $\left( E^0_{p,\,
*},\, d^0\right)$ are acyclic for $p>0$, which means that we have
$$  E^1_{p,\, q}=0 \quad \textrm{pour tout} \ q \ \textrm{si} \ p>0.$$
\end{proo}

We can now conclude the proof of the acyclicity of the augmented bar construction.

\begin{thm}
\label{acyclicitébaraugmentée} The homology of the chain complex
 $\left( \Po \boxtimes_c \bar{\B}(\Po),\,
d\right)$ is given by :
$$ \left\{ \begin{array}{l}
H_0\left( \Po \boxtimes_c \bar{\B}(\Po)\right)=I , \\
H_n\left( \Po \boxtimes_c \bar{\B}(\Po)\right)=0 \quad
\textrm{otherwise}.
\end{array} \right.$$
\end{thm}

\begin{proo}
Since the filtration $F_i$ is exhaustive ($\Po \boxtimes_c
\bar{\B}(\Po)=\bigcup_i F_i$) and bounded below ($F_{-1}\left(\Po
\boxtimes_c \bar{\B}(\Po)\right)=0$), we can apply the classical
theorem of convergence of spectral sequences (\emph{cf.}
\cite{Weibel} $5.5.1$). Hence, we get that the spectral sequence
$E^*_{p,\, q}$ converges to the homology of $\Po \boxtimes_c
\bar{\B}(\Po)$
$$E^1_{p,\, q} \Longrightarrow H_{p+q}\left(\Po \boxtimes_c \bar{\B}(\Po),\,
d\right) .$$ And we conclude using Lemma~\ref{degeneree}.
\end{proo}

\begin{cor}
The augmentation morphism
$$\xymatrix@C=50pt{\Po\boxtimes_c \bar{\B}(\Po) \ar[r]^{\varepsilon_\Po
\boxtimes_c \varepsilon_{\F^c(\Sigma \oPo)}} & I\boxtimes_c I = I}
$$ is a quasi-isomorphism.
\end{cor}

This last result can be generalized in the following way :

We define the \emph{augmentation morphism} $\varepsilon(\Po,\,
\Po,\, \Po)$ by the composition
$$\B(\Po,\, \Po,\,
\Po)=\Po\boxtimes_c \bar{\B}(\Po) \boxtimes_c \Po
\xrightarrow{\Po\boxtimes_c \varepsilon_{\F^c(\Sigma
\oPo)}\boxtimes_c \Po}  \Po\boxtimes_c I\boxtimes_c
\Po=\Po\boxtimes_c \Po \twoheadrightarrow \Po \boxtimes_\Po \Po=
\Po.$$

\begin{thm}
Let $R$ be a left differential $\Po$-module. The augmentation
morphism $\varepsilon(\Po,\, \Po,\, R)$
$$ \xymatrix@C=75pt{ \B(\Po,\Po,\,R)=\Po\boxtimes_\Po
\B(\Po,\, \Po,\, \Po) \boxtimes_\Po R \ar[r]^(0.6){\Po
\boxtimes_\Po \varepsilon(\Po,\, \Po,\, \Po) \boxtimes_\Po R }
 & \Po\boxtimes_\Po
 \Po \boxtimes_\Po R = R} $$ is a quasi-isomorphism.
\end{thm}

\begin{proo}
We define the same kind of filtration based on the number of
elements of $\oPo$ indexing the graphs.
\end{proo}

\begin{rem}
Using the same arguments, we have the same results for the
augmented bar construction on the right $\bar{\B}(\Po)\boxtimes_c
\Po$.
\end{rem}

\subsubsection{\bf Acyclicity of the augmented cobar construction}

Dually, we can prove that the augmented cobar construction
$\bar{\B}^c(\mathcal{C})\boxtimes_c \mathcal{C}$ is
acyclic. But in order to make the spectral sequence converge, we
need to work with a weight graded coproperad $\Co$.

\begin{thm}
\label{acyclicitécobarcoaugmentée} For every weight graded
coaugmented coproperad $\Co$, the homology of the augmented cobar
construction $\left(\bar{\B}^c(\mathcal{C})\boxtimes_c \mathcal{C}
,\, d\right)$ is given by
$$ \left\{ \begin{array}{l}
H_0\left( \bar{\B}^c(\mathcal{C})\boxtimes_c \mathcal{C}\right)=I , \\
H_n\left( \bar{\B}^c(\mathcal{C})\boxtimes_c \mathcal{C}\right)=0
\quad \textrm{otherwise}.
\end{array} \right.$$
\end{thm}

\begin{proo}
We define a filtration $F_{i}$ on
$\left(\bar{\B}^c(\mathcal{C})\boxtimes_c \mathcal{C} ,\,
d\right)$ by
$$F_{i}=\bigoplus_{s\ge -i} \left( \bar{\B}^c(\mathcal{C})\boxtimes_c
 \mathcal{C} \right)_{(s)} .$$
This filtration is stable under the differential $d$ :
\begin{enumerate}
\item Since the differential $\delta$ does not change the number
of elements of $\oC$, we have $\delta(F_i)\subset F_i$.

\item The derivation $d_{\theta'}$ of the reduced cobar
construction $\bar{\B}^c(\mathcal{C})$ splits into two parts each
element indexing a vertex. Therefore, the global number of
vertices indexed by elements of $\oC$ is raised by $1$. The
derivation $d_{\theta'}$ verifies $d_{\theta'}(F_i)\subset
F_{i-1}$.

\item The morphism $d_{\theta'_r}$ decomposes elements of $\Co$
into two parts and includes on part in $\bar{\B}^c(\mathcal{C})$.
The global number of elements of $\oC$ indexing vertices is
increasing. We have $d_{\theta'_r}(F_i)\subset F_i$.
\end{enumerate}
This filtration induces a spectral sequence $E^*_{p,\, q}$. Its
first term is given by
$$E^0_{p,\, q}=F_p\left((\bar{\B}^c(\mathcal{C})\boxtimes_c
\mathcal{C})_{p+q}\right)/
F_{p-1}\left((\bar{\B}^c(\mathcal{C})\boxtimes_c
\mathcal{C})_{p+q}\right)=\left((\bar{\B}^c(\mathcal{C})\boxtimes_c
\mathcal{C})_{(-p)} \right)_{p+q}.$$

The differential $d^0$ is the sum of two terms
$d^0=\delta+d'_{\theta'_r}$, where $d'_{\theta'_r}$ corresponds to
$d_{\theta'_r}$ when it does not strictly increase the number of
elements of $\oC$. The morphism $d'_{\theta'_r}$ takes elements of
$\oC$ in the first level, desuspends them and include them in
$\bar{\B}^c(\mathcal{C})$.

We show then that
$$ E^1_{p,\, q}=\left\{ \begin{array}{l}
I \quad \textrm{if} \quad p=q=0, \\
0 \quad \textrm{otherwise}.
\end{array} \right.$$

We define the image under the homotopy $h$ of an element
$(b_1,\ldots ,\,b_r)\, \sigma \,(c_1,\ldots,\, c_s)$ of
$\bar{\B}^c(\mathcal{C})\boxtimes \mathcal{C}$ by first
decomposing $b_1\in \bar{\B}^c(\mathcal{C})_{(n)}$ :
$$\xymatrix{b_1 \ar@{|->}[r]&\sum b'_1\otimes \Sigma^{-1}
c },$$ where $b'_1 \in \bar{\B}^c(\mathcal{C})_{(n-1)}$ and $c\in
\oC$. If the element $\Sigma^{-1}c$ is linked above to elements of
$I$, then we suspend it and include it in the line of elements of
$\Co$. Using the same arguments, one can show that
$hd^0+d^0h=i_d$. \\

Since the coproperad $\Co$ is weight graded, one can decompose the
chain complex $\bar{\B}^c(\mathcal{C})\boxtimes_c \mathcal{C}$
into a direct sum indexed by the total weight :
$$\bar{\B}^c(\mathcal{C})\boxtimes_c \mathcal{C}=
\bigoplus_{\rho \in \mathbb{N}} (\bar{\B}^c(\mathcal{C})\boxtimes
\mathcal{C})^{(\rho)}. $$ This decomposition is compatible with
the filtration $F_i$ and with the homotopy $h$. Therefore, we have
$$ E^1_{p,\, q}\big( (\bar{\B}^c(\mathcal{C})\boxtimes_c \mathcal{C})^{(\rho)}
\big)=0\quad \textrm{for every}\  p,\,q \ \textrm{when}\ \rho>0 \
\textrm{and}$$
$$E^1_{p,\, q}\big( (\bar{\B}^c(\mathcal{C})\boxtimes_c \mathcal{C})^{(0)}
\big)=\left\{ \begin{array}{l}
I \quad \textrm{if} \quad p=q=0, \\
0 \quad \textrm{otherwise}.
\end{array} \right.$$

Since the filtration $F_i$ is bounded on the sub-complex
$(\bar{\B}^c(\mathcal{C})\boxtimes_c \mathcal{C})^{(\rho)}$
$$F_0\big ((\bar{\B}^c(\mathcal{C})\boxtimes_c \mathcal{C})^{(\rho)}\big) =
(\bar{\B}^c(\mathcal{C})\boxtimes_c \mathcal{C})^{(\rho)} \quad
\textrm{et} \quad F_{-\rho-1}\big
((\bar{\B}^c(\mathcal{C})\boxtimes_c \mathcal{C})^{(\rho)}\big)
=0,$$ we can apply the classical theorem of convergence of
spectral sequences (\emph{cf.} \cite{Weibel} $5.5.1$). This gives
that the spectral sequence $E^*_{p,\, q}\big(
(\bar{\B}^c(\mathcal{C})\boxtimes_c \mathcal{C})^{(\rho)} \big)$
converges to the homology of $(\bar{\B}^c(\mathcal{C})\boxtimes_c
 \mathcal{C})^{(\rho)}$. We conclude by taking the direct sum on
 $\rho$.
\end{proo}

\begin{rem}
Using the same methods, we can show that the augmented bar
construction on the left is acyclic.
\end{rem}



\section{Comparison lemmas}
\label{ComparisonLemmas}

We prove here comparison lemmas for quasi-free $\Po$-modules
and quasi-free properads. These lemmas are generalizations to the composition
product $\boxtimes_c$ of classical lemmas for the tensor product $\otimes_k$.
Notice that B. Fresse has given in \cite{Fresse} such lemmas for the
composition product $\circ$ of operads. These technical lemmas allow us to
show that the bar-cobar construction of a weight graded properad $\Po$ is
a resolution of $\Po$.

\subsection{Comparison lemma for quasi-free modules}

\subsubsection{\bf Filtration and spectral sequence of an analytic
quasi-free module}

Let $\Po$ be a weight graded differential properad $\Po$ and let
$L=M\boxtimes_c \Po$ be an analytic quasi-free $\Po$-module on the right, where
$M=\Upsilon(V)$. Denote $M^{(\alpha)}_{d}$ be the sub-$\Sy$-bimodule of $M$
composed by elements of homological degree $d$ and global weight $\alpha$ (\emph{cf.}
\ref{analyticquasifreemodule}).

We define the following filtration on $L$ :
$$F_s(L):=\bigoplus_\Xi \bigoplus_{|\bar{d}|+|\bar{\alpha}|\le s}
M_{\bar{d}}^{(\bar{\alpha})} (\ol,\, \ok) \otimes_{\Sy_\ok}
k[\Sc_\kj ] \otimes_{\Sy_\oj} \Po_{\bar{e}}(\oj,\, \oi),$$
where the direct sum $(\Xi)$ is over the integers $m,\,n$ and $N$ and
the tuples $\ol,\, \ok,\, \oj,\, \oi$.

\begin{lem}
The filtration $F_s$ is stable under the differential $d$ of $L$.
\end{lem}

\begin{proo}
The differential $d$ is the sum of three terms.
\begin{itemize}
\item The differential $\delta_M$, induced by the one of M, decreases the homological degree
$|\bar{d}|$ by $1$. Therefore, we have $\delta_M(F_s)\subset
F_{s-1}$.
\item The differential $\delta_\Po$, induced by the one of $\Po$, decreases the homological
degree $|\bar{e}|$ by $1$. We have $\delta_\Po(F_s)\subset F_s$.
\item By definition of the analytic quasi-free $\Po$-modules, one has that the
morphism $d_\theta$ decreases the total graduation $|\bar{\alpha}|$  by at least $1$ and
the homological degree $|\bar{d}|$ by $1$. We have $d_\theta(F_s)\subset F_{s-2}$.
\end{itemize}
\end{proo}

This filtration gives a spectral sequence $E^*_{s,\, t}$. The first term of this spectral
sequence $E^0_{s,\, t}$ is isomorphic to the following
dg-$\Sy$-bimodule :
$$E^0_{s,\, t}=\bigoplus_\Xi \bigoplus_{r=0}^s \bigoplus_{|\bar{d}|=s-r, \,
|\bar{e}|=t+r \atop |\bar{\alpha}|=r}  M_{\bar{d}}^{(\bar{\alpha})} (\ol,\, \ok)
\otimes_{\Sy_\ok} k[\Sc_\kj ] \otimes_{\Sy_\oj}
\Po_{\bar{e}}(\oj,\, \oi).$$
Using the same notations as in section~\ref{compositionproductofdg}, we have
$$E^0_{s,\,t}=\bigoplus_{r=0}^s I^0_{s-r,\, t+r}\left(\bigoplus_{|\bar{\alpha}|=r}
 M^{(\bar{\alpha})}\boxtimes_c \Po\right).$$

On can decompose the weight graded differential $\Sy$-bimodule $L$ with its weight
$L=\bigoplus_{\rho} L^{(\rho)}$. Since this decomposition is stable under
the filtration $F_s$, the spectral sequence $E^*_{s,_, t}$ is isomorphic to
$$E^*_{s,\,t}(L)=\bigoplus_{\rho \in \mathbb{N}} E^*_{s,\, t}(L^{(\rho)}).$$
By definition of the analytic quasi-free $\Po$-modules, the weight of the
elements of $V$ is at least $1$. Hence, we have
$$E^0_{s,\,t}(L^{(\rho)})=\bigoplus_{r=0}^{\textrm{min}(s,\, \rho)}
I^0_{s-r,\, t+r}\left(\big( \bigoplus_{|\bar{\alpha}|=r}
M^{(\bar{\alpha})} \boxtimes_c \Po \big) \bigcap L^{(\rho)}\right)
.$$

\begin{pro}
\label{SuiteSpectrale}
Let $L$ be an analytic quasi-free $\Po$-module. The spectral sequence $E^*_{s,\,t}$
converges to the homology of $L$
$$E^*_{s,\, t}\Longrightarrow H_{s+t}(L,\,d) .$$
The differential $d^0$ on $E^0_{s,\, t}$ is given by $\delta_\Po$ and the differential
$d^1$ on $E^1_{s,\, t}$ is given by $\delta_M$. Therefore, we have
$$ E^2_{s,\,t}=\bigoplus_{r=0}^s I^2_{s-r,\, t+r}\left(\bigoplus_{|\bar{\alpha}|=r}
M^{(\bar{\alpha})}\boxtimes_c \Po\right),$$
which can be written
$$E^2_{s,\,t}=\bigoplus_{r=0}^s \bigoplus_{|\bar{\alpha}|=r}
\bigoplus_{ |\bar{d}|=s-r \atop |\bar{e}|=t-r} (
H_*(M))^{(\bar{\alpha})}_{\bar{d}} \boxtimes_c
(H_*(\Po))_{\bar{e}}.$$
\end{pro}

\begin{proo}
Since the filtration is exhaustive $\bigcup_s F_s = L$ and bounded below
$F_{-1}=0$, the classical theorem of convergence of spectral sequences
(\emph{cf.} \cite{Weibel} 5.5.1) shows that $E^*_{s,\,t}$ converges
to the homology of $L$.

The form of the differentials $d^0$ and $d^1$ comes from the study
of the action of $d=\delta_M+\delta_\Po+d_\theta$ on the
filtration $F_s$. And the formula of $E^2_{s,\,t}$ comes from
Proposition~\ref{I2II2}.
\end{proo}

\subsubsection{\bf Comparison lemma for analytic quasi-free modules}

To show the comparison lemma, we need to prove the following lemma first.

\begin{lem}
Let $\Psi\, :\, \Po \to \Po'$ be a morphism of weight graded connected dg-properads.
Let $(L,\, \lambda)$ be an analytic quasi-free $\Po$-module and $(L',\, \lambda')$ be
an analytic quasi-free $\Po'$-module. Denote $\bar{L}=M$ and $\bar{L}'=M'$ the
indecomposable quotients. We consider the action of $\Po$ on
$L'$ induced by the morphism $\Psi$. Let $\Phi \, :\, L \to L'$ be a morphism of
weight graded differential $\Po$-modules.
\begin{enumerate}

\item The morphism $\Phi$ preserves the filtration $F_s$ and leads to a morphism
of spectral sequences
$$E^*(\Phi)\ :\ E^*(L) \to E^*(L').$$

\item Let $\bar{\Phi} \, :\, M\to M'$ be the morphism of
dg-$\Sy$-bimodules induces by $\Phi$. The morphism
$E^0(\Phi)\, : \, M\boxtimes_c \Po \to M'\boxtimes_c \Po'$
is equal to $E^0=\bar{\Phi}\boxtimes_c \Psi$.
\end{enumerate}
\end{lem}

\begin{proo}
$\ $
\begin{enumerate}

\item Let $(m_1,\ldots,\, m_b)\, \sigma \,(p_1,\ldots,\, p_a)$ be an element of
$F_s(M\boxtimes_c \Po)$. It means that
$|\bar{d}|+|\bar{\alpha}|\le s$. Since, $\Phi$ is a morphism of
$\Po$-modules, we have
\begin{eqnarray*}
&& \Phi\big((m_1,\ldots,\, m_b)\, \sigma \,(p_1,\ldots,\,
p_a)\big) =
\Phi\circ \lambda\big((m_1,\ldots,\, m_b)\, \sigma \,(p_1,\ldots,\, p_a)\big) =\\
&& \lambda' \big( ( \Phi(m_1),\ldots,\, \Phi(m_b))\, \sigma
\,(\Psi(p_1),\ldots,\, \Psi(p_a))\big).
\end{eqnarray*}
Denote $$\Phi(m_i)=(m_1^i,\ldots,\, m_{b_i}^i)\, \sigma^i\,
(p_1^i,\ldots,\, p_{a_i}^i).$$
The application $\Phi$ is a morphism of
dg-$\Sy$-bimodules. Therefore, we have $d_i=|\bar{d^i}|+|\bar{e^i}|$ which implies
$|\bar{d^i}|\le d_i$ and $\sum_i |\bar{d^i}| \le |\bar{d}|$.
Since $\Phi$ is a morphism of weight graded $\Po$-modules, we have $\sum_i |\bar{\alpha^i}|\le
|\bar{\alpha}|$. We conclude that
$$\Phi\big((m_1,\ldots,\, m_b)\, \sigma \,(p_1,\ldots,\, p_a)\big)\in F_s(L').$$

\item The morphism $E^0_{s,\,t}(\Phi)$ is given by the following morphism on the quotients
$$E^0_{s,\,t} \ : \ F_s(L_{s+t})/F_{s-1}(L_{s+t}) \to F_s(L'_{s+t})/F_{s-1}(L'_{s+t}).$$
Let $(m_1,\ldots,\, m_b)\, \sigma \,(p_1,\ldots,\, p_a)$ be
an element of $E^0_{s,\,t}$, where
$|\bar{\alpha}|=r\le s$, $|\bar{d}|=s-r$ and $|\bar{e}|=t+r$. We have
$$E^0_{s,\, t}(\Phi)\big( (m_1,\ldots,\, m_b)\, \sigma
\,(p_1,\ldots,\, p_a)\big)=0$$
if and only if $\Phi\big(
(m_1,\ldots,\, m_b)\, \sigma \,(p_1,\ldots,\, p_a)\big)\in
F_{s-1}(L')$. We have seen that
$$ \Phi\big((m_1,\ldots,\, m_b)\, \sigma \,(p_1,\ldots,\, p_a)\big) =
\lambda' \big( ( \Phi(m_1),\ldots,\, \Phi(m_b))\, \sigma
\,(\Psi(p_1),\ldots,\, \Psi(p_a))\big).$$
Therefore $\Phi\big((m_1,\ldots,\, m_b)\, \sigma \,(p_1,\ldots,\, p_a)\big)$
belongs to $F_s(L'_{s+t})\backslash F_{s-1}(L'_{s+t})$ if and only if
$\Phi(m_i)=m_i=\bar{\Phi}(m_i)$. Consider
$\Phi(m_i)=(m_1^i,\ldots,\, m_{b_i}^i)$ $\sigma^i$
$(p_1^i,\ldots,\, p_{a_i}^i)$ where at least one $p_j^i$
does not belong to $I$. By definition of a weight graded properad, the weight of
these $p_j^i$ is at least $1$. Since $\Phi$ preserves the total weight , we have
$|\bar{\alpha^i}|<\alpha^i$. Finally, the morphism $E^0(\Phi)$
is equal to $\bar{\Phi}\boxtimes_c \Psi$.
\end{enumerate}
\end{proo}

\begin{thm}[Comparison lemma for analytic quasi-free modules]
\label{lemmecomparaisonquasi-freemodules}
With the same notations as in the previous lemma, if two of the following morphisms
$$\left\{ \begin{array}{l}
\Psi \ : \ \Po \to \Po', \\
\bar{\Phi} \ : \ M \to M', \\
\Phi \ : \ L \to L',
\end{array}\right.$$
are quasi-isomorphisms, then the third one is also a quasi-isomorphism.
\end{thm}

\begin{rem}
If the properads are not weight graded and if $\Phi=\bar{\Phi}\boxtimes_c \Psi$, when
$\bar{\Phi}$ and $\Psi$ are quasi-isomorphisms, then $\Phi$ is
also a quasi-isomorphism. The proof is straight. We have immediately
$E^0(\Phi)=\bar{\Phi}\boxtimes_c \Psi$. Since
$\bar{\Phi}$ and $\Psi$ are quasi-isomorphisms, we get that
$E^2(\Phi)$ is an isomorphism by
Proposition~\ref{SuiteSpectrale} ($d^0=\delta_\Po$ et
$d^1=\delta_M$). And with the same proposition, the convergence of the
spectral sequence $E^*$ shows that $\Phi$
is a quasi-isomorphism between $L$ and $L'$.
\end{rem}

\begin{proo}
$ $
\begin{enumerate}

\item  Suppose that $\Psi  \, : \, \Po \to \Po'$ and
$\bar{\Phi} \, : \, M \to M'$ are quasi-isomorphisms. The previous lemma
shows that $\Phi$ induces a morphism of spectral sequences
$E^*(\Phi) \, :\, E^*(L)\to E^*(L')$ and that $E^0(\Phi)=\bar{\Phi}\boxtimes_c \Psi$.
We conclude with the same arguments.\\

\item Suppose that $\Psi  \, : \, \Po \to \Po'$ and
$\Phi \, : \, L \to L'$ are quasi-isomorphisms. We have seen that the spectral
sequence $E^*_{s,\, t}$
preserves the decomposition $L^{(\rho)}$. Moreover, we have
$$E^2_{s,\,t}(L^{(\rho)})=\bigoplus_\Xi \bigoplus_{r=\textrm{max}
(0,\, -t)}^{\textrm{min}(s,\, \rho)}I^2_{s-r,\, t+r}\left(
\bigoplus_{\chi}  M_{\bar{d}}^{(\bar{\alpha})} (\ol,\, \ok)
\otimes_{\Sy_\ok} k[\Sc_\kj ] \otimes_{\Sy_\oj}
\Po_{\bar{e}}^{(\bar{\beta})}(\oj,\, \oi)  \right), $$
where the second sum $(\chi)$ is on $|\bar{\alpha}|=r$,
$|\bar{\beta}|=\rho-r$, $|\bar{d}|=s-r$ and $|\bar{e}|=t+r$. When $s\ge\rho$,
we have for $r=\rho$
\begin{eqnarray*}
 I^2_{s-\rho,\, t+\rho}\left( \bigoplus_{\chi}
M_{\bar{d}}^{(\bar{\alpha})} (\ol,\, \ok) \otimes_{\Sy_\ok}
k[\Sc_\kj ] \otimes_{\Sy_\oj} \Po_{\bar{e}}^{(\bar{\beta})}(\oj,\,
\oi)\right) &=& I^2_{s-\rho,\, t+\rho} \left(
\bigoplus_{|\bar{\alpha}|=\rho
\atop |\bar{d}|=s-\rho} M_{\bar{d}}^{(\bar{\alpha})} \right)\\
&=& \left\{ \begin{array}{ll} H_{s-\rho}(M^{(\rho)}) &
\textrm{if}\
 t=-\rho, \\
0 & \textrm{otherwise}.
\end{array} \right.
\end{eqnarray*}

Finally, we have
\begin{eqnarray*}
\left\{ \begin{array}{ll}
 E^2_{s,\,t}(L^{(\rho)})=0 &
\textrm{if}\
 t<-\rho, \\
E^2_{s,\, -\rho}(L^{(\rho)})=H_{s-\rho}(M^{(\rho)}) & \textrm{if}
\ s\ge \rho.
\end{array} \right.
\end{eqnarray*}

We now show by induction on $\rho$ that
$\bar{\Phi}^{(\rho)} \, :\, M^{(\rho)} \to M'^{(\rho)}$ is a quasi-isomorphism
$H_*(\bar{\Phi}^{(\rho)}) \, :\, H_*(M^{(\rho)}) \to H_*(M'^{(\rho)})$.\\

For $\rho=0$, since  $L^{(0)}=M^{(0)}$, we have
$\bar{\Phi}^{(0)}=\Phi^{(0)}$, which is a quasi-isomorphism.\\

Suppose now that the result is true for $r<\rho$ in every (homological) degree $*$
and for $r=\rho$ in degree $*<d$.
Let show that it is also true for $*=d$. (Since every chain complexes are
concentrated in non-negative degree, we have that $H_s(\bar{\Phi}^{(\rho)})$
is always an isomorphism for $s<0$).

By the preceding lemma, we have
$$E^2_{s,\,t}(\Phi^{(\rho)})=\bigoplus_{r=0}^{\textrm{min}(s,\, \rho)}
I^2_{s-r,\, t+r}\left( \bigoplus_{|\bar{\alpha}|=r}
\bar{\Phi}^{{(\bar{\alpha})}}\boxtimes_c \Psi \right).$$
With the induction hypothesis, we get that
$E^2_{s,\,t}(\Phi^{(\rho)})\, :\, E^2_{s,\,t}(L^{(\rho)})\to
E^2_{s,\,t}(L'^{(\rho)})$ is an isomorphism for $s<d+\rho$.
Let show that $$E^2_{d+\rho,\,-\rho}(\Phi^{(\rho)})\, :\,
E^2_{d+\rho,\,-\rho}(L^{(\rho)})\to
E^2_{d+\rho,\,-\rho}(L'^{(\rho)})
$$ is still an isomorphism. To do that, we consider the mapping cone of $\Phi^{(\rho)} :$
$\cone(\Phi^{(\rho)})=\Sigma^{-1}L^{(\rho)}\oplus L'^{(\rho)}$.
On this cone, we define the following filtration
$$F_s(\cone(\Phi^{(\rho)})):=F_{s-1}(\Sigma^{-1}L^{(\rho)})\oplus F_s(L'^{(\rho)}) .$$
This filtration induces a spectral sequence which verifies
$$E^1_{*,\, t}(\cone(\Phi^{(\rho)}))=\cone(E^1_{*,\,t}(\Phi^{(\rho)})).$$
The mapping cone of $E^1_{*,\,t}(\Phi^{(\rho)})$ verifies the long exact sequence
\begin{eqnarray*}
&& \cdots \to H_{s+1}\big(\cone(E^1_{*,\, t}(\Phi^{(\rho)}))\big) \to
H_s\big(E^1_{*,\,t}(L^{(\rho)})\big)  \xrightarrow{H_s\big(E^1_{*,\,t}(\Phi^{(\rho)}) \big)} \\
&& H_s\big(E^1_{*,\,t}(L'^{(\rho)})\big) \to
H_s\big(\cone(E^1_{*,\, t}(\Phi^{(\rho)}))\big) \to
H_{s-1}\big(E^1_{*,\,t}(L^{(\rho)})\big) \to \cdots
\end{eqnarray*}
which finally gives the following long exact sequence $(\xi)$  :
\begin{eqnarray*} && \xymatrix@C=35pt{\cdots
\ar[r]& E^2_{s+1,\, t}(\cone(\Phi^{(\rho)})) \ar[r]&
E^2_{s,\,t}(L^{(\rho)}) \ar[r]^{E^2_{s,\,t}(\Phi^{(\rho)})} & E^2_{s,\,t}(L'^{(\rho)}) }\\
&& \xymatrix@C=35pt{\quad  \ar[r]&E^2_{s,\,
t}(\cone(\Phi^{(\rho)})) \ar[r] & E^2_{s-1,\,t}(L^{(\rho)})
\ar[r]& \cdots}
\end{eqnarray*}
We have seen that for all $t<-\rho$,
$$E^2_{s,\,t}(L^{(\rho)})=E^2_{s,\,t}(L'^{(\rho)})=0.$$
Therefore, the long exact sequence $(\xi)$ shows that
$E^2_{s,\,t}(\cone(\Phi^{(\rho)}))=0$ for all $s$, when
$t<-\rho$.

In the same way, we have seen that
$E^2_{s,\,t}(\Phi^{(\rho)})$ is an isomorphism for
$s<d+\rho$. Thanks to the long exact sequence $(\xi)$, we get that
$E^2_{s,\,t}(\cone(\Phi^{(\rho)}))=0$ for all $t$, when $s<d+\rho$.

These two results prove that
$$
\left\{
\begin{array}{l}
E^2_{d+\rho,\, -\rho}(\cone(\Phi^{(\rho)}))=E^\infty_{d+\rho,\,
-\rho}(\cone(\Phi^{(\rho)})), \\
E^2_{d+\rho+1,\,
-\rho}(\cone(\Phi^{(\rho)}))=E^\infty_{d+\rho+1,\,
-\rho}(\cone(\Phi^{(\rho)})).
\end{array} \right.$$

Since the filtration $F_s$ of the mapping cone $\Phi^{(\rho)}$ is
bounded below and exhaustive, we know that the spectral sequence
$E^*_{s,\,t}(\cone(\Phi^{(\rho)}))$ converges to the homology of this mapping cone.
The morphism $\Phi^{(\rho)}$ is a
quasi-isomorphism. Hence this homology is null and we have the equalities  :
$$\left\{ \begin{array}{l}
E^2_{d+\rho,\, -\rho}(\cone(\Phi^{(\rho)}))=E^\infty_{d+\rho,\,
-\rho}(\cone(\Phi^{(\rho)}))=0, \\
E^2_{d+\rho+1,\,
-\rho}(\cone(\Phi^{(\rho)}))=E^\infty_{d+\rho+1,\,
-\rho}(\cone(\Phi^{(\rho)}))=0.
\end{array} \right. $$
If we inject these equalities in the long exact sequence
$(\xi)$, we show that $E^2_{d+\rho,\, -\rho}(\Phi^{(\rho)})$ is an
isomorphism.

We conclude using the fact that $E^2_{d+\rho,\,
-\rho}(L^{(\rho)})=H_d(M^{(\rho)})$ and that $E^2_{d+\rho,\,
-\rho}(L'^{(\rho)})=H_d(M'^{(\rho)})$.
Therefore, $E^2_{d+\rho,\,
-\rho}(\Phi^{(\rho)})$ is equal to $H_d\big(
\bar{\Phi}^{(\rho)}\big)\, : \, H_d(M^{(\rho)}) \to
H_d(M'^{(\rho)})$ which is an isomorphism.\\

\item Suppose that $\Phi \, : \, L \to L'$ et $\bar{\Phi} \, :\,
M \to M'$ are quasi-isomorphisms. We are going to use the same methods.
Since $M^{(0)}=I$, we get from the relation
$$E^2_{s,\,t}(L^{(\rho)})=\bigoplus_\Xi
\bigoplus_{r=\textrm{max}(0,\, -t)}^{\textrm{min}(s,\, \rho)}
I^2_{s-r,\, t+r}\left( \bigoplus_\chi M_{\bar{d}}^{(\bar{\alpha})}
(\ol,\, \ok) \otimes_{\Sy_\ok} k[\Sc_\kj ] \otimes_{\Sy_\oj}
\Po_{\bar{e}}^{(\bar{\beta})}(\oj,\, \oi)  \right), $$
where the second direct sum $\chi$ is on $|\bar{\alpha}|=r$,
$|\bar{\beta}|=\rho-r$, $|\bar{d}|=s-r$ and $|\bar{e}|=t+r$, that for $s=0$
$$E^2_{0,\,t}(L^{(\rho)})= I^2_{0,\, t}(\Po_t^{(\rho)})=H_t(\Po^{(\rho)}). $$
One can see that
$$E^2_{s,\,t}(L^{(\rho)})=0,$$
for $s<0$.

Let show that the morphism $\Psi^{(\rho)} \, : \, \Po^{(\rho)}\to \Po'^{(\rho)}$
are quasi-isomorphisms, by induction on $\rho$.\\

By definition of connected properads, we have that
$\Po^{(0)}=\Po'^{(0)}=I$. Hence, the morphism $\Psi^{(0)}=id_I$ is a
quasi-isomorphism.\\

Suppose that the result is true for every
$r<\rho$. We show by induction on $t$ that
$$H_t(\Psi^{(\rho)})\ :\  H_t(\Po^{(\rho)}) \to H_t(\Po'^{(\rho)}) $$
is an isomorphism. We know that it is true for
$t<0$. Suppose that it is still true for $t<e$. By the previous lemma, we have
$$E^2_{s,\,t}(\Phi^{(\rho)})=\bigoplus_{r=0}^{\textrm{min}(s,\, \rho)}
I^2_{s-r,\, t+r}\left( \bigoplus_{|\bar{\alpha}|=r \atop
|\bar{\beta}|=\rho-r} \bar{\Phi}^{(\bar{\alpha})}\boxtimes_c
\Psi^{(\bar{\beta})}
 \right).$$
By the induction hypothesis, we have that
$E^2_{s,\,t}(\Phi^{(\rho)}) \, : \, E^2_{s,\,t}(L^{(\rho)}) \to
E^2_{s,\,t}(L'^{(\rho)})$ is an isomorphism for every $s$,
when $t<e$. By injecting this in the long exact sequence
$(\xi)$, we show that $E^2_{s,\,t} (\cone(\Phi^{(\rho)}))=0$
for all $s$, when $t<e$.

The form of the filtration $F_s(\cone(\Phi^{(\rho)}))=F_{s-1}(\
\Sigma^{-1}L^{(\rho)})\oplus F_s(L'^{(\rho)})$ and the fact that
$F_{-1}(L^{(\rho)})=0$ show that
$E^2_{s,\,t}(\cone(\Phi^{(\rho)}))=0$ for $s<0$.

These two results with the convergence of the spectral sequence
$E^*_{s,\,t}(\cone(\Phi^{(\rho)}))$ to the homology of the cone
of $\Phi^{(\rho)}$, which is null, show that
$$
\left\{
\begin{array}{l}
E^2_{0,\, e}(\cone(\Phi^{(\rho)}))=E^\infty_{0,\,
e}(\cone(\Phi^{(\rho)}))=0 \\
E^2_{1,\, e}(\cone(\Phi^{(\rho)}))=E^\infty_{1,\,
e}(\cone(\Phi^{(\rho)}))=0.
\end{array} \right.
$$

When we inject these two equalities in the long exact sequence
$(\xi)$, we get that the morphism
$$\xymatrix@C=45pt{E^2_{0,\,e}(L^{(\rho)}) \ar[r]^{E^2_{0,\, e}(\Phi^{(\rho)})} &
E^2_{0,\,e}(L'^{(\rho)}) }$$
is an isomorphism. We conclude using the equalities
$E^2_{0,\,e}(L^{(\rho)})=H_e(\Po^{(\rho)})$,
 $E^2_{0,\,e}(L'^{(\rho)})=H_e(\Po'^{(\rho)})$ and
$E^2_{0,\, e}(\Phi^{(\rho)})=H_e(\Psi^{(\rho)}).$
\end{enumerate}
\end{proo}

\subsubsection{\bf Comparison lemma for analytic quasi-cofree comodules}

We can prove the same kind of result for analytic quasi-cofree comodule.
Let $\Co$ be a connected dg-coproperad and let
$L=M\boxtimes_c \mathcal{C}$ be an analytic quasi-cofree $\mathcal{C}$-comodule
on the right, where $M=\Upsilon(V)$. On this chain complex, we define the following filtration
$$F'_s(L):=\bigoplus_{(\Xi)} \bigoplus_{|\bar{e}|+|\bar{\beta}|\le s}
M_{\bar{d}}(\ol,\, \ok)\otimes_{\Sy_\ok} k[\Sc_\kj]
\otimes_{\Sy_\oj}\mathcal{C}_{\bar{e}}^{(\bar{\beta})}(\oj,\,\oi),
$$
where the direct sum $(\Xi)$ is on the integers $m$, $n$ and
$N$ and on the tuples $\ol$, $\ok$, $\oj$, $\oi$, $\bar{d}$,
$\bar{e}$ and $\bar{\beta}$ such that $|\bar{e}|+|\bar{\beta}|\le
s$.

Applying the same arguments to this filtration and the associated spectral sequence, we
get the following theorem.

\begin{thm}[Comparison lemma for analytic quasi-cofree comodules]
\label{lemmecomparaisoncomodule}

Let $\Psi\, :\, \mathcal{C} \to \mathcal{C}'$ be a morphism
of weight graded connected dg-coproperad and let $(L,\, \lambda)$
and $(L',\, \lambda')$ be two analytic quasi-cofree comodules
on $\mathcal{C}$ and $\mathcal{C}'$. Denote $\bar{L}=M$ and
$\bar{L'}=M'$ the indecomposable quotients. Let
$\Phi \, :\, L \to L'$ be a morphism of analytic $\mathcal{C}'$-comodules,
where the structure of $\mathcal{C}'$-comodule on  $L$
is given by the functor $\Psi$. Denote
$\bar{\Phi} \, :\, M\to
M'$ the morphism of dg-$\Sy$-bimodules induced by $\Phi$. \\

If two of the three morphisms $\left\{ \begin{array}{l}
\Psi \ : \ \mathcal{C} \to \mathcal{C}', \\
\bar{\Phi} \ : \ M \to M,' \\
\Phi \ : \ L \to L',
\end{array} \right. $
are quasi-isomorphisms, then the third one is a quasi-isomorphism.
\end{thm}

\subsubsection{\bf Bar-cobar resolution}

In this section, we use the two last theorems to prove that
the bar-cobar construction on a weight graded properad $\Po$ is a
resolution of $\Po$. Notice that this result is a generalization to properads
of one given by V. Ginzburg and M.M. Kapranov for operads in \cite{GK}.\\

We define the morphism $\zeta \, :\, \bar{\B}^c(\bar{\B}(\Po))\to \Po$ by the
following composition
\begin{eqnarray*}
\xymatrix{\bar{\B}^c(\bar{\B}(\Po))=\F\big(\Sigma^{-1}
\bar{\F}^c(\Sigma \oPo)\big) \ar@{>>}[r] & \F\big(\Sigma^{-1}
\bar{\F}_{(1)}^c(\Sigma \oPo)\big)=\F(\oPo) \ar[r]^(0.75){c_\Po} &
\Po,}
\end{eqnarray*}
where the morphism $c_\Po$ corresponds to the compositions $\mu$
of all the elements of $\Po$ that are indexing a graph of
$\F(\oPo)$.

\begin{pro}
The morphism $\zeta$ is a morphism of weight graded differential properads.
\end{pro}

\begin{proo}
Since the definition of $\zeta$ is based on the composition $\mu$ of $\Po$,
one can see that $\zeta\circ
\mu_{\F(\Sigma^{-1} \bar{\F}^c(\Sigma \Po))} =\mu_\Po\circ
(\zeta\boxtimes_c\zeta)$. Hence, the morphism $\zeta$
is a morphism of properads.

The morphism $\zeta$ is either null or given by compositions
of operations which preserves the weight. Therefore, the morphism $\zeta$
is a morphism of weight graded properads.

It remains to show that $\zeta$ commutes with the differentials.

On $\bar{\B}^c(\bar{\B}(\Po))=\F\big(\Sigma^{-1}
\bar{\F}^c(\Sigma \oPo)\big)$  the differential $d$ is the sum
of the derivation  $d_{\theta'}$ of the cobar
construction induced by the partial coproduct of $\F^c(\Sigma
\oPo)$ with the canonical differential
$\delta_{\bar{\B}(\Po)}$. And the canonical differential
$\delta_{\bar{\B}(\Po)}$ is the sum of the coderivation
$d_\theta$ of the bar construction $\bar{\B}(\Po)$ induced by the
partial product of $\Po$ with the canonical differential
$\delta_\Po$. Let $\xi$ be an element of $\F\big(\Sigma^{-1}
\bar{\F}^c(\Sigma \oPo)\big)$ represented by an indexed graph $g$.
Three different cases are possible.

\begin{itemize}

\item When all the vertices of $g$ are indexed by elements of
$\F^c_{(1)}(\Sigma \oPo)=\Sigma \oPo$, the image of $\xi$ under the differential $d$
is equal to the image under the canonical differential $\delta_\Po$.
Since the composition $\mu$ of $\Po$ is a morphism of differential properads,
we have $\zeta\circ d(\xi)= \zeta \circ
\delta_\Po(\xi)=\delta_\Po\circ \zeta(\xi)$.

\item When at least one vertex of $g$ is indexed by a element of
$\F^c_{(s)}(\Sigma \oPo)$, with $s\ge 3$, the image of $\xi$ under
$\zeta$ is null. Moreover, the image of $\xi$ under the differential $d$
is a sum of graphs where at least one vertex is indexed by an element of
$\F^c_{(s)} (\Sigma \oPo)$, with
$s\ge 2$. Therefore, we have $d\circ \zeta(\xi) + \zeta\circ d(\xi)=0$.

\item When the vertices of $g$ are indexed by elements of
$\F^c_{(s)}(\Sigma \oPo)$ where $s=1,\, 2$, with at least one of
$\F^c_{(2)}(\Sigma \oPo)$, the image of $\xi$ under $\zeta$ is null.
It remains to show that $\zeta\circ d(\xi)=0$. Since
$\delta_\Po$ preserves the natural graduation of $\F^c$, one has
$\zeta\circ \delta_\Po (\xi)=0$. Study the image under
$d_\theta+d_{\theta'}$ of an element of $\F^c_{(2)}(\Sigma
\oPo)$. We denote $\xi=X\otimes \Sigma^{-1}(\Sigma p_1\otimes \Sigma
p_2)\otimes Y$, where $(\Sigma p_1\otimes \Sigma p_2)$ belongs to
$\F_{(2)}^c(\Sigma \oPo)$ and $X$, $Y$ to $\F(\Sigma^-1
\bar{\F}^c(\Sigma \oPo))$. Applying $d_{\theta'}$, we get
one term of the form
$$(-1)^{|X|+1}(-1)^{|p_1|+1}X\otimes \Sigma^{-1}\Sigma p_1\otimes
 \Sigma^{-1}\Sigma p_2\otimes Y, $$
where $\Sigma^{-1}\Sigma p_1 $ and $\Sigma^{-1}\Sigma p_2 $ are elements of
$\Sigma^{-1} \F^c_{(1)}(\Sigma \oPo)$. And applying $d_\theta$, we get
one term of the form
$$(-1)^{|X|+1}(-1)^{|p_1|}X\otimes \Sigma^{-1} \Sigma \mu(p_1\otimes p_2)
\otimes Y, $$ where  $\Sigma^{-1} \Sigma \mu(p_1\otimes p_2) $
belongs to $\Sigma^{-1} \F^c_{(1)}(\Sigma \oPo)$. Therefore,
$\zeta\circ (d_\theta +d_{\theta'}) (\xi)$ is a sum of terms of the form
$$\big( (-1)^{|X|+|p_1|} +(-1)^{|X|+|p_1|+1}\big) X\otimes \mu(p_1\otimes p_2)
\otimes Y=0.$$
\end{itemize}
\end{proo}

\begin{rem}
With this proof, one can understand the use of the sign $-1$ in the definition of the
derivation $d_{\theta'}$ of the cobar construction.
\end{rem}

Denote $\mathcal{C}:=\bar{\B}(\Po)$ the connected dg-coproperad
defined by the bar construction
on $\Po$. We consider the augmented cobar construction
on the right $\bar{\B}^c(\mathcal{C})\boxtimes_c
\mathcal{C}=\bar{\B}^c(\bar{\B}(\Po))\boxtimes_c \bar{\B}(\Po)$.
Remark that the coproperad $\bar{\B}(\Po)=\F^c(\Sigma \oPo)$
is weight graded (by the global weight) and that the augmented
cobar construction
$\bar{\B}^c(\mathcal{C})\boxtimes_c \mathcal{C}$ is an analytic quasi-cofree
$\mathcal{C}$-comodule.

\begin{lem}
The morphism $ \zeta\boxtimes_c \bar{\B}(\Po) \ : \
\bar{\B}^c(\bar{\B}(\Po))\boxtimes_c \bar{\B}(\Po) \to
\Po\boxtimes_c \bar{\B}(\Po)$ between the augmented cobar construction
$\bar{\B}^c(\mathcal{C})\boxtimes_c \mathcal{C}$ and
the augmented bar construction $\Po\boxtimes_c \bar{\B}(\Po)$ is
a morphism of dg-$\Sy$-bimodules. Moreover, it is a
morphism of analytic quasi-free $\bar{\B}(\Po)$-comodules.
\end{lem}

\begin{proo}
Since $\zeta$ and $id_{\bar{\B}(\Po)}$ are
morphisms of dg-$\Sy$-bimodules, the morphism $\zeta\boxtimes_c
id_{\bar{\B}(\Po)}$ preserves the canonical differentials
$\delta_{\bar{\B}^c(\mathcal{C})}+\delta_{\mathcal{C}}$ and
$\delta_\Po + \delta_{\bar{\B}(\Po)}$. One can see that the morphism
$d_{\theta_r}$ that appears in the definition of the
differential of the cobar construction
$\bar{\B}^c(\mathcal{C})\boxtimes_c \mathcal{C}$ (\emph{cf.}
\ref{Defsbarcobar}) corresponds  via $\zeta
\boxtimes_c \bar{\B}(\Po)$ to the morphism $d_{\theta_l}$ of
the differential of the bar construction
$\Po\boxtimes_c \bar{\B}(\Po)$. (The morphism $d_{\theta_r}$
corresponds to extract an operation of $\Po$ from $\Co=\bar{\B}(\Po)$ from
the bottom and compose it, from the top, in
$\bar{\B}^c(\bar{\B}(\Po))$. The morphism $d_{\theta_l}$ corresponds to
extract an operation of $\Po$ in $\bar{\B}(\Po)$ from the bottom to compose
it at the top of $\Po$.)

Since $\zeta$ and $id_{\bar{\B}(\Po)}$ preserve the
weight graduation coming from the one of $\Po$, we have that
$\zeta\boxtimes_c \bar{\B}(\Po)$ is a morphism of
analytic quasi-free $\bar{\B}(\Po)$-comodules.
\end{proo}

\begin{thm}[Bar-cobar resolution]
\label{barcobarresolution}
For every connected dg-properad $\Po$,
the morphism $$\zeta \ :\ \bar{\B}^c(\bar{\B}(\Po))\to \Po$$
is a quasi-isomorphism of weight graded dg-properads.
\end{thm}

\begin{proo}
We know from Theorem~\ref{acyclicitécobarcoaugmentée} that the
augmented cobar construction $
\bar{\B}^c(\bar{\B}(\Po))\boxtimes_c \bar{\B}(\Po)$ is acyclic and
from Theorem~\ref{acyclicitébaraugmentée} that the augmented bar
construction $\Po\boxtimes_c \bar{\B}(\Po)$ is acyclic too.
Therefore, the morphism $\zeta\boxtimes_c\bar{\B}(\Po)$ is a
quasi-isomorphism. We apply the comparison lemma for analytic
quasi-cofree comodules (Theorem~\ref{lemmecomparaisoncomodule})
 to the quasi-isomorphisms $\Psi=id_{\bar{\B}(\Po)}$ and $\Phi=\zeta\boxtimes_c
\bar{\B}(\Po)$. It gives that $\zeta=\bar{\Phi}$ is a quasi-isomorphism.
\end{proo}

\begin{rem} In the rest of the text, we will mainly apply this theorem
to quadratic properads (which are naturally graded). The
comparison lemma for quasi-free properads, proved in the next section, will allow us
to simplify the bar-cobar resolution to get a quasi-free resolution, namely the
Koszul resolutions (\emph{cf.} Section~\ref{KoszulDuality}).
\end{rem}

\subsection{Comparison lemma for quasi-free properad}

We show here a comparison lemma of the same type but for
quasi-free properads.

\begin{thm}[Comparison lemma for quasi-free properads]
\label{lemmecomparaisonproperades}
Let $M$ and $M'$ be two weight graded
dg-$\Sy$-bimodules such that $M^{(0)}=M'^{(0)}=0$.
Let $\Po$ and $\Po'$ be two quasi-free properads of the form
$\Po=\F(M)$ and $\Po'=\F(M')$, endowed with derivations
$d_\theta$ and  $d_{\theta'}$ coming from decomposable morphisms of
weight graded dg-$\Sy$-bimodules $\theta \, :\,
M \to \bigoplus_{s\ge 2}\F_{(s)}(M)$ and
 $\theta' \, :\, M' \to \bigoplus_{s\ge 2}\F_{(s)}(M')$.
Given a morphism of weight graded dg-properads $\Phi \, :\, \Po \to \Po'$, we consider
the induced morphism $\bar{\Phi} \, :\, M=\F_{(1)}(M) \to
M'=\F_{(1)}(M')$.\\

The morphism $\Phi$ is a quasi-isomorphism if and only if
$\bar{\Phi}$ is a quasi-isomorphism.
\end{thm}

\begin{proo}
Since the morphisms $\theta$ and $\theta'$ preserve the weight,
we have that the derivations $d_\theta$ and $d_{\theta'}$
preserve the weight. Therefore, the
differentials $\delta_\theta=\delta_M+d_\theta$ and
 $\delta_{\theta'}=\delta_{M'}+d_{\theta'}$ preserve the weight
and we can decompose the chain complexes $\F(M)=\bigoplus_{\rho \in
\mathbb{N}} \F(M)^ {(\rho)}$ and $\F(M')=\bigoplus_{\rho \in
\mathbb{N}} \F(M')^ {(\rho)}$ with the weight.

We introduce the following filtration
$$F_s(\F(M)):=\bigoplus_{r\ge -s} \F_{(r)}(M). $$
This filtration is compatible with the preceding decomposition.
To make the spectral sequences converge, we consider the
sub-complexes $\F(M)^
{(\rho)}$ and the filtration
$$F_s(\F(M)^{(\rho)})=\bigoplus_{r\ge -s} \F_{(r)}(M)^{(\rho)}.$$
Since the weight graded dg-$\Sy$-bimodules $M$ and $M'$ are concentrated
in weight $>0$, we have
$\F_{(r)}(M)^{(\rho)}= \F_{(r)}(M')^{(\rho)}=0$ when $r>\rho$ and
$$F_s(\F(M)^{(\rho)})=\bigoplus_{-s \le r \le \rho} \F_{(r)}(M)^{(\rho)}.$$

One can see that the filtration $F_s$ is stable under the
differential $\delta_\theta=\delta_M+d_\theta$. (We have
$\delta_M(F_s)\subset F_s$ and $d_\theta(F_s) \subset F_{s-1}$).
Therefore, this filtration induces a spectral sequence $E^*_{s,\,t}$.
The first term of this spectral sequence is given by the formula
$$E^0_{s,\,t}=F_s(\F(M)^{(\rho)}_{s+t})/F_{s-1}(\F(M)^{(\rho)}_{s+t})
=\F_{-s}(M)^{(\rho)}_{s+t},$$
and the differential $d^0$ is equal to the canonical differential
$\delta_M$. From this description and the properties of $\Phi$,
we get that $E^0_{s,\,t}(\Phi)=\F_{(-s)}(\bar{\Phi})_{s+t}$.

The filtration $F_s$ of $\F(M)^{(\rho)}$ is bounded below and exhaustive
($F_{-\rho-1}=0$ and $F_0=\F(M)^{(\rho)}$). Hence, the classical convergence theorem
of spectral sequences
(\emph{cf.} \cite{Weibel} 5.5.1) shows that the spectral sequence
$E^*_{s,\,t}$ converges to the homology of $\F(M)^{(\rho)}$.\\

We can now prove the equivalence :\\

($\Leftarrow $)   If $\bar{\Phi}\, :\, M \to M'$ is a
quasi-isomorphism, since the functor $\F$ is exact
(\emph{cf.} Proposition~\ref{freefuncteurexact}), we have that
$$\F_{(-s)}(\bar{\Phi})_{s+t}=E^0_{s,\, t}(\Phi^{(\rho)}) \, : \,
E^0_{s,\,t}(\F(M)^{(\rho)}) \to E^0_{s,\,t}(\F(M')^{(\rho)})$$ is
a quasi-isomorphism. Hence, $E^1(\Phi^{(\rho)})$ is an
isomorphism. The convergence of the spectral sequence
shows that $\Phi^{(\rho)} \, :\, \F(M)^{(\rho)} \to
\F(M')^{(\rho)}$ is a quasi-isomorphism. \\

($\Rightarrow$)  In the other way, if $\Phi$ is a
quasi-isomorphism, then each $\Phi^{(\rho)}$ is a
quasi-isomorphism. We will show by induction on
$\rho$ that $\bar{\Phi}^{(\rho)}\, :\, M^{(\rho)} \to M'^{(\rho)}$
is a quasi-isomorphism.

For $\rho=0$, we have $M^{(0)}=M'^{(0)}=0$ and the morphism
$\bar{\Phi}^{(0)}$ is a quasi-isomorphism.

Suppose that the result is true up to $\rho$ and show it for
$\rho+1$. In this case, since $E^0_{s,\,t}(\F(M)^{(\rho
+1)})= \F_{(-s)}(M)^{(\rho+1)}_{s+t}$,
one can see that $E^1_{s,\, t}(\Phi^{(\rho +1)})$ is an
isomorphism for every $t$, so long as $s<-1$.

Once again, we consider the mapping cone of the application
$\Phi^{(\rho+1)}$. We define the same filtration on the mapping cone
$\cone(\Phi^{(\rho+1)})$ as on $\F(M)^{(\rho+1)}$ :
$$F_s(\cone(\Phi^{(\rho+1)})):=F_s(\Sigma^{-1} \F(M)^{(\rho+1)})\oplus
F_s(\F(M')^{(\rho+1)}).$$

This filtration induces a spectral sequence such that
$$E^0_{s,\, *}(\cone(\Phi^{(\rho+1)}))=\cone(E^0_{s,\,*}(\Phi^{\rho+1})).$$

As in the proof of Theorem~\ref{lemmecomparaisonquasi-freemodules},
we have the following long exact sequence
\begin{eqnarray*}
\xymatrix@C=20pt{\cdots \ar[r]& E^1_{s,\,t+1}(\F(M)^{(\rho+1)})
\ar[r] & E^1_{s,\,t+1}(\F(M')^{(\rho+1)}) \ar[r] &
E^1_{s,\,t}(\cone(\Phi^{(\rho+1)})) }\\
\xymatrix@C=20pt{\ar[r] &  E^1_{s,\,t}(\F(M)^{(\rho+1)})\ar[r] &
E^1_{s,\,t} (\F(M')^{(\rho+1)}) \ar[r] &
E^1_{s,\,t-1}(\cone(\Phi^{(\rho+1)}))\ar[r]& \cdots}
\end{eqnarray*}

Since $E^1_{s,\,t}(\Phi^{(\rho+1)})$ is an isomorphism, for
$s<-1$, the long exact sequence shows that
$E^1_{s,\,t}(\cone(\Phi^{(\rho+1)})=0$ for every $t$ so long as
$s<-1$. The form of the filtration gives that
$E^1_{s,\,t}(\cone(\Phi^{(\rho+1)}))=0$ for every $t$ when $s\ge0$.
The spectral sequence $E^*_{s,\,t}(\cone(\Phi^{(\rho+1)}))$ collapses at
rank $E^1$ (only the row $s=-1$ is not null). It implies that
$$E^\infty_{-1,\,t}(\cone(\Phi^{(\rho+1)}))=
E^1_{-1,\,t}(\cone(\Phi^{(\rho+1)})).$$ Applying the classical convergence theorem
of spectral sequences, we get that the spectra sequence
$E^*_{s,\,t}(\cone(\Phi^{(\rho+1)}))$ converges
to the homology of the mapping cone of $\Phi^{(\rho+1)}$ which is null
since $\Phi^{(\rho+1)}$ is a quasi-isomorphism. Therefore we have
$E^1_{-1,\,t}(\cone(\Phi^{(\rho+1)}))=0$, for every  $t$, which injected in
the long exact sequence gives that
$$ E^0_{-1,\,t}(\Phi^{(\rho+1)})= \bar{\Phi}^{(\rho+1)} \ : \ M^{(\rho+1)} \to
 M'^{(\rho+1)}.$$
is a quasi-isomorphism.
\end{proo}

\begin{rem}
This theorem generalizes to properads a result of B. Fresse \cite{Fresse}
for ``connected'' operads ($\Po(0)=0$ and $\Po(1)=k$). This hypothesis ensures
the convergence of the spectral sequence introduced by the author.
To be able to generalize this result, we have introduced the graduation by the weight.
Therefore, Theorem~\ref{lemmecomparaisonproperades} can be applied to weight
graded operads non-necessarily ``connected'' in the sense of \cite{Fresse},
for instance operads with unitary operations or
weight graded algebras.
\end{rem}

\begin{cor}
\label{BarDeQuasiFree}
Let $\Po=\F(M)$ be a quasi-free properad verifying the hypothesis of the previous theorem.
The natural projection $\bar{\B}(\F(M))\to \Sigma M$ is a quasi-isomorphism.
\end{cor}

\begin{proo}
Denote $M':=\Sigma^{-1}\bar{\B}(\F(M))$. The bar-cobar construction of the weight graded connected
dg-properad $\Po$ is a quasi-free properad of the form $\B^c\B(\F(M)))=\F(M')$, which is quasi-isomorphic to
$\Po=\F(M)$ by Theorem~\ref{barcobarresolution}. The comparison lemma for quasi-free properad
(Theorem~\ref{lemmecomparaisonproperades}) concludes the proof.
\end{proo}

\section{Simplicial bar construction}

In this section, we generalize the simplicial bar construction of
associative algebras and operads to properads. We show that the simplicial
bar construction of a properad is quasi-isomorphic to its differential
bar construction.

\subsection{Definitions}

In every monoidal category $(\mathcal{A},\, \Box,\,
I)$, one can associate to every monoid $M$ a simplicial chain complex
on the objects $M^{\Box n}$,
for $n\in \mathbb{N}$. In the category of $k$-modules, one gets the classical
bar construction of associative algebras. In the case of monads, one
finds the triple bar construction of J. Beck \cite{Beck}. We apply this construction
to the monoidal category of $\Sy$-bimodules.

\subsubsection{\bf Simplicial bar construction with coefficients}

\begin{dei}[Simplicial bar construction with coefficients]

Let $\Po$ be a differential properad and
let $(L,\, l)$ and $(R,\, r)$ be two differential $\Po$-modules
on the right and on the left.
We denote
$$\mathcal{C}_n(L,\, \Po,\, R):=\Sigma^n L\boxtimes_c \underbrace{\Po\boxtimes_c \cdots
\boxtimes_c \Po}_n \boxtimes_c R.$$
We define the face maps $d_i\,
:\, \mathcal{C}_n(L,\,\Po,\, R) \to \mathcal{C}_{n-1}(L,\,\Po,\,
R)$ by \vspace{6pt}

\begin{tabular}{l}
$\bullet$ the right action $l$ if $i=0$, \\

$\bullet$ the composition $\mu$ of the $i^{\textrm{th}}$ with the ${(i+1)}^{\textrm{th}}$ line, composed
by elements of $\Po$,\\

$\bullet$ the left action $r$ if $i=n$.

\end{tabular}\\

The degeneracy maps $s_j \, : \, \mathcal{C}_n(L,\, \Po,\, R)
\to \mathcal{C}_{n+1}(L,\, \Po,\, R)$ are given by the insertion of
the unit $\eta$ of the properad : $L\boxtimes_c \Po^{\boxtimes_c j}
\boxtimes_c \eta \boxtimes_c
\Po^{\boxtimes_c n-j} \boxtimes_c R$.\\

The differential $d_\mathcal{C}$ on $\mathcal{C}(L,\, \Po,\, R)$
is defined by the sum of the simplicial differential with the canonical differentials
$$ d_\mathcal{C}:=\delta_L + \delta_\Po +\delta_R +\sum_{i=0}^n
(-1)^{i+1} d_i.$$
One can check that $d_\mathcal{C}^2=0$ thanks to the suspension $\Sigma^n$ in the definition
of the $\mathcal{C}_n(L,\, \Po,\, R)$.

This chain complex is called the \emph{simplicial bar construction with
coefficients in $L$ and $R$}.
\end{dei}

\begin{rem}
When the $\Sy$-modules $\Po$, $L$ and $R$ are concentrated in degree $0$,
this construction is strictly simplicial.
\end{rem}

Like the differential bar construction (\emph{cf.} Proposition~\ref{BarFoncteurExact}),
the simplicial bar construction is an exact functor.

\begin{pro}
\label{SimplicialBarFoncteurExact}
The simplicial bar construction $\Co(L,\, \Po,\, R)$ preserves quasi-isomorphisms.
\end{pro}

\begin{proo}
We use the same proof as for Proposition~\ref{BarFoncteurExact} with the
filtration $F_r:=\oplus_{s\leq r} \Co_{(s)}(L,\, \Po,\, R)$, where
$s$ is the number of vertices of the underlying level-graph indexed
by elements of $\Po$ (including the trivial
vertices indexed by $I$).
\end{proo}

We define the \emph{augmentation morphism} $\varepsilon \,  :\,
\mathcal{C}(L,\, \Po,\, R) \to L{\boxtimes_c}_\Po R$ by the canonical projection
$$\mathcal{C}_0(L,\, \Po,\, R)=L\boxtimes_c R \twoheadrightarrow L{\boxtimes_c}_\Po R. $$

As in the case of the differential bar construction, one has a notion of
\emph{reduced simplicial bar construction}.

\begin{dei}[Reduced simplicial bar construction]

The simplicial bar construction with trivial coefficients
$\mathcal{C}(I,\, \Po ,\, I)$ is called the \emph{reduced
simplicial bar construction}. We denote it by
$\bar{\mathcal{C}}(\Po)$.
\end{dei}

\begin{rem}
We have seen previously that the elements of the bar construction
$\B(L,\, \Po,\, R)$ can be represented by graphs and that the
coderivation $d_\theta$ corresponds to the notion of ``adjacent
pair composition''. The elements of the simplicial bar
construction can be represented by level graphs and the faces
$d_i$ correspond to compositions of two levels of operations.
\end{rem}

\subsubsection{\bf Augmented simplicial bar construction}

\begin{dei}[Augmented simplicial bar construction]
The chain complexes $\mathcal{C}(I,\, \Po,\, \Po)$ $=$
$\bar{\mathcal{C}}(\Po)\boxtimes_c \Po$
and $ \mathcal{C}(\Po,\,
\Po,\, I)=\Po \boxtimes_c \bar{\mathcal{C}}(\Po)$ are called
\emph{augmented simplicial bar construction on the right and on the left}.
\end{dei}

\begin{pro}
For every properad $\Po$ and every right $\Po$-module $L$, the augmentation morphism
$$\varepsilon \ : \ \mathcal{C}(L,\, \Po,\, \Po) \to L $$
is a quasi-isomorphism of right simplicial $\Po$-modules.
\end{pro}

We have the same result on the left for every left $\Po$-module $R$.

\begin{proo}
The proof is the same as in the case of associative algebras
(\emph{cf.} S\'eminaire Cartan \cite{Cartan}). We consider the following
filtration $$F_r:=\oplus_{s\leq r} \Co(L,\, \underbrace{\Po,\, \Po}_{s}),$$
 where $s$ denotes the
total homological degree of elements of $\Po$. This filtration is stable under the differential
$d_\Co$. It induces a spectral sequence $E^*_{s,\, t}$ such that
$E^0_{s,\, t}=\Co(L,\, \underbrace{\Po,\, \Po}_{s})_{s+t}$ and $d^0=\delta_L+d$.
We introduce
an extra degeneracy map $s_{n+1}\ :$
$$ \Sigma^n L \boxtimes_c \Po^{\boxtimes_c n} \boxtimes_c \Po \simeq \Sigma^n L
\boxtimes_c \Po^{\boxtimes_c n+1}\boxtimes_c I
\xrightarrow{\Sigma^n L\boxtimes_c \Po^{\boxtimes_c n+1} \boxtimes_c \eta}
\Sigma^{n+1} L \boxtimes_c \Po^{\boxtimes_c n+1} \boxtimes_c \Po, $$ which induces
a contracting homotopy on $E^0$. Once again, we conclude with the classical convergence theorem for
spectral sequences.
\end{proo}

\begin{cor}
The augmented simplicial bar construction on the left and on the right are
acyclic.
\end{cor}

\subsubsection{\bf Normalized bar construction}

To any simplicial chain complex, one can associate its normalized chain complex, which is
the quotient by the images of the degeneracy maps.

\begin{dei}[Normalized bar construction]

The \emph{normalized bar construction} is given by the quotient of the
simplicial bar construction by the images of the degeneracies maps.
We denote $\N_n(L,\, \Po,\, R)$ the following dg-$\Sy$-module
$$\Sigma^n \coker\left( L\boxtimes_c \Po^{\boxtimes_c (n-1)}
\boxtimes_c R \xrightarrow{\sum_{i=0}^{n-1} L\boxtimes_c
\Po^{\boxtimes_c i} \boxtimes_c \eta \boxtimes_c \Po^{\boxtimes_c
(n-i-1)} \boxtimes_c R} L \boxtimes_c \Po^{\boxtimes_c n}
\boxtimes_c R  \right). $$
The chain complex $\N(I,\,
\Po,\, I)$, denoted $\bar{\N}(\Po)$, is called the \emph{reduced normalized bar
construction}.
\end{dei}

\begin{rem}
The normalized bar construction $\N(L,\, \Po,\, R)$ is isomorphic to the sub-complex of the simplicial
bar construction $\Co(L,\, \Po,\, R)$ spanned by the level-graphs indexed by at least one
element of $\oPo$ on each level.
\end{rem}

\begin{pro}
\label{NormalizedFuncteurExact}
The normalized bar construction $\N(L,\, \Po,\, R)$ preserves quasi-isomorphisms.
\end{pro}

\begin{proo}
The proof is the same as Proposition~\ref{SimplicialBarFoncteurExact}.
\end{proo}

We compute the homology of the reduced normalized bar construction of a quasi-free
properad as we did for the bar construction in corollary~\ref{BarDeQuasiFree} .

\begin{pro}
\label{SimplicialBarDeQuasiFree}
Let $\Po=\F(M)$ be a quasi-free properad on a weight graded dg-$\Sy$-bimodule $M$ such that
$M^{(0)}=0$ and where the derivation $d_\theta$ is decomposable, $\theta \,  :\,  M \to \oplus_{s\ge 2}
\F_{(s)}(M)$.
The natural projection $\bar{\N}(\F(M))\to \Sigma  M$ is a quasi-isomorphism.
\end{pro}

\begin{proo}
On the sub-complex $\bar{\N}(\F(M))^{(\rho)}$, composed by elements of total
weight $(\rho)$, we consider the filtration $F_r:=\oplus_{s\ge -r} \bar{\N}_{(s)}(\F(M))^{(\rho)}$,
where $s$ denotes the number of vertices indexed by elements of $M$. The differential
$d_\N=d+\delta_M+d_\theta$ preserves this filtration, which induces a spectral sequence
$E^*_{s,\, t}$. We have $E^0_{s,\, t}=\bar{\N}_{(s)}(\F(M))^{(\rho)}_{s+t}$ and
$d^0=d+\delta_M$. Therefore, $E^0(\bar{\N}(\F(M))^{(\rho))}=\bar{\N}(E^0(\F(M)))^{(\rho)}$ and
the lemma is true if it is true for a weight graded free dg-properad.

If $\Po=\F(M)$ is a weight graded free dg-properad. Consider the filtration
$F'_r:=\oplus_{s\le r} \bar{\N}(\F(M))^s$, where $s$ denotes the sum of the homological degree
of the elements of $M$. One can see that the filtration is stable under the differential
$d_\N=d+\delta_M$. We have $E^0_{s,\, t}= \bar{\N}_t(\F(M))_{s+t}$ and $d^0=d$.
For every $s$, we define a contracting homotopy $h$ on $E^0_{s,\, *}$, which concludes the proof.
Let $\xi$ be an element of $\bar{\N}_t(\F(M))$ represented by a graph $g$ with $t$ levels.
Consider the operations of $\bar{\F}(M)$ indexing the last level and for each
of them, consider the operations of $M$ that have no operation below (such that
there are no outgoing edge attached to another operation). Call extremal such operations.
Every outgoing edge of an extremal operation is an outgoing edge of $g$.
Denote $m$ the extremal operation of $\xi$ that has among its outgoing edges the one indexed by the
smallest integer. We define the image of $\xi$ under $h$ by the level graph with $t+1$ levels
indexed by $m$ on the last level and by the other operations of $\xi$ on the first levels (except
the one without $m$).
\end{proo}

\subsection{Levelization morphism}

The \emph{levelization morphism} is an injective morphism
between the bar construction and the normalized bar construction which induces an isomorphism
in homology.

\subsubsection{\bf Definition}

\begin{figure}[h]
$$\xymatrix{& \ar[d] &  \ar[d]&  \\
 \ar@{--}[r] & *+[F-,]{\nu_1} \ar@{--}[r] \ar@/_1pc/[d] \ar[dr]   &
*+[F-,]{\nu_2} \ar@/^1pc/[d]  \ar[dl] | \hole  \ar@{--}[r] & \\
 \ar@{--}[r] & *+[F-,]{\mu_1}  \ar[d]  \ar@{--}[r] &
*+[F-,]{\mu_2}  \ar[d]  \ar@{--}[r] &  \\
 & & &  } $$
\caption{An example of $g_\xi$.}  \label{xi}
\end{figure}
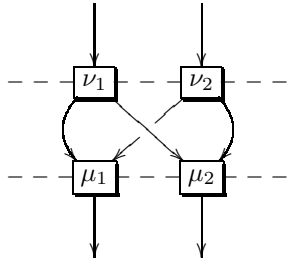

Let $\xi$ be an element of $\bar{\B}_{(n)}(\Po)=\F^c_{(n)}(\Sigma \oPo)$
represented by a graph $g_\xi$ with $n$ vertices (\emph{cf.} figure~\ref{xi}).

\begin{figure}[h]
$$\xymatrix{& \ar[dd] &  \ar[d]&  \\
 \ar@{--}[rr] &    &
*+[F-,]{\nu_2} \ar@/^1pc/[dd]  \ar[dddl] | \hole  \ar@{--}[r] & \\
 \ar@{--}[r] & *+[F-,]{\nu_1} \ar@/_1pc/[dd]    \ar[dr] \ar@{--}[rr] &
& \\
 \ar@{--}[rr] & &
*+[F-,]{\mu_2}  \ar[dd]  \ar@{--}[r] &  \\
 \ar@{--}[r] & *+[F-,]{\mu_1}  \ar[d]  \ar@{--}[rr] &
 &  \\
 & & &  } $$
\caption{An example of $g\in \mathcal{G}_4(\xi)$.} \label{echel}
\end{figure}
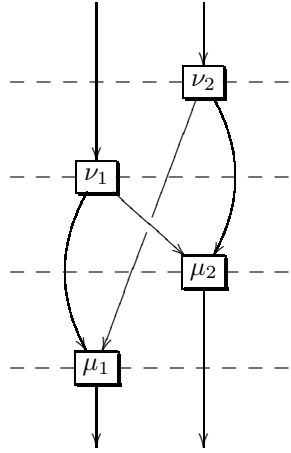

The element $\xi$ comes from a graph $g_r$ with $r$ levels
under the equivalence relation $\equiv$ (\emph{cf.} \cite{BVfreemonoid}).
Recall that the relation $\equiv$ allows to change the level of an operation,
up to sign. With this equivalence relation, one can represent $\xi$ with at least
one $n$-level-graph with $n$ non-trivial vertices (that's-to-say with one non-trivial
vertex by level, (\emph{cf.} figure~\ref{echel})). Denote,
$\mathcal{G}_n(\xi)$ the set of graphs with $n$ levels such that
there is only one non-trivial vertex on each level, which give $g_\xi$ under the
equivalence relation $\equiv$.

Let $g$ be a graph of $\mathcal{G}_n(\xi)$. The element $\xi$ is determined by $g$
and by the sequence of operations
$\Sigma p_1 \otimes \cdots \otimes \Sigma p_n$ which index the vertices of $g$.
To associate an element of $\bar{\N}(\Po)$ from $\xi$, one needs
to permutate the suspensions. This makes signs appear. We denote
$$g(\xi):=(-1)^{\sum_{i=1}^n(n-i)|p_i|} \Sigma^n g(p_1\otimes \cdots \otimes p_n),$$
where $g(p_1\otimes \cdots \otimes p_n)$ corresponds to the graph $g$
whose vertices are indexed by the operations $p_1,\, \ldots,\, p_n$.

\begin{dei}[Levelization morphism]
We define the \emph{levelization morphism} $$e \ : \
\bar{\B}(\Po) \to \bar{\N}(\Po)$$ by the following formula $e(\xi):=\sum_{g\in
\mathcal{G}_n(\xi)} g(\xi)$, for $\xi$ in
$\bar{\B}_{(n)}(\Po)$.
\end{dei}

\subsubsection{\bf Homological properties}

We prove that the levelization morphism is an injective
quasi-isomorphism of dg-$\Sy$-bimodules.

\begin{lem}
Let $\Po$ be a differential augmented properad.

The levelization morphism is an injective morphism of
dg-$\Sy$-bimodules
$$ e \ : \ \bar{\B}(\Po) \to \bar{\N}(\Po).$$
\end{lem}

\begin{proo}
Denote $d_\B$ the differential of $\bar{\B}(\Po)$. It is the sum of $2$ terms :
$ d_\B = \delta_\Po +d_\theta$.

Call $d_\N$ the differential on $\bar{\N}(\Po)$. It is the sum of $2$ terms.
The image of an element $\tau'=\Sigma^n \tau$ of
$\bar{\N}(\Po)$ under $d_\N$ is
$$ d_\N(\tau') = (-1)^n\Sigma^n \delta_\Po(\tau)  +
\sum_{i=1}^{n-1}(-1)^{i+1} \Sigma^{n-1} d_i(\tau).$$

Show that $d_\N\circ e(\xi)= e\circ d_\B(\xi)$.
\begin{itemize}
\item The commutativity of the canonical differentials
$\delta_\Po \circ e(\xi)= e\circ
\delta_\Po (\xi)$
comes from the good choice of signs and permutation of suspensions.

\item The coderivation $d_\theta$ composes pairs of adjacent operations
of $\Sigma \oPo$ of the graph representing $\xi$. And, $d_i\circ e$ corresponds to the composition
of $2$ levels of operations of $\oPo$. Two cases must be distinguished. Denote
$\xi=X\otimes \Sigma p \otimes \Sigma q \otimes Y$.

\begin{enumerate}
\item If the operations $\Sigma p$ and $\Sigma q$ are linked by at least of edge,
then the part of $\sum_{i=1}^{n-1}
(-1)^{i+1}d_i \circ e(\xi)$ involving the composition of
$p$ with $q$ is of the form $$\sum \varepsilon(-1)^{j+2}
\Sigma^{n-1} X'\otimes \mu(p\otimes q) \otimes Y',$$ where
$X'=p_1\otimes \cdots \otimes p_j$, $Y'=q_1\otimes \cdots \otimes
q_{n-j-2} \otimes \rho_1\otimes \cdots \otimes \rho_r$ and
$$\varepsilon= (-1)^{\sum_{i=1}^j
(n-i)|p_i|+(n-j-1)|p|+(n-j-2)|q|+\sum_{k=1}^{n-j-2}
(n-j-k-2)|q_k|}.$$ And the part of $d_\theta(\xi)$ involving
the composition of $\Sigma p$ with $\Sigma q$ is of the form
$$(-1)^{|X|+|p|} X\otimes \Sigma \mu(p\otimes q) \otimes Y.$$
The image of such an element under $e$ is
\begin{eqnarray*}
&& \sum (-1)^{|X|+|p|}\varepsilon(-1)^{|X'|+|p|}\Sigma^{n-1}X'\otimes
\mu(p\otimes q) \otimes Y'=\\
&& \sum \varepsilon (-1)^j \Sigma^{n-1}X'\otimes
\mu(p\otimes q) \otimes Y'.
\end{eqnarray*}

\item If the operations $\Sigma p$ and $\Sigma q$ are not linked
by an edge, then $d_\theta$ does not compose $\Sigma p$ with $\Sigma q$. And,
the part of $\sum_{i=0}^n  (-1)^{i+1} d_i\circ e(\xi)$ that comes from the composition
of the levels where lie $p$ and $q$ is the sum of the $2$ following terms :
\begin{eqnarray*}
\sum (-1)^j \Sigma^n d_{j+1} \left(\varepsilon X'\otimes p \otimes q \otimes Y' +
\varepsilon (-1)^{(|p|+1)(|q|+1)+|p|+|q|}
X'\otimes q \otimes p \otimes Y' \right)=\\
\sum (-1)^j \Sigma^n\left( \varepsilon X'\otimes p \otimes q \otimes Y'
+ \varepsilon (-1)^{(|p|+1)(|q|+1) +|p|+|q|+|p||q|} X'\otimes
p\otimes q \otimes Y' \right) =0.
\end{eqnarray*}
\end{enumerate}

\noindent We conclude by remarking that the image of $\xi$ under $\sum_{i=1}^{n-1} (-1)^{i+1}d_i \circ e - e
\circ d_\theta$ is a sum of terms of the form $(1)$ or $(2)$.
\end{itemize}

Show now the injectivity of $e$. Let $\xi$ be an
element of $\bar{\B}_{(n)}(\Po)=\F_{(n)}^c(\Sigma \oPo)$. By definition
of the cofree (connected) coproperad, one knows that
$\xi=\xi_1+\cdots +\xi_r$ where each $\xi_i$ is a finite sum
of elements of $\F_{(n)}^c(\Sigma \oPo)$ coming from indexing
the vertices of the same graph $g_{\xi_i}$. We introduce a morphism
$\pi \, : \, \bar{\N}_n(\Po) \to
\bar{\B}_{(n)}(\Po)$. To any element $\tau$ of $\bar{\N}_n(\Po)$
represented by a graph with $n$ levels, we
associate its image $\pi(\tau)$ under the equivalent relation $\equiv$.
Therefore, we have $\pi
\circ e(\xi)=n_1\xi_1+\cdots n_r\xi_r$ where the $n_i$ are integers $>0$.
The equation $e(\xi)=0$ gives
$n_1\xi_1+\cdots n_r\xi_r=0$. By identification of the underlying graphs, we have
$\xi_i=0$, for all $i$, hence $\xi=0$.
\end{proo}

The following theorem answers a question raised by M. Markl and A.A. Voronov
in \cite{MV} and generalizes
a result of B. Fresse for operads in \cite{Fresse}.

\begin{thm}\footnote{The proof of this theorem was given to the author by Benoit Fresse.}
\label{echelonnement}
Let $\Po$ be a weight graded augmented dg-properad.

The levelization morphism
$$ e \ : \ \bar{\B}(\Po) \to \bar{\N}(\Po) $$
is a quasi-isomorphism.
\end{thm}

\begin{proo}
Since the reduced bar construction (\emph{cf.}
Proposition~\ref{BarFoncteurExact}) and the reduced simplicial bar
construction (\emph{cf.}
Proposition~\ref{SimplicialBarFoncteurExact}) preserve
quasi-isomorphisms, it is enough to prove this theorem for a
resolution of $\Po$. The properad $\Po$ is weight graded and
augmented. Therefore, Theorem~\ref{barcobarresolution} shows that
the bar-cobar construction $\bar{\B^c}(\bar{\B}(\Po))$ of $\Po$ is
a quasi-free resolution of $\Po$. Denote
$M=\Sigma^{-1}\bar{\B}(\Po)$ and
$\Po'=\bar{\B^c}(\bar{\B}(\Po))=\F(M)$.
Proposition~\ref{BarDeQuasiFree} shows that
$\bar{\B}(\Po')=\bar{\B}(\F(M))$ is quasi-isomorphic to $\Sigma M$
and Proposition~\ref{SimplicialBarDeQuasiFree} shows that
$\bar{\N}(\Po')=\bar{\N}(\F(M))$ is quasi-isomorphic to $\Sigma
M$. We have the following commutative diagram of quasi-isomorphism
$$\xymatrix{ \bar{\B}(\Po) \ar[r]^{e} & \bar{\N}(\Po)\\
\bar{\B}(\F(M)) \ar[r]^{e} \ar[u]^{q.i.} \ar[d]_{q.i.} &  \bar{\N}(\F(M)) \ar[u]_{q.i.}
\ar[d]^{q.i.}\\
\Sigma M  \ar@{=}[r] & \Sigma M.} $$
\end{proo}

\begin{rem}
Using the same theorem for operads (\emph{cf.} \cite{Fresse}), the author has proved
in \cite{BV} a correspondence between operads and posets (partially ordered sets). To a
set operad $\Po$, one can associate a family of posets $\Pi_\Po$ of partition type.
The main theorem of \cite{BV} asserts that the operad is Koszul over $\mathbb{Z}$ if and only if
each poset of $\Pi_\Po$ is Cohen-Macaulay. It seems that the same kind of work can be
done in the case of properad with the previous theorem.
\end{rem}

\section{Koszul duality}
\label{KoszulDuality}

We deduce the Koszul duality theory of properads from the
comparison lemmas of the previous section.

\subsection{Koszul dual}
\label{constructionduale}
To a weight graded properad $\Po$, we define its \emph{Koszul dual}
as a sub-coproperad of its reduced bar construction.
Dually, to a weight graded coproperad $\mathcal{C}$, we define its Koszul dual
as a quotient properad of its reduced cobar construction.\\

\subsubsection{\bf Koszul dual of a properad}
\label{Koszuldualofproperad}
Let $\Po$ be a weight graded augmented dg-properad.
Recall that the reduced bar construction $\bar{\B}(\Po)=\F^c(\Sigma \oPo)$
of $\Po$ is bigraded by the number $(s)$ of non-trivial indexed vertices and by
the total weight $(\rho)$
$$\bar{\B}_{(s)}(\Po)=\bigoplus_{\rho \in \mathbb{N}}
\bar{\B}_{(s)}(\Po)^{(\rho)}.$$
Since the coderivation $d_\theta$, induced by the partial product, composes pairs of
adjacent vertices, we have
$$d_\theta\big( \bar{\B}_{(s)}(\Po)^{(\rho)}\big) \subset
\big( \bar{\B}_{(s-1)}(\Po)^{(\rho)}\big).$$

When $\Po$ is connected ($\Po^{(0)}=I$), the bar construction of $\Po$ has the following form

\begin{lem}
Let $\Po$ be a weight graded connected dg-properad.
We have the equalities
$$\left\{
\begin{array}{l}
\bar{\B}_{(\rho)}(\Po)^{(\rho)}=\F^c_{(\rho)}(\Sigma \oPo^{(1)}), \\
\bar{\B}_{(s)}(\Po)^{(\rho)}=0 \quad \textrm{if} \quad s>\rho.
\end{array} \right.$$
\end{lem}

\begin{rem}
The same result holds for the reduced cobar construction of a weight graded connected
dg-coproperad.
\end{rem}

\begin{dei}[Koszul dual of a properad]
Let $\Po$ be a weight graded connected dg-properad. We define the
\emph{Koszul dual of $\Po$} by
the weight graded dg-$\Sy$-bimodule
$${\Po^{\ac}}_{(\rho)}:=H_{(\rho)}\big( \bar{\B}_*(\Po)^{(\rho)},\,
d_\theta \big). $$
\end{dei}

The preceding lemma gives the form of the chain complex $\big(
\bar{\B}_*(\Po)^{(\rho)},\, d_\theta \big)$ :
$$\xymatrix{\cdots \ar[r]^{d_\theta} & 0 \ar[r]^(0.3){d_\theta} &
\bar{\B}_{(\rho)}(\Po)^{(\rho)} \ar[r]^(0.45){d_\theta} &
\bar{\B}_{(\rho-1)}(\Po)^{(\rho)} \ar[r]^(0.6){d_\theta}& \cdots \
.}$$

Hence, we have the formula
$${\Po^{\ac}}_{(\rho)}=\textrm{ker}\big(d_\theta
\, : \, \bar{\B}_{(\rho)}(\Po)^{(\rho)} \to
\bar{\B}_{(\rho-1)}(\Po)^{(\rho)} \big).$$

This formula shows that ${\Po^{\ac}}_{(\rho)}$ is a
weight graded dg-$\Sy$-bimodule, where the differential is induced
by the one of $\Po$ : $\delta_\Po$.

If the properad is concentrated in homological degree $0$, we have

$$\big(\bar{\B}_{(s)}(\Po)^{(\rho)}   \big)_d =\left\{
\begin{array}{ll}
\bar{\B}_{(s)}(\Po)^{(\rho)}  & \textrm{if} \quad d=s, \\
0 & \textrm{otherwise}.
\end{array} \right.$$

The dual coproperad is not necessarily concentrated in degree $0$. It verifies
$$\big( {\Po^{\ac}}_{(\rho)}\big)_d =\left\{\begin{array}{ll}
 {\Po^{\ac}}_{(\rho)}& \textrm{if} \quad d=\rho, \\
0 & \textrm{otherwise}.
\end{array} \right. $$

\subsubsection{\bf Koszul dual of a coproperad}
In the same way, we define the Koszul dual of a coproperad.

\begin{dei}[Koszul dual of a coproperad]
Let $\mathcal{C}$ be a weight graded connected dg-coproperad. We define the
 \emph{Koszul dual of $\mathcal{C}$} by the weight graded dg-$\Sy$-bimodule
$${\mathcal{C}^{\ac}}_{(\rho)}:=H_{(\rho)}\big( \bar{\B}^c_*
(\mathcal{C})^{(\rho)},\,d_{\theta'} \big). $$
\end{dei}

The cobar construction endowed with the derivation $d_{\theta'}$ is a chain complex
of the following form
$$\xymatrix{\cdots  \ar[r]^(0.4){d_{\theta'}} &
\bar{\B}_{(\rho-1)}(\mathcal{C})^{(\rho)}
\ar[r]^(0.55){d_{\theta'}} &
\bar{\B}_{(\rho)}(\mathcal{C})^{(\rho)}
\ar[r]^(0.65){d_{\theta'}}& 0
 \ .}$$

Therefore, the Koszul dual of a coproperad is isomorphic to

$${\mathcal{C}^{\ac}}_{(\rho)}=\textrm{coker}\big(d_{\theta'}
\, : \, \bar{\B}^c_{(\rho-1)}(\mathcal{C})^{(\rho)} \to
\bar{\B}^c_{(\rho)}(\mathcal{C})^{(\rho)} \big). $$

The module ${\mathcal{C}^{\ac}}_{(\rho)}$ is a weight graded
dg-$\Sy$-bimodule, where the differential is induced by the one of
$\mathcal{C}$ : $\delta_\mathcal{C}$.

\begin{pro}
Let $\Po$ be weight graded connected dg-properad.
The Koszul dual of $\Po$, denoted $\Po^{\ac}$, is a weight graded sub-dg-coproperad
of $\F^c(\Sigma \Po^{(1)})$.\\

Let $\mathcal{C}$ be a weight graded connected dg-coproperad.
The Koszul dual of $\mathcal{C}$, denoted
$\mathcal{C}^{\ac}$, is a weight graded dg-properad
 quotient of $\F(\Sigma^{-1} \mathcal{C}^{(1)})$.
\end{pro}

\begin{proo}
Since ${\Po^{\ac}}_{(\rho)}=\textrm{ker}\big(d_\theta \, : \,
\bar{\B}_{(\rho)}(\Po)^{(\rho)} \to
\bar{\B}_{(\rho-1)}(\Po)^{(\rho)} \big)$, we have that
${\Po^{\ac}}_{(\rho)}$ is a weight graded sub-dg-$\Sy$-bimodule
of $\F^c_{(\rho)}(\Sigma \Po^{(1)})$. It remains to show that
$\Po^{\ac}$ is stable under the coproduct
$\Delta$ of $\F^c(\Sigma \Po^{(1)})$. Let $\xi$ be an element of
${\Po^{\ac}}_{(\rho)}$, that's-to-say $\xi\in
\F_{(\rho)}^c(\Sigma\Po^{(1)})$ and $d_\theta(\xi)=0$. Denote
$$\Delta(\xi)=\sum_\Xi (\xi^1_1,\, \ldots ,\, \xi^1_{a_1})\, \sigma
\, (\xi^2_1,\, \ldots ,\, \xi^2_{a_2}), $$
where the sum is on a family  $\Xi$ of 2-level graphs, with $\xi_i^j \in
\F_{(s_{ij})}^c(\Sigma \Po^{(1)})$. The fact that $d_\theta$ is a
coderivation means that $d_\theta \circ \Delta (\xi)=
\Delta \circ d_\theta (\xi)$. Therefore, we have
\begin{eqnarray*}
d_\theta \circ \Delta (\xi)&=& \sum_\Xi  \Big( \sum_{k=1}^{a_1}
\pm (\xi^1_1,\, \ldots ,\, d_\theta(\xi^1_k),\, \ldots ,\,
\xi^1_{a_1})\, \sigma
\, (\xi^2_1,\, \ldots ,\, \xi^2_{a_2}) \\
&+&  \sum_{k=1}^{a_2}
   \pm (\xi^1_1,\, \ldots ,\,
\xi^1_{a_1})\, \sigma \, (\xi^2_1,\, \ldots ,\, d_\theta(\xi^2_k)
\, \ldots ,\,\xi^2_{a_2}) \Big)=0,
\end{eqnarray*}
where the $d_\theta(\xi^j_i)$ belong to $\F^c_{(s_{ij}-1)}
(\Sigma (\Po^{(1)}\oplus\Po^{(2)}))$ with only one vertex indexed by
$\Sigma \Po^{(2)}$. By identification, we have $d_\theta(\xi^j_i)=0$
for every $i,\,j$.\\

In the same way, one can see that the Koszul dual of $\mathcal{C}$ is
a quotient of
$\F(\Sigma^{-1}\mathcal{C}^{(1)})$. It remains to show that the product $\mu$
on $\F(\Sigma^{-1} \mathcal{C}^{(1)})$ defines a product on the quotient.
We consider the product
$$\mu \big( (c_1^1,\, \ldots ,\, d_{\theta'}(c),\, \ldots ,\,
c^1_{a_1}) \, \sigma \, (c^2_1,\, \ldots ,\, c^2_{a_2}) \big), $$
where the $c^j_i$ belong to $\F_{(s_{ij})}(\Sigma^{-1}
\mathcal{C}^{(1)})$ and  $c$ to $\F_{(s)}(\Sigma^{-1}
(\mathcal{C}^{(1)} \oplus \mathcal{C}^{(2)}))$ with one vertex indexed by
$\mathcal{C}^{(2)}$. Since the $c^j_i$ are elements of
$\F_{(s_{ij})}(\Sigma^{-1} \mathcal{C}^{(1)})$, we have
$d_{\theta'}(c^j_i)=0$. The fact that $d_{\theta'}$ is a
derivation gives
\begin{eqnarray*}
&&\mu \big( (c_1^1,\, \ldots ,\, d_{\theta'}(c),\, \ldots ,\,
c^1_{a_1}) \, \sigma \, (c^2_1,\, \ldots ,\, c^2_{a_2}) \big)\\
&=& d_{\theta'}\big( \mu \big(  (c_1^1,\, \ldots ,\, c,\, \ldots ,\,
c^1_{a_1}) \, \sigma \, (c^2_1,\, \ldots ,\, c^2_{a_2})  \big)
\big).
\end{eqnarray*} Therefore the product $\mu \big( (c_1^1,\, \ldots ,\, d_{\theta'}(c),\,
\ldots ,\, c^1_{a_1}) \, \sigma \, (c^2_1,\, \ldots ,\, c^2_{a_2})
\big) $ is null on the cokernel of $d_{\theta'}$.
\end{proo}

\subsubsection{\bf Quadratic properad}

We show here that the dual construction is a quadratic construction.\\

We recall from \ref{quadratic} that a quadratical properad is a properad of the form
$\F(V)/(R)$ where $R$ belongs to $\F_{(2)}(V)$. If we consider the weight
given by the number of indexed vertices, the ideal $(R)$ is homogenous
and the quadratic properad is weight graded.

Every quadratic properad is completely determined
by its graduation $\Po=\bigoplus_{\rho \in
\mathbb{N}} \Po_{(\rho)}$, the dg-$\Sy$-bimodule $\Po_{(1)}$ and the
dg-$\Sy$-bimodule $\Po_{(2)}$. One has $V=\Po_{(1)}$ and
$R=\textrm{ker}\big(\F_{(2)}(\Po_{(1)}) \xrightarrow{\mu}
\Po_{(2)} \big)$.

\begin{lem}
\label{dualequadratique}
Let $\mathcal{C}$ be a weight graded connected dg-coproperad.
The Koszul dual $\mathcal{C}^{\ac}$ is a quadratic properad generated by
$$ {\mathcal{C}^{\ac}}_{(1)}=\Sigma^{-1} \mathcal{C}^{(1)}
\quad \textrm{and} \quad  {\mathcal{C}^{\ac}}_{(2)}= \coker
\big(\theta' \, :\, \Sigma^{-1}\mathcal{C}^{(2)} \to
\F_{(2)}(\Sigma^{-1} \mathcal{C}^{(1)})\big).$$
\end{lem}

\begin{proo}
Denote $R=\textrm{ker}\big(\F_{(2)}(\mathcal{C}^{\ac}_{(1)}) \to
\mathcal{C}_{(2)}^{\ac}\big)=\textrm{Im}(\theta'_{\Sigma^{-1}\mathcal{C}^{(2)}})$.
By definition, $\mathcal{C}^{\ac}_{(\rho)}$ is the quotient of
$\F_{(\rho)}(\Sigma^{-1} \mathcal{C}^{(1)})$ by the image of
$d_{\theta'}$ on $\bar{\B}^c_{(\rho-1)}(\mathcal{C}^{(\rho)})$.
This image corresponds to graphs with $\rho$ vertices such that at least
one pair of adjacent vertices is the image of an element of
$\Sigma^{-1} \mathcal{C}^{(2)}$ under $\theta'$. This is isomorphic to the
part of weight $\rho$ of the ideal generated
by $R$. Therefore, we have $\mathcal{C}^{\ac}=\F(\Sigma^{-1}
\mathcal{C}^{(1)})/(R)$.
\end{proo}

\begin{rem}
Dually, one can define the notion of quadratic coproperad and prove
that the Koszul dual of a properad is a quadratic coproperad.
\end{rem}

\subsection{Koszul properads and Koszul resolutions}

We give here the notions of \emph{Koszul properad and coproperad}. When
a properad $\Po$ is a Koszul properad, it is quadratic, its dual is Koszul and
its bidual is isomorphic to $\Po$. We prove
that a properad $\Po$ is Koszul if and only if the reduced cobar construction of
$\Po^{\ac}$ is a resolution of $\Po$.

\subsubsection{\bf Definitions}

\begin{dei}[Koszul properad]
Let $\Po$ be a weight graded connected dg-properad.
The properad $\Po$ is a \emph{Koszul properad} if
the natural inclusion $\Po^{\ac}\hookrightarrow \bar{\B}(\Po)$ is a
quasi-isomorphism.
\end{dei}

In the same way, we have the definition of a \emph{Koszul coproperad}.

\begin{dei}[Koszul coproperad]
Let $\mathcal{C}$ be a weight graded connected dg-coproperad.
We call $\Co$ a \emph{Koszul coproperad} if the natural projection
$\bar{\B}^c(\mathcal{C}) \twoheadrightarrow \mathcal{C}^{\ac}$ is
a quasi-isomorphism.
\end{dei}

\begin{pro}
If $\Po$ is a weight graded Koszul properad, then
its Koszul dual $\Po^{\ac}$ is a Koszul coproperad and ${\Po^{\ac}}^{\ac}=\Po$.
\end{pro}

\begin{proo}
The properad $\Po$ is concentrated in (homological) degree $0$
($\delta_\Po=0$ and $\Po_0=\Po$). In this case, $\Po^{\ac}_{(\rho)}$
is a homogenous $\Sy$-bimodule of degree $\rho$ and
the elements of degree $0$ of
$\bar{\B}^c(\Po^{\ac})^{(\rho)}$ correspond to the elements of
$\bar{\B}^c_{(\rho)} ({\Po^{\ac}})^{(\rho)}$. This shows that
$$H_0\big(\bar{\B}^c_{(\rho)}(\Po^{\ac})^{(\rho)} \big)=H_\rho\big(
\bar{\B}^c_*(\Po^{\ac})^{(\rho)},\, d_{\theta'}
\big)={\Po^{\ac}}_{(\rho)}^{\ac}.$$

Since the properad $\Po$ is Koszul, the inclusion
$\Po^{\ac}\hookrightarrow \bar{\B}(\Po)$ is a
quasi-isomorphism. Using the comparison lemma
for quasi-free properads
(Theorem~\ref{lemmecomparaisonproperades}), we show that
the induced morphism $\bar{\B}^c(\Po^{\ac})\to
\bar{\B}^c(\bar{\B}(\Po))$ is a quasi-isomorphism.
Since the bar-cobar construction is a
resolution of $\Po$ (theorem~\ref{barcobarresolution}), the
cobar construction $\bar{\B}^c(\Po^{\ac})$ is quasi-isomorphic
to $\Po$. The $\Sy$-bimodule $\Po$ is homogenous of degree
$0$, therefore we have
$$ H_*\big(\bar{\B}^c(\Po^{\ac})^{(\rho)}\big) =
\left\{ \begin{array}{ll}
\Po^{(\rho)} & \textrm{if} \quad *=0,\\
0 & \textrm{otherwise}.
\end{array} \right.$$

Finally, we get the $\Po^{\ac}$ is a Koszul coproperad and that
${\Po^{\ac}}^{\ac}=\Po$.
\end{proo}

\begin{cor}
\label{Koszul2quadratic}
Let $\Po$ a weight graded connected properad. If
$\Po$ is a Koszul properad then $\Po$ is quadratic.
\end{cor}

\begin{proo}
If $\Po$ is a Koszul properad, by the preceding proposition, we know that
$\Po={\Po^{\ac}}^{\ac}$. And Lemma~\ref{dualequadratique}
shows that ${\Po^{\ac}}^{\ac}$ is quadratic.
\end{proo}

\subsubsection{\bf Koszul resolution and minimal model}

We have seen in the proof of the previous proposition that, when
$\Po$ is a Koszul properad, the quasi-isomorphism $\Po^{\ac}
\hookrightarrow \bar{\B}(\Po)$ induces a quasi-isomorphism
$\bar{\B}^c(\Po^{\ac}) \to \bar{\B}^c(\bar{\B}(\Po))$. When we compose
this quasi-isomorphism with the bar-cobar resolution $\zeta\, :\,
\bar{\B}^c(\bar{\B}(\Po)) \to \Po$ (\emph{cf.} Theorem~\ref{barcobarresolution}),
we get a quasi-isomorphism of weight graded dg-properads
$$ \bar{\B}^c(\Po^{\ac}) \to \Po.$$

\begin{dei}[Koszul resolution]
When $\Po$ is a Koszul properad, the resolution
$\bar{\B}^c(\Po^{\ac}) \to \Po$ is called a \emph{Koszul resolution}.
\end{dei}

Moreover, we have the following equivalence.

\begin{thm}[Koszul resolution]
\label{resolutionKoszul}
Let $\Po$ be a weight graded connected dg-properad.
The properad $\Po$ is a Koszul properad if and only if the morphism of
weight graded dg-properads
$\bar{\B}^c(\Po^{\ac}) \to \Po$ is a quasi-isomorphism.
\end{thm}

\begin{proo}
On the other way, if the morphism
$\bar{\B}^c (\Po^{\ac}) \to \Po $ is a quasi-isomorphism, then
$\bar{\B}^c(\Po^{\ac}) \to
 \bar{\B}^c(\bar{\B}(\Po))$ is also a quasi-isomorphism. Since we work with weight
 graded quasi-free properads, we can apply the comparison lemma for
 quasi-free properads (Theorem~\ref{lemmecomparaisonproperades}) which gives that the inclusion
 $\Po^{\ac}\hookrightarrow \bar{\B}(\Po)$ is a quasi-isomorphism.
\end{proo}

When $\Po$ is a Koszul properad concentrated in degree $0$, we have a
resolution $\bar{\B}^c(\Po^{\ac})=\F^c(\Sigma^{-1} \oPo^{\ac}) \to \Po$
build on a quasi-free properad and such that the differential is quadratic
$d_{\theta'}(\Sigma^{-1} \oPo^{\ac})\subset \F_{(2)}(\Sigma^{-1}
\oPo^{\ac})$. By analogy with associative algebras and operads, we call
this resolution the \emph{minimal model} of $\Po$.

\subsection{Koszul complex}

It is as difficult to show that
a properad $\Po$ is a Koszul properad from the definition ($\Po^{\ac}\to
\bar{\B}(\Po)$ quasi-isomorphism) as to show that the
cobar construction of $\Po^{\ac}$ is a resolution of $\Po$
($\bar{\B}^c(\Po^{\ac} \to \Po$ quasi-isomorphism). We introduce a smaller
chain complex whose acyclicity is a criterion which proves that
$\Po$ is a Koszul properad.

\subsubsection{\bf Koszul complex with coefficients}

Following the same method as in the definition of the bar construction
with coefficients (\emph{cf.} \ref{barconstructrioncoef}), we introduce here
the notion of \emph{Koszul complex with coefficients}.

\begin{dei}[Koszul complex with coefficients]

Let $\Po$ be a weight graded connected dg-properad
and let $L$ and $R$ be two modules (on the right and on the left)
on $\Po$. We call \emph{Koszul complex with
coefficients in the modules $L$ and $R$}, the chain complex defined
on the $\Sy$-bimodule $L\boxtimes_c\Po^{\ac}\boxtimes_c R$ by the
differential $d$, sum of the three following terms :
\vspace{6pt}

\begin{enumerate}
\item the canonical differentials $\delta_\Po$, $\delta_R$ and $\delta_L$
induced by the ones of $\Po$,$R$ and $L$.

\item the homogenous morphism $d_{\theta_L}$ of degree $-1$ that comes from
the structure of $\Po$-comodule on the left on $\Po^{\ac}$ :
$$ \theta_l\ : \ \Po^{\ac} \xrightarrow{\Delta}  \Po^{\ac}
\boxtimes_c \Po^{\ac} \twoheadrightarrow (I \oplus
\underbrace{\oPo^{\ac}}_1)  \boxtimes_c  \Po^{\ac} \to (I\oplus
\underbrace{\Po^{(1)}}_1)\boxtimes_c \Po^{\ac},$$

\item the homogenous morphism $d_{\theta_R}$ of degree $-1$ that comes from
the structure of $\Po$-comodule on the right on $\Po^{\ac}$ :
$$ \theta_r\ : \ \Po^{\ac} \xrightarrow{\Delta} \Po^{\ac}
\boxtimes_c \Po^{\ac} \twoheadrightarrow \Po^{\ac} \boxtimes_c (I
\oplus \underbrace{\oPo^{\ac}}_1) \to \Po^{\ac}\boxtimes_c
(I\oplus \underbrace{\Po^{(1)}}_1).$$
\end{enumerate}

We denote it $\mathcal{K}(L,\, \Po,\, R)$.
\end{dei}

Since $\Po^{\ac}$ can be injected into the bar construction
$\bar{\B}(\Po)$, the Koszul complex with coefficients is a
sub-$\Sy$-bimodule of the bar construction with coefficients.

\begin{pro}
Let $\Po$ be a weight graded connected dg-properad
and let $L$ and $R$ be two modules (on the right and on the left)
on $\Po$. The dg-$\Sy$-bimodule
$\mathcal{K}(L,\, \Po, \, R)=L\boxtimes_c \Po^{\ac} \boxtimes_c R$
is a sub-complex of the bar construction with coefficients
$\B(L,\, \Po,\, R)=L\boxtimes_c \bar{\B}(\Po)_c \boxtimes R$.
\end{pro}

\begin{proo}
This proposition relies on the fact that differential of the bar
construction is defined with the coproduct $\Delta$ of
$\bar{\B}(\Po)=\F^c(\Sigma \oPo)$ and that the differential of the
Koszul complex is also defined with the coproduct
on $\Po^{\ac}$. The Koszul dual $\Po^{\ac}$ is a sub-coproperad
of $\bar{\B}(\Po)$, which concludes the proof.
\end{proo}

The inclusion $\Po^{\ac}\hookrightarrow \bar{\B}(\Po)$ induces a
morphism of dg-$\Sy$-bimodules
$$L\boxtimes_c\Po^{\ac}\boxtimes_c R \hookrightarrow \B(L,\, \Po,\, R).$$

\subsubsection{\bf Koszul complex and minimal model}

We have studied a particular bar construction, namely the augmented
bar construction $\B(I,\, \Po,\,
\Po)=\bar{\B}(\Po)\boxtimes_c \Po$. We consider here the equivalent chain complex
in the framework of Koszul complexes.

\begin{dei}[Koszul complex]

We call \emph{Koszul complex} the chain complex
$\mathcal{K}(I,\, \Po,\, \Po)$ $=$ $\Po^{\ac}\boxtimes_c \Po$ (and
$\mathcal{K}(\Po,\, \Po,\, I)=\Po\boxtimes_c \Po^{\ac}$).
\end{dei}

The differential of this complex is defined by the sum of the morphism
$d_{\theta_r}$ with the
canonical differential $\delta_\Po$ induced by the one of
$\Po$. Remark that the morphism $d_{\theta_r}$ corresponds to
extract one operation $\nu\in\Po^{\ac}_{(1)}$ of $\Po^{\ac}$
from the top, identify this operation $\nu$ as an
operation of $\Po_{(1)}$ (since $\Po^{\ac}_{(1)}\cong\Po_{(1)}$)
and finally compose it in $\Po$. This can be graphically resumed by the
same kind of diagrams as the ones given by J.-L. Loday in \cite{Loday1}
for associative algebras.\\

The main theorem of this paper is the following criterion.

\begin{thm}[Koszul criterion]
\label{criteredeKoszul}
Let $\Po$ be a weight graded connected dg-properad.
The following assertions are equivalent.
\vspace{6pt}

\begin{tabular}{ll}
$(1)$ & The properad $\Po$ is a Koszul properad \\
 & (the inclusion $\Po^{\ac}\hookrightarrow
\bar{\B}(\Po)$
is a quasi-isomorphism).\\

$(2)$  &The Koszul complex $\Po^{\ac}\boxtimes_c \Po$ is acyclic.\\

$(2')$ &The Koszul complex $\Po \boxtimes_c \Po^{\ac}$ is acyclic.\\

$(3)$ &The morphism of weight graded dg-properads \\
& $\bar{\B}^c(\Po^{\ac})\twoheadrightarrow \Po $ is a quasi-isomorphism.
\end{tabular}
\end{thm}

\begin{proo}
We have already seen the equivalence $(1)\iff(3)$ (\emph{cf.}
Theorem~\ref{resolutionKoszul}).

By the comparison lemma for quasi-free $\Po$-modules (on the right),
the assertion $(1)$ is equivalent to the fact that the
morphism $\Po^{\ac}\boxtimes_c \Po \to \bar{\B}(\Po)\boxtimes_c
\Po$ is a quasi-isomorphism. The acyclicity of the augmented bar
construction (Theorem~\ref{acyclicitébaraugmentée}) shows
that the properad $\Po$ is Koszul $(1)$ if and only if
the Koszul complex $\Po^{\ac}\boxtimes_c \Po$ is acyclic $(2)$.

We use the comparison lemma for quasi-free $\Po$-modules on the left, to prove the
equivalence $(1)\iff(2')$.
\end{proo}

\subsection{Koszul duality for props}
\label{KDprops}

One can write the same theory for quadratic props. We are going to show how
to retrieve Koszul duality for quadratic props defined by connected relations
(\emph{cf.} \ref{quadratic}) from Koszul duality of properads.\\

Let $\widetilde{\Po}$ be quadratic prop $\widetilde{\F}(V)/(R)$
defined by connected relations $R\subset
\F_{(2)}(V)\subset\widetilde{\F}_{(2)}(V)$. Recall from
Proposition~\ref{quadraticconnectedprop} that
$\widetilde{\Po}=\widetilde{\F}(V)/(R)=\Ss(\F(V)/(R))$ is
isomorphic to the free prop on the properad $\Po=\F(V)/(R)$. With
the same methods, one can generalize the bar and cobar
constructions, the Koszul duals and complexes to props. Denote by
$\ \widetilde{}\ $ the constructions related to props. Since the
concatenation of elements of $\widetilde{\Po}$ is free, the
following diagram is commutative
$$\xymatrix{{\widetilde{\Po}}^{\ac} \ar@{^{(}->}[r]   &
\widetilde{\F}^c(\Sigma \bar{\widetilde{\Po}})\\
\Po^{\ac} \ar@{^{(}->}[r] \ar@{=}[u] & \F^c(\Sigma \bar{\Po}). \ar@{^{(}->}[u]}$$
In other words, the Koszul dual coprop $\widetilde{\Po}^{\ac}$ is a coproperad isomorphic
to the Koszul dual $\Po^{\ac}$ of $\Po$. One can also see that the Koszul complex
$\widetilde{\Po}^{\ac}\boxtimes \widetilde{\Po}$
associated to the prop $\widetilde{\Po}$
is isomorphic to $\Ss(\Po^{\ac} \boxtimes_c \Po)$. Hence, one is acyclic if and only if
the other is acyclic. The cobar construction $\widetilde{\F}(\Sigma^{-1}
\bar{\widetilde{\Po}^{\ac}})$ on
$\widetilde{\Po}^{\ac}$ is isomorphic to
$\Ss(\F(\Sigma^{-1} \bar{\Po^{\ac}}))=\Ss(\bar{\B}^c(\Po^{\ac})$. Since
$\widetilde{\Po}=\Ss(\Po)$, the cobar construction on $\widetilde{\Po}^{\ac}$ is a
resolution of $\widetilde{\Po}$ (morphism of props) if and only if the cobar construction on $\Po^{\ac}$
is resolution of $\Po$ (morphism of properads).

\begin{cor}
Let $\widetilde{\Po}$ be a quadratic connected dg-prop defined by connected relations.
The following assertions are equivalent.

\begin{enumerate}
\item The Koszul complex $\widetilde{\Po}^{\ac} \boxtimes \widetilde{\Po}$ is acyclic.

\item The morphism of weight graded dg-props
$\bar{\widetilde{\B}^c}(\widetilde{\Po}^{\ac}) \to \widetilde{\Po}$ is a quasi-isomorphism.
\end{enumerate}
\end{cor}

\begin{rem}
Since the author does not know any type of gebras defined on a
quadratic prop with non-confected relations, we are reluctant to
write Koszul duality for props in whole generality. The Koszul
resolutions for props are based on the free prop with is much
bigger than the free properad. For the applications of Koszul
duality studied in this paper (for instance notion of $\Po$-gebra
up to homotopy \emph{cf.} \ref{Pgebrahomotopy}), the resolution of
the underlying properad is enough and contains all the information
about the algebraic structure.
\end{rem}

\subsection{Relation with the Koszul duality theories
for associative algebras and for operads}

The notion of Koszul dual of a properad introduced here
is a coproperad (\emph{cf.}
\ref{constructionduale}). On the contrary, the Koszul dual of
an associative algebra is an algebra in \cite{Priddy} and
the Koszul dual of an operad is an operad in \cite{GK}.
To recover these classical constructions, we consider their
linear dual : the Czech dual of $\Po^{\ac}$.

\begin{dei}[Czech dual of an $\Sy$-bimodule]
Let $\Po$ be an $\Sy$-bimodule, we define the \emph{Czech dual
$\Po^\vee$ of $\Po$} by the $\Sy$-bimodule $\Po^\vee := \bigoplus_{\rho,\,
m,\, n} \Po_{(\rho)}^\vee(m,\, n)$, where
$$\Po_{(\rho)}^\vee (m,\,n):=sgn_{\Sy_m}\otimes_k \Po_{(\rho)}(m,\, n)^* \otimes_k
sgn_{\Sy_n}.$$
\end{dei}

The Czech dual is the linear dual twisted by the signature representations.

\begin{lem}
Let $(\Co,\, \Delta,\, \varepsilon)$ be a weight graded coproperad
such that the dimensions of the
modules $\Co_{(\rho)}(m,\,n)$ are finite over $k$,
for every $m$, $n$ and $\rho$. The $\Sy$-bimodule $\Co^\vee$
is naturally endowed with a structure of weight graded properad.
\end{lem}

\begin{proo}
Denote $\Po=\Co^\vee$. The coproduct $\Delta$ can be decomposed by the weight
$\Delta=\bigoplus_{\rho \in \mathbb{N}} \Delta_{(\rho)}$. To
define la multiplication $\mu$ on $\Co^\vee$, we dualize
$$\Delta_{(\rho)}(m,\,n) \, :\,
\Co_{(\rho)}(m,\,n)\to (\Co\boxtimes_c\Co)_{(\rho)}(m,\,n).$$

More precisely, we have
\begin{eqnarray*}
(\Po\boxtimes_c \Po)_{(\rho)}(m,\, n)&=&\bigoplus_{g\in
\mathcal{G}^2_c(m,\, n)} \left( \bigotimes_{\nu \, \textrm{vertex
of} \, g \atop \sum \rho_\nu =\rho}
\Co_{(\rho_\nu)}^\vee(|Out(\nu)|,\, |In(\nu)|)
\right)\Bigg/_\approx \\
&=&\bigoplus_{g\in \mathcal{G}^2_c(m,\, n)} \left( \left(
\bigotimes_{\nu \, \textrm{vertex of} \, g \atop \sum \rho_\nu
=\rho} \Co_{(\rho_\nu)}(|Out(\nu)|,\, |In(\nu)|)
\right)\Bigg/_\approx\right)^\vee,
\end{eqnarray*}
using the finite dimensional hypothesis of the $\Co_{(\rho)}(m,\,
n)$ and identifying the invariants with the coinvariants (we work
over a field $k$ of characteristic $0$). Hence, we define the composition
$\mu_{(\rho)}$ by the formula :
\begin{eqnarray*}
(\Po\boxtimes_c \Po)_{(\rho)}(m,\, n) &=&\bigoplus_{g\in
\mathcal{G}^2_c(m,\, n)} \left( \left( \bigotimes_{\nu \,
\textrm{vertex of} \, g \atop \sum \rho_\nu =\rho}
\Co_{(\rho_\nu)}(|Out(\nu)|,\, |In(\nu)|)
\right)\Bigg/_\approx\right)^\vee \\ &\longrightarrow & \left(
\bigoplus_{g\in \mathcal{G}^2_c(m,\, n)}  \left( \bigotimes_{\nu \,
\textrm{vertex of} \, g \atop \sum \rho_\nu =\rho}
\Co_{(\rho_\nu)}(|Out(\nu)|,\, |In(\nu)|)
\right)\Bigg/_\approx\right)^\vee\\
&=&(\Co\boxtimes_c\Co)_{(\rho)}^\vee(m,\,n) \xrightarrow{
^t\Delta_{(\rho)}} \Co^\vee_{(\rho)}(m,\, n)=\Po_{(\rho)}(m,\, n).
\end{eqnarray*}

The coassociativity of the coproduct $\Delta$ induces
the associativity of the product $\mu$. And the counit of
$\Co$ $\varepsilon \, :\, \Co \to I$ gives, after passage to the dual,
the unit of $\Po$ : $\varepsilon \, : \, I \to \Po$.
\end{proo}

\begin{pro}
\label{dimfinie}
Let $\Po$ be a weight graded properad (for instance a quadratic properad).
Denote $V=\Po^{(1)}$, the $\Sy$-bimodule generated by the elements of weight $1$
of $\Po$.

If there exists two integers $M$ and $N$ such that $V(m,\, n)=0$
for $m>M$ or $n>N$ and if the dimension of every module $V(m,\,
n)$ is finite on $k$, then the cobar construction
$\bar{\B}^c(\Po)$ on $\Po$ and the Koszul dual $\Po^{\ac}$ verify
the hypothesis of the preceding lemma.
\end{pro}

\begin{proo}
Since $\Po^{\ac}$ is a weight graded sub-coproperad of the cobar
construction $\bar{\B}^c(\Po)$, it is enough to show that the
dimension of the modules $\bar{\B}^c(\Po)_{(\rho)}(m,\,n)$ is
finite.

We have seen that
$\bar{\B}^c(\Po)_{(\rho)}(m,\,n)=\F^c_{(\rho)}(V)(m,\, n)$. In the
case where $V$ verifies the hypotheses of the proposition,
this module is given by a sum indexed by the set of graphs
with $\rho$ vertices and such that each vertex has at most
 $N$ inputs and $M$ outputs. Since this set is finite,
the dimension of the modules $V(m,\, n)$ is finite.
\end{proo}

\begin{cor}
\label{dualdimensionfinie}
For every quadratic properad
generated by an $\Sy$-bimodule $V$ such that the sum of the dimensions
$\sum_{m,\,n} \dim_k V(m,\,n)$ is finite, the Czech dual
${\Po^{\ac}}^\vee$ of the Koszul dual of $\Po$ is endowed with a natural
structure of properad .

Moreover, if $\Po=\F(V)/(R)$, then the
properad ${\Po^{\ac}}^\vee$ is quadratic and
${\Po^{\ac}}^\vee =\F(\Sigma V)/(\Sigma^2 R^\perp)$. We follow the classical notation
and denote this properad $\Po^!$.
\end{cor}

One can notice that the properad $\Po^!$ is
determined by $\Po^{\ac}_{(1)}=\Sigma V$ and
$\Po^{\ac}_{(2)}=\Sigma^2 R$.

Recall that in the case of quadratic algebra
$A=T(V)/(R)$, S. Priddy defined the dual algebra $A^!$ by
$T(V^*)/(R^\perp)$, when $V$ is finite dimensional. In the same way,
V. Ginzburg and M. M. Kapranov define the dual
of a quadratic operad $\Po=\F(V)/(R)$ by
$\Po^!=\F(V^\vee)/(R^\perp)$, when $V$
is finite dimensional.\\

In the finite dimensional case, the conceptual definitions of the
dual objects given here correspond to the classical definitions
 up to suspension  (\emph{cf.} \cite{BGS} and
\cite{Fresse}). In particular, the dual $A^{\ac}_{(n)}$ of an algebra
$A$ is isomorphic to $\Sigma^n {A^!}_n^*$ and the dual $\Po^{\ac}_{(n)}$ of an operad
$\Po$ is isomorphic to  $\Sigma^n{\Po^!}_n^\vee$. \\

When $\Po$ is an algebra, the Koszul complex and the Koszul resolution, coming from the cobar
construction,  introduced here correspond to the ones given by \cite{Priddy}.
When $\Po$ is an operad, they correspond to the ones given by \cite{GK}.

\section{Examples and $\Po$-gebra up to homotopy}

The last difficulty is to be able to prove that the Koszul complex
of a properad is acyclic. We consider a large class of properads
(properads defined by a \emph{replacement rule}) and we show that
they are Koszul properads. The examples of the properad $\BLi$ of
Lie bialgebras and the properad $\IBi$ of infinitesimal Hopf
bialgebra fall into this case. This method is a generalization to
properads of works of T. Fox, M. Markl \cite{FM}, M. Markl
\cite{MarklDistributive} and W. L. Gan \cite{Gan}. We define the
notion of \emph{$\Po$-gebra up to homotopy} and give the examples
of several $\Po$-gebras up to homotopy.

\subsection{Replacement rule}

Let $\Po$ be a quadratic properad of the form
$$\Po=\F(V,\, W)/(R\oplus D \oplus S),$$ where $R\subset \F_{(2)}(V)$, $S\subset
\F_{(2)}(W)$ and
$$D\subset (I\oplus \underbrace{W}_1)\boxtimes_c
(I\oplus \underbrace{V}_1)  \bigoplus (I\oplus
\underbrace{V}_1)\boxtimes_c (I\oplus \underbrace{W}_1).$$
The two pairs of $\Sy$-bimodules $(V,\, R )$ and $(W,\, S)$
generates two properads denoted $\Ao=\F(V)/(R)$ and
$\Bo=\F(W)/(S)$.

\begin{dei}[Replacement rule]
Let $\lambda$ be a morphism of $\Sy$-bimodules
$$\lambda \ : \   (I\oplus \underbrace{W}_1)\boxtimes_c
(I\oplus \underbrace{V}_1)  \to (I\oplus
\underbrace{V}_1)\boxtimes_c (I\oplus \underbrace{W}_1).$$
such that the $\Sy$-bimodule $D$ is defined by the image
$$(id,\, -\lambda) :\  (I\oplus \underbrace{W}_1)\boxtimes_c
(I\oplus \underbrace{V}_1) \to  (I\oplus
\underbrace{W}_1)\boxtimes_c (I\oplus \underbrace{V}_1)  \bigoplus
(I\oplus \underbrace{V}_1)\boxtimes_c (I\oplus
\underbrace{W}_1).$$
We call $\lambda$ a \emph{replacement rule}
and denote $D$ by $D_\lambda$ if the two following morphisms are injective
$$ \left\lbrace
\begin{array}{c}
 \underbrace{\Ao}_1 \boxtimes_c \underbrace{\Bo}_2 \to \Po  \\
\underbrace{\Ao}_2 \boxtimes_c \underbrace{\Bo}_1 \to \Po.
\end{array}
\right.$$
\end{dei}

\begin{rem}
The last condition must be seen as a coherence axiom. It ensures that the natural
morphism $\Ao \boxtimes \Bo \to \Po$ is injective (\emph{cf.} \cite{FM}).
\end{rem}

The examples treated here are of this form.
They come from the same kind of ``distributive law'' between two operads
described in \cite{FM} and in \cite{MarklDistributive}.

\begin{dei}[Reversed $\Sy$-bimodule]
Let $\Po$ be an $\Sy$-bimodule. We define its \emph{reversed $\Sy$-bimodule}
by
$$\Po^{op}(m,\, n):=\Po(n,\, m) .$$
\end{dei}

The action to take the reversed $\Sy$-bimodule $\Po^{op}$ corresponds to inverse the
direction of the operations of $\Po$. Let $\Po$ and $\Qo$ be two
$\Sy$-bimodules, we have
$$(\Po\boxtimes_c \Qo )^{op}=\Qo^{op} \boxtimes_c \Po^{op}.$$
When the $\Sy$-bimodule $\Po$ is endowed with a structure of
properad, the reversed $\Sy$-bimodule $\Po^{op}$
is also a properad.  \\

The properads treated are ($\BLi$ and $\IBi$) are generated by
operations ($n$ inputs and one output) and cooperations (one input
and $m$ outputs). One can write them $\F(V\oplus W)/(R\oplus
D_\lambda \oplus S)$, where $V$ represents the generating
operations ($V(m,\, n)=0$ for $m>1$) and where $W$ represents the
generating cooperations ($W(m,\, n)=0$ for $n>1$). The properad
$\Ao=\F(V)/(R)$ is an operad and the properad $\Bo^{op}=
\F(W^{op})/(S^{op})$ is an operad too. In this cases, the
replacement rule has the following form
$$\lambda \ : \   W\otimes_k V \to (I\oplus \underbrace{V}_1)\boxtimes_c
(I\oplus \underbrace{W}_1).$$

For the properad $\BLi$ associated to Lie bialgebras (\emph{cf.}
\ref{BiLie}), one has $\Ao=\Bo^{op}=\Li$ and the replacement rule
$\lambda$ is given by
$$\lambda \ : \
\vcenter{\xymatrix@M=0pt@R=6pt@C=4pt{1 & & 2 \\\ar@{-}[dr] &  &\ar@{-}[dl]  \\  &\ar@{-}[d] & \\
& \ar@{-}[dl]\ar@{-}[dr]& \\ & & \\ 1 & & 2}} \mapsto
 \vcenter{\xymatrix@M=0pt@R=6pt@C=2pt{ & 1 & & 2 \\ & \ar@{-}[d] &  &\ar@{-}[d]  \\  & \ar@{-}[dl]
 \ar@{-}[dr] &  & \ar@{-}[dl]\\
\ar@{-}[d]&  & \ar@{-}[d] & \\ & & & \\ 1 & & 2 &}} -
 \vcenter{\xymatrix@M=0pt@R=6pt@C=2pt{ & 2 & & 1 \\ & \ar@{-}[d] &  &\ar@{-}[d]  \\  & \ar@{-}[dl]
 \ar@{-}[dr] &  & \ar@{-}[dl]\\
\ar@{-}[d]&  & \ar@{-}[d] & \\ & & & \\ 1 & & 2 &}} +
\vcenter{\xymatrix@M=0pt@R=6pt@C=2pt{1 & & 2 & \\ \ar@{-}[d] &
&\ar@{-}[d] &  \\ \ar@{-}[dr] &
& \ar@{-}[dl] \ar@{-}[dr]& \\
& \ar@{-}[d]&  & \ar@{-}[d] \\ & & & \\ & 1 & & 2}}-
\vcenter{\xymatrix@M=0pt@R=6pt@C=2pt{2 & & 1 & \\ \ar@{-}[d] &
&\ar@{-}[d] &  \\ \ar@{-}[dr] &
& \ar@{-}[dl] \ar@{-}[dr]& \\
& \ar@{-}[d]&  & \ar@{-}[d] \\ & & & \\ & 1 & & 2}}.$$ In the case
of the properad $\IBi$ of infinitesimal Hopf bialgebra
 (\emph{cf.} \ref{InfBi}), the operads $\Ao$ and $\Bo^{op}$ are equal to the
 operad $\A$ of associative algebras and the replacement rule is given by
$$\lambda \ : \
\vcenter{\xymatrix@M=0pt@R=6pt@C=6pt{\ar@{-}[dr] &  &\ar@{-}[dl]  \\  &\ar@{-}[d] & \\
& \ar@{-}[dl]\ar@{-}[dr]& \\ & & }} \mapsto
 \vcenter{\xymatrix@M=0pt@R=6pt@C=6pt{ & \ar@{-}[d] &  &\ar@{-}[d]  \\  & \ar@{-}[dl]
 \ar@{-}[dr] &  & \ar@{-}[dl]\\
\ar@{-}[d]&  & \ar@{-}[d] & \\ & & & }}  +
\vcenter{\xymatrix@M=0pt@R=6pt@C=6pt{ \ar@{-}[d] & &\ar@{-}[d] &
\\ \ar@{-}[dr] &
& \ar@{-}[dl] \ar@{-}[dr]& \\
& \ar@{-}[d]&  & \ar@{-}[d] \\ & & & }} \ .$$

A replacement rule allows to permutate vertically operations and cooperations.

\begin{lem}
\label{FormeProperade}
Every properad of the form $\Po=\F(V,\,
W)/(R\oplus D_\lambda \oplus S)$ defined by a replacement rule
 $\lambda$ is isomorphic to the $\Sy$-bimodule
$$\Po \cong A\boxtimes_c B .$$
\end{lem}

\begin{proo}
As we have seen in the previous remark, the last condition defining a replacement rule
ensures that the natural morphism  $A\boxtimes_c B \to \Po$ is injective.

To show that it is also surjective, we choose a representative in $\F(V\oplus W)$
of every element of $\Po$. It can be written as a finite sum of graphs indexed by operations
of $V$ and $W$. For every graph, we fix each vertex on a level and we permutate
vertically elements of $V$ and elements of $W$ with the replacement rule to have finally
all the elements of $V$ under the ones of $W$. Therefore, we get an element of
$\F(V)\boxtimes_c\F(W)$ which projected in $A\boxtimes_c B$ gives the wanted element.
\end{proo}

In the two previous examples, this lemma shows that we can write
every element of $\Po$ as a sum of elements of $A\boxtimes_c B$.
This corresponds to put the cooperations of $\Bo$ at the top and
the operations of $\Ao$ under. In the case of Lie bialgebras, this
result has yet been proved be B. Enriquez and P. Etingov in
\cite{EE} (section 6.4).

\subsection{Koszul dual of a properad defined by a replacement rule}

\begin{pro}
\label{FormeProperadeDuale} Let $\Po$ be a properad of the form
$\Po=\F(V,\, W)/(R\oplus D_\lambda \oplus S)$ defined by a
replacement rule
 $\lambda$ and such that the sum of the dimensions of $V$ and
 $W$ on $k$, $\sum_{m,\, n} \dim_k (V\oplus
W)(m,\, n)$, is finite.

The dual properad is given by
$$\Po^!=\F( \Sigma W^\vee \oplus \Sigma V^\vee)/(\Sigma^2 S^\perp
 \oplus \Sigma^2
D_{ ^t\lambda} \oplus  \Sigma^2 R^\perp),$$ where the replacement rule is
$^t\lambda$. And the dual coproperad is isomorphic to the $\Sy$-bimodules
$$\Po^{\ac}\cong B^{\ac}\boxtimes_c A^{\ac}.$$
\end{pro}

\begin{proo}
Corollary~\ref{dualdimensionfinie} shows that ${\Po^!}^\vee
\cong \Po^{\ac}$.

Denote $R^\perp$ the orthogonal of $R$ in $\F_{(2)}(V^\vee)$,
$S^\perp$ the orthogonal of $S$ in $\F_{(2)}(W^\vee)$ and
$D_\lambda^\perp$ the orthogonal of $D_\lambda$ in $(I\oplus
\underbrace{W^\vee}_1)\boxtimes_c (I\oplus \underbrace{V^\vee}_1)
\bigoplus (I\oplus \underbrace{V^\vee}_1)\boxtimes_c (I\oplus
\underbrace{W^\vee}_1)$. Therefore, the dual properad $\Po^!$ is
given by
$$\Po^!=\F(\Sigma V^\vee \oplus \Sigma W^\vee)/(\Sigma^2 R^\perp \oplus \Sigma^2
D_\lambda^\perp \oplus \Sigma^2 S^\perp).$$ It remains to see that
$\Sy$-bimodule $D_\lambda^\perp$ is isomorphic to the image of the morphism
$$(id,- {^t\lambda}) \ : \  (I\oplus \underbrace{V^\vee}_1)\boxtimes_c
(I\oplus \underbrace{W^\vee}_1) \to (I\oplus
\underbrace{V^\vee}_1)\boxtimes_c (I\oplus \underbrace{W^\vee}_1)
\bigoplus (I\oplus \underbrace{W^\vee}_1)\boxtimes_c (I\oplus
\underbrace{V^\vee}_1),$$ which means that
$D_\lambda^\perp=D_{^t\lambda}$.

Applying the previous lemma to the properad $\Po^!$, we get
$\Po^!\cong B^!\boxtimes_c A^!$ and $\Po^{\ac}\cong B^{\ac}
\boxtimes_c A^{\ac}$ by Czech duality.
\end{proo}

The two examples given by the properads $\BLi$ and $\IBi$ of Lie
bialgebras and infinitesimal Hopf bialgebras verify the hypotheses
of the proposition. They are generated by a finite number of
operations and cooperations.

Since the composition of the form
$\vcenter{\xymatrix@M=0pt@R=6pt@C=6pt{& \ar@{-}[d] &   \\
&\ar@{-}[dl]
\ar@{-}[dr] & \\
\ar@{-}[dr]& & \ar@{-}[dl]\\ &\ar@{-}[d] & \\ & & }}$
does not appear in their definition, we must have it in the relations
of the Koszul duals.

\begin{cor}
$ $
\begin{enumerate}
\item The Koszul dual properad $\BLi^!$ of the properad of Lie bialgebras
is given by
$$\BLi^!=\F(V)/(R),$$
where $V=\vcenter{\xymatrix@M=0pt@R=6pt@C=4pt{{\scriptstyle 1} & &
{\scriptstyle 2}
\\\ar@{-}[dr] &
&\ar@{-}[dl]  \\  &\ar@{-}[d] & \\
& & \\ & & }} \oplus \vcenter{\xymatrix@M=0pt@R=6pt@C=4pt{  & \ar@{-}[d]&  \\
&\ar@{-}[dl] \ar@{-}[dr]& \\ & & \\{\scriptstyle 1}&
&{\scriptstyle 2} }}$, that's-to-say a commutative operation and a cocommutative cooperation,
and
\begin{eqnarray*}
R&=&\left\{ \begin{array}{l} \left(
\vcenter{\xymatrix@M=0pt@R=6pt@C=4pt{{\scriptstyle 1} & &
{\scriptstyle 2}& & {\scriptstyle 3}\\
\ar@{-}[dr] &
&\ar@{-}[dl] & & \ar@{-}[dl]  \\
& \ar@{-}[dr] & &\ar@{-}[dl]  & \\
& &\ar@{-}[d] & & \\
& & \\ & & }} - \vcenter{\xymatrix@M=0pt@R=6pt@C=4pt{{\scriptstyle
1} & &
{\scriptstyle 2}& & {\scriptstyle 3}\\
\ar@{-}[dr] &
&\ar@{-}[dr] & & \ar@{-}[dl]  \\
& \ar@{-}[dr] & &\ar@{-}[dl]  & \\
& &\ar@{-}[d] & & \\
& & \\ & & }} \right).k[\Sy_3]
 \\
\oplus \ k[\Sy_3].\left( \vcenter{\xymatrix@M=0pt@R=6pt@C=4pt{ & & & & \\
 & &\ar@{-}[d] & & \\
& & \ar@{-}[dl]  \ar@{-}[dr] & & \\
& \ar@{-}[dl]  \ar@{-}[dr] & &\ar@{-}[dr] & \\
& & & &\\ {\scriptstyle 1} & & {\scriptstyle 2}& &{\scriptstyle 3}
}} -
\vcenter{\xymatrix@M=0pt@R=6pt@C=4pt{ & & & & \\
 & &\ar@{-}[d] & & \\
& & \ar@{-}[dl]  \ar@{-}[dr] & & \\
& \ar@{-}[dl]   & &\ar@{-}[dl] \ar@{-}[dr] & \\
& & & &\\ {\scriptstyle 1} & & {\scriptstyle 2}& &{\scriptstyle 3}
}} \right)
\\ \ \\
\oplus \ \vcenter{\xymatrix@M=0pt@R=6pt@C=4pt{{\scriptstyle 1} & &
{\scriptstyle 2}
\\\ar@{-}[dr] &  &\ar@{-}[dl]  \\  &\ar@{-}[d] & \\
& \ar@{-}[dl]\ar@{-}[dr]& \\ & & \\ {\scriptstyle 1} & &
{\scriptstyle 2}}} -
 \vcenter{\xymatrix@M=0pt@R=6pt@C=2pt{ & {\scriptstyle 1} & &{\scriptstyle 2}
 \\ & \ar@{-}[d] &  &\ar@{-}[d]  \\  & \ar@{-}[dl]
 \ar@{-}[dr] &  & \ar@{-}[dl]\\
\ar@{-}[d]&  & \ar@{-}[d] & \\ & & & \\ {\scriptstyle 1} & &
{\scriptstyle 2} &}} \oplus
\vcenter{\xymatrix@M=0pt@R=6pt@C=4pt{{\scriptstyle 1} & &
{\scriptstyle 2}
\\\ar@{-}[dr] &  &\ar@{-}[dl]  \\  &\ar@{-}[d] & \\
& \ar@{-}[dl]\ar@{-}[dr]& \\ & & \\ {\scriptstyle 1} & &
{\scriptstyle 2}}} -
 \vcenter{\xymatrix@M=0pt@R=6pt@C=2pt{ & {\scriptstyle 2} & & {\scriptstyle 1}
  \\ & \ar@{-}[d] &  &\ar@{-}[d]  \\  & \ar@{-}[dl]
 \ar@{-}[dr] &  & \ar@{-}[dl]\\
\ar@{-}[d]&  & \ar@{-}[d] & \\ & & & \\ {\scriptstyle 1} & &
{\scriptstyle 2} &}} \oplus
\vcenter{\xymatrix@M=0pt@R=6pt@C=4pt{{\scriptstyle 1} & &
{\scriptstyle 2}
\\\ar@{-}[dr] &  &\ar@{-}[dl]  \\  &\ar@{-}[d] & \\
& \ar@{-}[dl]\ar@{-}[dr]& \\ & & \\ {\scriptstyle 1} & &
{\scriptstyle 2}}} -
\vcenter{\xymatrix@M=0pt@R=6pt@C=2pt{{\scriptstyle 1} & & {\scriptstyle 2} & \\
\ar@{-}[d] & &\ar@{-}[d] &  \\ \ar@{-}[dr] &
& \ar@{-}[dl] \ar@{-}[dr]& \\
& \ar@{-}[d]&  & \ar@{-}[d] \\ & & & \\ &{\scriptstyle 1} &
&{\scriptstyle 2}}} \oplus
\vcenter{\xymatrix@M=0pt@R=6pt@C=4pt{{\scriptstyle 1} & &
{\scriptstyle 2}
\\\ar@{-}[dr] &  &\ar@{-}[dl]  \\  &\ar@{-}[d] & \\
& \ar@{-}[dl]\ar@{-}[dr]& \\ & & \\ {\scriptstyle 1} & &
{\scriptstyle 2}}} -
\vcenter{\xymatrix@M=0pt@R=6pt@C=2pt{{\scriptstyle 2} & &
{\scriptstyle 1} &
\\ \ar@{-}[d] & &\ar@{-}[d] &  \\ \ar@{-}[dr] &
& \ar@{-}[dl] \ar@{-}[dr]& \\
& \ar@{-}[d]&  & \ar@{-}[d] \\ & & & \\ & {\scriptstyle 1} & &
{\scriptstyle 2}}} \\ \ \\
\oplus \ \vcenter{\xymatrix@M=0pt@R=7pt@C=7pt{& \ar@{-}[d] &   \\
&\ar@{-}[dl]
\ar@{-}[dr] & \\
\ar@{-}[dr]& & \ar@{-}[dl]\\ &\ar@{-}[d] & \\ & & }} \ .
\end{array}
\right. \end{eqnarray*}

It corresponds to the dioperad of non-unitary Frobenius algebras.
As an $\Sy$-bimodule, it is isomorphic to
$\BLi^!(m,\, n)=k,$ with trivial actions of $\Sy_m$ and $\Sy_n$.

\item The Koszul dual properad $\IBi^!$ of the properad of
infinitesimal Hopf algebras is given by
$$\IBi^!=\F(V)/(R),$$
o\`u $V=\vcenter{\xymatrix@M=0pt@R=6pt@C=6pt{\ar@{-}[dr] &
&\ar@{-}[dl]  \\  &\ar@{-}[d] & \\
& & \\ & & }} \oplus \vcenter{\xymatrix@M=0pt@R=6pt@C=6pt{  & \ar@{-}[d]&  \\
&\ar@{-}[dl] \ar@{-}[dr]& \\ & &  \\ & &}}$  et
\begin{eqnarray*}
R&=&\left\{ \begin{array}{l} \vcenter{\xymatrix@M=0pt@R=6pt@C=6pt{
\ar@{-}[dr] &
&\ar@{-}[dl] & & \ar@{-}[dl]  \\
& \ar@{-}[dr] & &\ar@{-}[dl]  & \\
& &\ar@{-}[d] & & \\
& & \\ & & }} - \vcenter{\xymatrix@M=0pt@R=6pt@C=6pt{
 \ar@{-}[dr] &
&\ar@{-}[dr] & & \ar@{-}[dl]  \\
& \ar@{-}[dr] & &\ar@{-}[dl]  & \\
& &\ar@{-}[d] & & \\
& & \\ & & }}
 \\
\oplus \  \vcenter{\xymatrix@M=0pt@R=6pt@C=6pt{ & & & & \\
 & &\ar@{-}[d] & & \\
& & \ar@{-}[dl]  \ar@{-}[dr] & & \\
& \ar@{-}[dl]  \ar@{-}[dr] & &\ar@{-}[dr] & \\
& & & & }} -
\vcenter{\xymatrix@M=0pt@R=6pt@C=6pt{ & & & & \\
 & &\ar@{-}[d] & & \\
& & \ar@{-}[dl]  \ar@{-}[dr] & & \\
& \ar@{-}[dl]   & &\ar@{-}[dl] \ar@{-}[dr] & \\
& & & & }}
\\ \ \\
\oplus \  \vcenter{\xymatrix@M=0pt@R=6pt@C=6pt{\ar@{-}[dr] &  &\ar@{-}[dl]  \\  &\ar@{-}[d] & \\
& \ar@{-}[dl]\ar@{-}[dr]& \\ & & }} -
 \vcenter{\xymatrix@M=0pt@R=6pt@C=6pt{ & \ar@{-}[d] &  &\ar@{-}[d]  \\  & \ar@{-}[dl]
 \ar@{-}[dr] &  & \ar@{-}[dl]\\
\ar@{-}[d]&  & \ar@{-}[d] & \\ & & & }}  \oplus
\vcenter{\xymatrix@M=0pt@R=6pt@C=6pt{\ar@{-}[dr] &  &\ar@{-}[dl]  \\  &\ar@{-}[d] & \\
& \ar@{-}[dl]\ar@{-}[dr]& \\ & & }} -
\vcenter{\xymatrix@M=0pt@R=6pt@C=6pt{\ar@{-}[d] & &\ar@{-}[d] &
\\ \ar@{-}[dr] &
& \ar@{-}[dl] \ar@{-}[dr]& \\
& \ar@{-}[d]&  & \ar@{-}[d] \\ & & & }} \\ \ \\
\oplus \ \vcenter{\xymatrix@M=0pt@R=7pt@C=7pt{& \ar@{-}[d] &   \\
&\ar@{-}[dl]
\ar@{-}[dr] & \\
\ar@{-}[dr]& & \ar@{-}[dl]\\ &\ar@{-}[d] & \\ & & }} \

\oplus
\vcenter{\xymatrix@R=1pt@C=1pt{ &\ar@{-}[d] & & &\ar@{-}[d] \\
 & *=0{}\ar@{-}[dl]\ar@{-}[dr] & & & *=0{}\ar@{-}[d] \\
 *=0{}\ar@{-}[ddrr] & & *=0{}\ar@{-}[ddll] | \hole  & & *=0{}\ar@{-}[d] \\
 & & & & *=0{}\ar@{-}[d] \\
 *=0{}\ar@{-}[d] & & *=0{}\ar@{-}[dr] & & *=0{}\ar@{-}[dl] \\
 *=0{}\ar@{-}[d] & & & *=0{}\ar@{-}[d] & \\
 & & & & }}
\oplus
\vcenter{\xymatrix@R=1pt@C=1pt{ *=0{}\ar@{-}[dddd]  & & & *=0{}\ar@{-}[d]  & \\
 & & & *=0{}\ar@{-}[dr] \ar@{-}[dl]  & \\
 & & *=0{}\ar@{-}[ddrr]  & & *=0{}\ar@{-}[ddll] | \hole  \\
 & & & & \\
 *=0{}\ar@{-}[dr] & & *=0{}\ar@{-}[dl] & & *=0{}\ar@{-}[d]  \\
 &  *=0{}\ar@{-}[d]& & & *=0{}\ar@{-}[d]  \\
 & & & &  }}
\oplus
\vcenter{\xymatrix@R=1pt@C=1pt{ *=0{}\ar@{-}[dd] &  &  &  *=0{}\ar@{-}[d] &   \\
 &  &  &  *=0{}\ar@{-}[dr]\ar@{-}[dl] &   \\
 *=0{}\ar@{-}[ddrr] &  & *=0{}\ar@{-}[ddll] | \hole  &  &   *=0{}\ar@{-}[dddd] \\
 &  &  &  &   \\
 *=0{}\ar@{-}[dr] &  &  *=0{}\ar@{-}[dl] &  &   \\
 &  *=0{}\ar@{-}[d] &  &  &   \\
 &  &  &  &   }}
\oplus
\vcenter{\xymatrix@R=1pt@C=1pt{ &  *=0{}\ar@{-}[d] &  & &   *=0{}\ar@{-}[dd] \\
 &  *=0{}\ar@{-}[dr]\ar@{-}[dl] &  &  &   \\
 *=0{}\ar@{-}[dddd] &  & *=0{}\ar@{-}[ddrr]  &  &  *=0{}\ar@{-}[ddll] | \hole  \\
 &  &  &  &   \\
 &  &  *=0{}\ar@{-}[dr] &  & *=0{}\ar@{-}[dl]   \\
 &  &  &  *=0{}\ar@{-}[d] &   \\
 &  &  &  &   }}.
\end{array}
\right. \end{eqnarray*}

It is isomorphic, as an $\Sy$-bimodule, to
$$\IBi^!(m,\, n)=k[\Sy_m]\otimes_k k[\Sy_n].$$
\end{enumerate}
\end{cor}

\subsection{Koszul complex of a properad defined by a replacement rule}

\begin{pro}
Let $\Po$ be a properad of the form $\Po=\F(V,\, W)/(R\oplus
D_\lambda \oplus S)$ defined by a distributive law $\lambda$ and
such that the sum of the dimensions of $V$ and $W$ on $k$,
$\sum_{m,\, n} \dim_k (V\oplus W)(m,\, n)$, is finite. We define
the two properads $\Ao$ and $\Bo$ by $\Ao=\F(V)/(R)$ and
$\Bo=\F(W)/(S)$. We suppose that $\Bo$ is a properad concentrated
in homological degree $0$. If $\Ao$ and $\Bo$ are Koszul
properads, then $\Po$ is a Koszul properad too.
\end{pro}

\begin{proo}
Lemma~\ref{FormeProperade} and Proposition~\ref{FormeProperadeDuale} show that
the Koszul complex of $\Po$ is of the following form
$$\Po^{\ac}\boxtimes_c \Po =(\Bo^{\ac}\boxtimes_c \Ao^{\ac})\boxtimes_c(\Ao\boxtimes_c \Bo)=
\Bo^{\ac}\boxtimes_c (\Ao^{\ac}\boxtimes_c \Ao )\boxtimes_c \Bo.$$
We introduce the following filtration of the Koszul complex
$F_n(\Po^{\ac}\boxtimes_c\Po)$ which corresponds to the sub-$\Sy$-bimodule
of $\Bo^{\ac}\boxtimes_c (\Ao^{\ac}\boxtimes_c \Ao )\boxtimes_c \Bo$
generated by 4-levels graphs where the sum of the weight of the
elements of $\Bo^{\ac}$
on the fourth level is less then $n$.
This filtration is stable under the differential of
the Koszul complex. Therefore, it induces a spectral sequence denoted
$E^*_{p,\, q}$. The first term of this spectral sequence
$E^0_{p,\, q}$ is isomorphic to the sub-$\Sy$-bimodule of
$\Bo^{\ac}\boxtimes_c (\Ao^{\ac}\boxtimes_c \Ao
)\boxtimes_c \Bo$ composed by elements of homological degree $p+q$ and that
such that the sum of the weight of the elements of $\Bo^{\ac}$ on the
fourth level is equal to $p$. Since $\Bo$ is concentrated in homological degree
$0$, we have $E^0_{p,\, q}=\underbrace{\Bo^{\ac}}_{p} \boxtimes_c \underbrace{(\Ao^{\ac}
\boxtimes_c \Ao)}_{q} \boxtimes_c \Bo$.
The differential $d_0$ is the differential of the Koszul complex of $\Ao$.
Since $\Ao$ is a Koszul properad, we get
$E^1_{p,\,q}=0$ if $q\ne 0$ and $E^1_{p,\,0}=\underbrace{B^{\ac}}_{p}\boxtimes_c B$.
the differential $d_1$ is the differential of the Koszul complex of $\Bo$.
Once again, since $\Bo$ is a Koszul properad, we have
$$E^2_{p,\,q}=\left\{
\begin{array}{ll}
I & \textrm{if} \ p=q=0, \\
0 & \textrm{otherwise}.
\end{array}\right. $$
The filtration is exhaustive and bounded below. Therefore, we can apply the classical
theorem of convergence of spectral sequences(\emph{cf.} \cite{Weibel}
5.5.1) which gives that the spectral sequences converges to the homology of the
Koszul complex of $\Po$. This complex is acyclic and the properad $\Po$ is Koszul.
\end{proo}

\begin{cor}
The properads $\BLi$ of Lie bialgebras and $\IBi$ of infinitesimal
Hopf algebras.
\end{cor}

\begin{proo}
In the case of $\BLi$, the properad $\Ao$ is the Koszul operad of
Lie algebras $\Li$ and the properad $\Bo$ is the reversed properad of the Koszul
operad $\Li$, $\Bo=\Li^{op}$, which is also a Koszul properad.

In the case of $\IBi$, the properad $\Ao$ is the Koszul operad $\A$ of associative algebras
and the properad $\Bo$ is the reversed properad of the Koszul operad $\A$, $\Bo=\A^{op}$,
which is also a Koszul properad.
\end{proo}

\textsc{application :} The cobar construction of the coproperad
$\BLi^{\ac}$ is a resolution of the properad $\BLi$. We can
interpret this homological result in terms of graphs cohomology.
Since the Koszul dual properad of $\BLi$ is a dioperad, the
differential of the cobar construction of $\BLi^{\ac}$ is equal to
the boundary map defined by M. Kontsevich in the context of graph
cohomology (vertex expansion). Therefore, we find the results of
M. Markl and A.A. Voronov in \cite{MV} : the cohomology of the
directed connected commutative graphs is equal to the properad
$\BLi$ and the cohomology of the directed connected ribbon graphs
is equal to the properad $\IBi$.

\subsection{$\Po$-gebra up to homotopy}
\label{Pgebrahomotopy}

\begin{dei}[$\Po$-gebra up to homotopy]

When $\Po$ is a Koszul properad, we call
\emph{$\Po$-gebra up to homotopy} any gebra on
$\bar{\B}^c(\Po^{\ac})$. We denote the category
of $\Po$-gebras up to homotopy by
\emph{$\Po_\infty$-gebras}.
\end{dei}

This notion generalizes to gebras the notion of algebra up to
homotopy (on an operad). For instance, we have the notion of Lie
bialgebras up to homotopy and the notion of infinitesimal Hopf
algebras up to homotopy.

\section{Poincar\'e series}

We generalize here the Poincar\'e series of quadratic associative algebras (\emph{cf.}
J.-L. Loday \cite{Loday1}) and binary quadratic operads (\emph{cf.}
V. Ginzburg and M.M. Kapranov \cite{GK}) to properads. When a properad $\Po$
is Koszul, we show a functional equation between the Poincar\'e series of $\Po$ and $\Po^{\ac}$.\\

\subsection{Poincar\'e series for Koszul properads}

Let $\Po$ be a weight graded connected dg-properad. Denote by $\K$ its Koszul complex
$\K(I,\, \Po,\, \Po)=\Po^{\ac}\boxtimes_c \Po$. This complex is the directed sum of its
sub-complexes indexed by the global weight $ \K = \bigoplus_{m,\, n,\, d \, \ge 0} \K_{(d)}(m,\, n), $
where $\K_{(d)}(m,\, n)$  has the following form :
$$ 0 \to \Po_{(d)}^{\ac}(m,\, n) \to \underbrace{\Po^{\ac}}_{(d-1)}\boxtimes_c
\underbrace{\Po}_{(1)}(m,\, n) \to \cdots \to
\underbrace{\Po^{\ac}}_{(1)}\boxtimes_c
\underbrace{\Po}_{(d-1)}(m,\, n) \to \Po_{(d)}(m,\, n) \to 0.$$

The main theorem of this paper (Theorem~\ref{criteredeKoszul}) asserts that a properad $\Po$ is Koszul
if and only if each chain complexes $\K_{(d)}(m,\, n)$ is acyclic for $d>0$. In this case,
we know from Corollary~\ref{Koszul2quadratic} that the properad $\Po$ is quadratic.
Suppose that $\Po$ is a quadratic properad generated
by an $\Sy$-module $V$ such that the total dimension of
$\bigoplus_{m,\, n \in \mathbb{N}^*} V(m, \,n)$ is finite.
In Proposition~\ref{dimfinie} we have seen that each module
$\F_{(d)}(V)(m,\,n)$ has finite dimension. Therefore, the dimension of each
$\Po_{(d)}(m,\, n)$ and $\Po^{\ac}_{(d)}(m,\, n)$ is finite.
The Euler-Poincar\'e characteristic of $\K_{(d)}(m,\, n)$ is null
$$ \sum_{k=0}^{d}(-1)^k dim (\underbrace{\Po^{\ac}}_{(k)} \boxtimes_c
\underbrace{\Po}_{(d-k)})(m,\, n)=0.$$
The Euler-Poincar\'e characteristic is given by the following formula
\begin{eqnarray*}
&& \sum_{k=0}^{d}(-1)^k dim (\underbrace{\Po^{\ac}}_{(k)}
\boxtimes_c
\underbrace{\Po}_{(d-k)})(m,\, n) = \\
&& \sum_\Xi \sharp \mathcal{S}_c^{\bar{k},\, \bar{j}}
\frac{n!}{\oi !\, \oj !}\, \frac{m!}{\ok ! \, \ol !} \,
dim\Po_{(o_1)}^{\ac} (l_1,\, k_1)\ldots dim\Po_{(o_b)}^{\ac}(l_b,\,
k_b).\\
&& \quad \quad \quad \quad \quad  \quad \quad  \quad
dim\Po_{(q_1)}(j_1,\, i_1) \ldots dim\Po_{(q_a)}(j_a,\,
i_a),
\end{eqnarray*}
where the sum $\Xi$ runs over the $n$-tuples $\oi$, $\oj$, $\ok$,
$\ol$, $\bar{o}$ and $\bar{q}$ such that $|\oi|=n$, $|\oj|=|\ok|$,
$|\ol|=m$, $|\bar{o}|=k$ and $|\bar{q}|=d-k$.

\begin{dei}[Poincar\'e series associated to an $\Sy$-bimodule]
To a weight graded $\Sy$-bimodule $\Po$, we associate the
\emph{Poincar\'e series} defined by
$$f_\Po(y,\, x,\, z):=\sum_{m,n \ge 1 \atop d\ge 0} \frac{\textrm{dim}\Po_{(d)}(m,\, n)}{m!\,
n!}y^mx^nz^d,$$
when the dimension of every $k$-modules $\Po_{(d)}(m,\, n)$ is finite.
\end{dei}

We define the function $\Psi$ by the formula
$$\Psi\big(g(y,\, X,\, z),\, f(Y,\, x,\,
z')\big):=\sum_{\Xi'} \sharp S_c^{\bar{k},\, \bar{j}}
\prod_{\beta=1}^b \frac{1}{k_\beta!}  \frac{\partial^{k_\beta
}g}{{\partial X}^{k_\beta}}(y,\, 0,\, z)
\prod_{\alpha=1}^a\frac{1}{j_\alpha!}\frac{\partial^{j_\alpha}f}{{\partial
Y}^{j_\alpha}}(0,\, x,\, z'),$$
where the sum $\Xi'$ runs over the $n$-tuples $\ok$ and $\oj$ such that
$|\ok|=|\oj|$.

\begin{thm}
\label{seriedepoincare}
Every Koszul properad $\Po$ generated by an $\Sy$-module
of global finite dimension verifies the following equation
$$\Psi\big(f_{\Po^{\ac}}(y,\, X,\, -z),\, f_\Po(Y,\, x,\, z)\big)=xy.$$
\end{thm}

\begin{proo}
We have
\begin{eqnarray*}
&& \Psi\big(f_{\Po^{\ac}}(y,\, X,\, -z),\, f_\Po(Y,\, x,\, z)\big) = \\
&&  \sum_{\Xi'}  \sharp S_c^{\bar{k},\, \bar{j}} \prod_{\beta=1}^b
\frac{1}{k_\beta!}  \frac{\partial^{k_\beta
}f_{\Po^{\ac}}}{{\partial X}^{k_\beta}}(y,\, 0,\, -z)
\prod_{\alpha=1}^a\frac{1}{j_\alpha!}\frac{\partial^{j_\alpha}f_\Po}{{\partial
Y}^{j_\alpha}}(0,\, x,\, z)  \\
& &\sum_{m,\, n \ge 1 \atop d\ge 0} \Big( \underbrace{\sum_{k=0}^d
(-1)^k \Big(\sum_\Xi \sharp S_c^{\bar{k},\, \bar{j}}
\prod_{\beta=1}^b \frac{dim\Po_{q_\beta}^{\ac}(l_\beta,\,
k_\beta)}{l_\beta!\, k_\beta!} \prod_{\alpha=1}^a
\frac{dim\Po_{o_\alpha}(j_\alpha,\, i_\alpha)}{j_\alpha!\,
i_\alpha!} \Big)}_{=0 \quad \textrm{for} \quad d>0} \Big) y^mx^nz^d    \\
& &=xy.
\end{eqnarray*}
\end{proo}

\begin{rem}
Since the Koszul complex $\K(I,\, \Po,\, \Po)$ is acyclic if and only if the Koszul complex
$\K(\Po,\, \Po,\, I)$ is acyclic, we also have the symmetric relation
$$ \Psi\big(f_\Po (y,\, X,\, -z),\, f_{\Po^{\ac}}(Y,\, x,\, z)\big)=xy.$$
\end{rem}

\subsection{Poincar\'e series for associative Koszul algebras}
When the quadratic properad $\Po$ is a quadratic algebra $A$, its
Poincar\'e series is equal to $\sum_{d\ge 0} dim (A_{(d)}) z^d \,
xy.$ If we denote
$$f_A(z):=\sum_{d\ge 0} dim(A_{(d)})z^d,$$
then the functional equation of the previous theorem is equivalent to the well
known relation (\emph{cf.} \cite{Loday1})
$$f_A(x).f_{A^{\ac}}(-x)=1.$$

\subsection{Poincar\'e series for Koszul operads}

When the quadratic pro\-pe\-rad $\Po$ is a quadratic operad, its
Poincar\'e series is equal to
$$\sum_{d\ge 0,\, n\ge 1}\frac{dim\Po_{(d)}(n)}{n !}\, x^n\, y\, z^d.$$
We denote $$f_{\Po}(x,\, z):=\sum_{d\ge 0,\, n\ge 1}\frac{dim\Po_{(d)}(n)}{n !}\, x^n\, z^d.$$

\begin{cor}
\label{EFKoszuloperad} The Poincar\'e series of a Koszul operad
$\Po$ verifies the relation
$$f_{\Po^{\ac}}(f_{\Po}(x,\, z),\, -z)=x.$$
\end{cor}

In the binary case (when $V=V(2)$), we have $\Po_{(n-1)}=\Po(n)$ and
$$f_{\Po}(x,\, z)=\sum_{n\ge 1}\frac{dim\Po(n)}{n !}\, x^n\, z^{n-1}.$$
We denote $$f_\Po(x):=\sum_{n\ge 1} (-1)^n
\frac{dim\Po(n)}{n!} x^n,$$
and if $\Po$ is Koszul, we have the preceding theorem gives the well known relation
(\emph{cf.} \cite{Loday2} Appendix B)
$$f_{\Po^{\ac}}(f_\Po(x))=x.$$

\begin{ex1}
Consider the free operad $\Po=\F(V)$ generated by the $\Sy$-module
$V=k\{\,
\vcenter{\xymatrix@M=0pt@R=6pt@C=6pt{\ar@{-}[dr] &  &\ar@{-}[dl]  \\  &\ar@{-}[d] & \\
& &}} ,\,
\vcenter{\xymatrix@M=0pt@R=6pt@C=6pt{\ar@{-}[dr] & \ar@{-}[d] &\ar@{-}[dl]  \\
&\ar@{-}[d] & \\  & &}}
 ,\, \ldots \}$
composed by one $n$-ary operation for every $n$.
Since it is a free operad, it is a Koszul operad. Whereas the total dimension
of $V$ is infinite, the dimension of every $\F(V)(n)$ is finite. Therefore, we can consider
the Euler-Poincar\'e characteristic of the Koszul complex of $\Po$.
Since $\Po^{\ac}=k\oplus V$, we have
 $$ f_{\Po^{\ac}}(x,\, z)= \sum_{d\ge 0,\, n\ge 1} dim\Po_{(d)}^{\ac}(n)x^n z^d=x+
 \sum_{n\ge 2} x^n z=x+z \frac{x^2}{1-x}.$$
And Corollary~\ref{EFKoszuloperad} gives the following equation
$$(z+1)f_\Po^2(x,\, z)-(1+x)f_\Po(x,\, z)+x=0.$$

Let $P_n(z)$ be the Poincar\'e polynomial of the Stasheff polytope
of dimension $n$, also called the \emph{associahedra} and denoted
$K^n$ or $K_{n+2}$. This polynomial is equal to
$P_n(z):=\sum_{k=0}^n \sharp Cel_{k}^{n}.z^k$, where $Cel_{k}^{n}$
represents the set of the cells of dimension $k$ of the polytope
of dimension $n$. Denote $f_K(x,\, z)=\sum_{n\ge 0}P_n(z)x^n$ the
generating series.

The cells of dimension $k$ of  $K^n$ are indexed by
planar trees with $n+2$ leaves and $n+1-k$ vertices.
This bijection implies that $\sharp Cel_{k}^{n}=dim
\Po_{n+1-k}(n+2)$, which gives
\begin{eqnarray}
f_\Po(x,\, z) &=& x+ \sum_{n\ge 2} \sum_{k=1}^{n-1} dim \Po_k^n \, z^kx^n \nonumber \\
&=& x+\sum_{n\ge2} \Big( \sum_{k=1}^{n-1} \sharp
Cel_{n-2-(k-1)}^{l-2}\, z^k \Big) x^n
\nonumber \\
&=& x+zx^2\sum_{n\ge 0} \Big( \sum_{k=0}^{n}\sharp Cel_{n-k}^{n}\,
z^k \Big) x^n
\nonumber \\
&=& x +zx^2\sum_{n\ge 0} P_n\Big(\frac{1}{y}\Big)\, (xz)^n =
x+zx^2f_K\Big(xz,\, \frac{1}{z}\Big). \nonumber
\end{eqnarray}

Therefore, we get the following relation verified by $f_K$
$$ ((1+z)x^2)f_K^2(x,\, z)+(-1+(2+z)x)f_K(x,\, z)+1=0.$$

Finally, the generating series of the Stasheff polytopes verifies
$$ f_K(x,\, z)=\frac{1+(2+z)x-\sqrt{1-2(2+z)x+z^2x^2}}{2(1+z)x^2}.$$
\end{ex1}

\bigskip

\begin{center}
\textsc{Acknowledgements}
\end{center}
$ $

This paper is part of the author Ph.D. thesis. I would like to
thank my advisor Jean-Louis Loday for his confidence. I wish also
to thank Benoit Fresse for his help and for the proof of
Theorem~\ref{echelonnement}. I am grateful to Birgit Richter and
Martin Markl who have carefully read my thesis and to Wee Liang
Gan who pointed some inaccuracies.



\tableofcontents

\medskip

{\small \textsc{Laboratoire J.A. Dieudonn\'e, Universit\'e de Nice, Parc Valrose, 06108 Nice
Cedex 02, France}\\
E-mail address : \texttt{brunov@math.unice.fr}\\
URL : \texttt{http://math.unice.fr/$\sim$brunov}}

\end{document}